\documentclass[12pt]{amsbook}
\usepackage{amsbsy,amssymb,amscd,amsfonts,latexsym,amstext,delarray,
amsmath,epsfig}
\input xypic

\newtheorem{thm}{Theorem}[section]
\newtheorem{prop}[thm]{Proposition}

\newtheorem{lem}[thm]{Lemma}

\newtheorem{defn}{Definition}[section]

\numberwithin{equation}{chapter}

\def\R{{\mathbb R}}
\def\A{{\mathbb A}}
\def\H{{\mathbb H}}

\def\Z{{\mathbb Z}}
\def\Q{{\mathbb Q}}
\def\C{{\mathbb C}}
\def\K{{\mathbb K}}
\def\P{{\mathbb P}}
\def\N{{\mathbb N}}
\def\F{{\mathbb F}}
\def\PSL{{\rm PSL}}
\def\PGL{{\rm PGL}}
\def\SL{{\rm SL}}
\def\GL{{\rm GL}}
\def\Tr{{\rm Tr}}

\def\AdS{{\rm AdS}}
\def\Aut{{\rm Aut}}
\def\End{{\rm End}}
\def\Hom{{\rm Hom}}
\def\Ker{{\rm Ker}}
\def\Coker{{\rm Coker}}
\def\sD{{\mathcal D}}
\def\sO{{\mathcal O}}
\def\cO{{\mathcal O}}

\def\cA{{\mathcal A}}

\def\cB{{\mathcal B}}
\def\cH{{\mathcal H}}
\def\cE{{\mathcal E}}
\def\cU{{\mathcal U}}
\def\cF{{\mathcal F}}
\def\cL{{\mathcal L}}
\def\cM{{\mathcal M}}
\def\Sp{{\rm Spec}}
\def\T{{\mathcal T}}
\def\Gal{{\rm Gal}}
\def\Ker{{\rm Ker}}
\def\Coker{{\rm Coker}}
\def\sE{{\mathcal E}}
\def\mX{{\mathfrak X}}
\def\fC{{\mathfrak{C}}}
\newcommand{\hL}{\mathbb{L}}
\def\m{{\mathfrak m}}

\def\ma{{\mathfrak a}}
\def\O{{\mathcal O}}

\def\cW{{\mathcal W}}
\def\cP{{\mathcal P}}

\def\cS{{\mathcal S}}
\def\cV{{\mathcal V}}

\newcommand{\ie}{{\it i.e.\/}\ }
\newcommand{\eg}{{\it e.g.\/}\ }
\newcommand{\cf}{{\it cf.\/}\ }

\title{Lectures on Arithmetic Noncommutative Geometry}
\author{Matilde Marcolli}
\date{}

\begin{document}

\maketitle

{\small
\begin{verse}
And indeed there will be time \\
To wonder ``Do I dare?'' and, ``Do I dare?'' \\
Time to turn back and descend the stair. \\
\smallskip
...\\
\smallskip
Do I dare \\
Disturb the Universe? \\
\smallskip
...\\
\smallskip
For I have known them all already, known them all; \\
Have known the evenings, mornings, afternoons, \\
I have measured out my life with coffee spoons. \\
\smallskip
...\\
\smallskip
I should have been a pair of ragged claws \\
Scuttling across the floors of silent seas.\\
...\\
No! I am not Prince Hamlet, nor was meant to be;\\
Am an attendant lord, one that will do\\
To swell a progress, start a scene or two\\
...\\
At times, indeed, almost ridiculous--\\
Almost, at times, the Fool.\\ 
...\\
We have lingered in the chambers of the sea\\
By sea-girls wreathed with seaweed red and brown\\
Till human voices wake us, and we drown.\\
\bigskip
(T.S.~Eliot, {\em ``The Love Song of J. Alfred Prufrock''})
\end{verse}
}

\newpage

\tableofcontents

\chapter{Ouverture}

Noncommutative geometry, as developed by Connes starting in the
early '80s (\cite{ConnesCR}, \cite{Co}, \cite{Co94}), extends the
tools of ordinary geometry to treat spaces that are quotients, for
which the usual ``ring of functions'', defined as functions
invariant with respect to the equivalence relation, is too small
to capture the information on the ``inner structure'' of
points in the quotient space. Typically, for such spaces functions
on the quotients are just constants, while a nontrivial ring of
functions, which remembers the structure of the equivalence
relation, can be defined using a noncommutative
algebra of coordinates, analogous to the non-commuting variables of quantum
mechanics. These ``quantum spaces'' are defined by extending the
Gelfan'd--Naimark correspondence
$$ X \text{ loc.comp. Hausdorff space} \Leftrightarrow C_0(X)
\text{ abelian $C^*$-algebra } $$ by dropping the commutativity
hypothesis in the right hand side. The correspondence then becomes a
definition of what is on the left hand side: a noncommutative space.

Such quotients are abundant in nature. They arise, for instance, from
foliations. Several recent results also show that noncommutative
spaces arise naturally in number theory and arithmetic geometry.
The first instance of such connections between noncommutative geometry
and number theory emerged in the work of Bost and
Connes \cite{BC}, which exhibits a very interesting noncommutative
space with remarkable arithmetic properties related to class field
theory. This reveals a very useful dictionary that relates the
phenomena of spontaneous symmetry breaking in quantum statistical
mechanics to the mathematics of Galois theory.
This space can be viewed as the space of
1-dimensional $\Q$-lattices up to scale, modulo the equivalence relation of
commensurability (\cf \cite{CoMar}).
This space is closely related to the noncommutative space used by
by Connes to obtain a spectral realization
of the zeros of the Riemann zeta function, \cite{Connes-Zeta}.
In fact, this is again the space of
commensurability classes of 1-dimensional $\Q$-lattices, but with the
scale factor also taken into account.

More recently, other results that point to deep connections between
noncommutative geometry and number theory appeared in the work of
Connes and Moscovici \cite{CoMo1} \cite{CoMo2} on the modular Hecke
algebras. This shows that
the Rankin--Cohen brackets, an important algebraic structure on
modular forms \cite{Zagier}, have a natural interpretation in the
language of noncommutative geometry, in terms of the Hopf algebra of
the transverse geometry of codimension one foliations.
The modular Hecke algebras, which naturally combine products and action
of Hecle operators on modular forms, can be viewed as the ``holomorphic
part'' of the algebra of coordinates on
the space of commensurability classes of 2-dimensional $\Q$-lattices
constructed in joint work of Connes and the author \cite{CoMar}.

Cases of occurrences of interesting number theory within
noncommutative geometry can be found in the classification of
noncommutative three spheres by Connes and Dubois--Violette
\cite{CDV1} \cite{CDV2}. Here the
corresponding moduli space has a ramified cover by a noncommutative
nilmanifold, where the noncommutative analog of the Jacobian of this
covering map is expressed naturally in terms of the ninth power of
the Dedekind eta function. Another such case occurs in Connes'
calculation \cite{CoQgr} of the explicit
cyclic cohomology Chern character of a spectral triple on
$SU_q(2)$ defined by Chakraborty and Pal \cite{ChaPal}.

Other instances of noncommutative spaces that arise in the context of
number theory and arithmetic geometry can be found in the
noncommutative compactification of modular curves of \cite{CDS},
\cite{ManMar}. This noncommutative space is again related to the
noncommutative geometry of $\Q$-lattices. In fact, it can be seen as a
stratum in the compactification of the space of commensurability
classes of 2-dimensional $\Q$-lattices (\cf \cite{CoMar}).

Another context in which noncommutative geometry provides a useful
tool for arithmetic geometry is in the description of the totally
degenerate fibers at ``arithmetic infinity'' of arithmetic
varieties over number fields, analyzed in joint work of the
author with Katia Consani (\cite{CM}, \cite{CM1}, \cite{CM2},
\cite{CM3}).

The present text is based on a series of lectures given by the
author at Vanderbilt University in May 2004, as well as on  
previous series of lectures given at the Fields Institute 
in Toronto (2002), at the University of Nottingham (2003), and 
at CIRM in Luminy (2004). 

The main focus of the lectures is the noncommutative geometry of
modular curves (following \cite{ManMar}) and of the archimedean fibers
of arithmetic varieties (following \cite{CM}). A chapter on the
noncommutative space of commensurability 
classes of 2-dimensional $\Q$-lattices is also included (following
\cite{CoMar}). The text reflects very closely the style of the
lectures. In particular, we have tried more to convey the general
picture than the details of the proofs of the specific results. Though
many proofs have not been included in the text, the reader will find
references to the relevant literature, where complete proofs are
provided (in particular \cite{CoMar}, \cite{CM}, and \cite{ManMar}). 

More explicitly, the text is organized as follows:
\begin{itemize}
\item We start by recalling a few preliminary notions of
noncommutative geometry (following \cite{Co94}).
\item The second chapter describes how various arithmetic properties of
modular curves can be seen by their ``noncommutative boundary''. This
part is based on the joint work of Yuri Manin and the author. The main
references are \cite{ManMar}, \cite{Mar-lyap},
\cite{Mar-cosm}.
\item We review briefly the work of Connes and the author \cite{CoMar}
on the noncommutative geometry of commensurability classes of
$\Q$-lattices. The relation of the
noncommutative space of commensurability classes of $\Q$-lattices to
the Hilbert 12th problem of explicit class field theory is based on
ongoing work of Connes, Ramachandran and the author \cite{CMR}, on the
original work of Bost and Connes \cite{BC} and on Manin's real
multiplication program \cite{Man3} \cite{Man5}.
\item The noncommutative geometry of the fibers at ``arithmetic
infinity'' of varieties over number fields is based on joint work of
Consani and the author, for which the references are \cite{CM},
\cite{CM1}, \cite{CM2}, \cite{CM3}, \cite{CMrev}. This chapter
also contains a detailed account of Manin's formula for the Green
function of Arakelov geometry for arithmetic surfaces, based on
\cite{Man-hyp}, and a proposed physical interpretation of this
formula, as in \cite{ManMar2}.
\end{itemize}

\section{The NCG dictionary}

There is a dictionary (\cf \cite{Co94}) relating concepts of
ordinary geometry to the corresponding counterparts in
noncommutative geometry. The entries can be arranged according to
the finer structures considered on the underlying space, roughly
according to the following table.

\bigskip

\begin{tabular}{|c|c|}\hline
measure theory & von Neumann algebras \\ \hline topology &
$C^*$--algebras \\ \hline smooth structures & smooth subalgebras
\\ \hline Riemannian geometry &  spectral triples \\ \hline
\end{tabular}

\bigskip

It is important to notice that, usually, the notions of
noncommutative geometry are ``richer'' than the corresponding
entries of the dictionary on the commutative side. For instance,
as Connes discovered, noncommutative measure spaces (von Neumann
algebras) come endowed with a natural time evolution which is
trivial in the commutative case. Similarly, at the level of
topology one often sees phenomena that are closer to rigid
analytic geometry. This is the case, for instance, with the
noncommutative tori $T_\theta$, which already at the $C^*$-algebra
level exhibit moduli that behave much like moduli of
one-dimensional complex tori (elliptic curves) in the commutative
case.

In the context we are going to discuss this richer structure of
noncommutative spaces is crucial, as it permits us to use tools
like $C^*$-algebras (topology) to study the properties more rigid
spaces like algebraic or arithmetic varieties.

\section{Noncommutative spaces}

The way to assign the algebra of coordinates to a quotient space
$X=Y/\sim$ can be explained in a short slogan as follows:

\begin{itemize}
\item Functions on $Y$ with $f(a)=f(b)$ for $a\sim b$. {\bf Poor!}
\item Functions $f_{ab}$ on the graph of the
equivalence relation. {\bf Good!}
\end{itemize}

The second description leads to a noncommutative algebra, as the
product, determined by the groupoid law of the equivalence
relation, has the form of a convolution product (like the product
of matrices).

For sufficiently {\em nice} quotients, even though the two notions are
not the same, they are related by Morita equivalence, which is the
suitable notion of ``isomorphism'' between noncommutative spaces.
For more general quotients, however, the two notions truly differ and the
second one is the only one that allows one to continue to make
sense of geometry on the quotient space.

A very simple example illustrating the above situation is the
following (\cf \cite{Co2}). Consider the topological space
$Y=[0,1] \times \{ 0,1\}$ with the equivalence relation $(x,0)\sim
(x,1)$ for $x\in (0,1)$. By the first method one only obtains
constant functions $\C$, while by the second method one obtains
$$
\{ f\in C([0,1])\otimes M_2(\C) : \, f(0) \text{ and } f(1) \text{
diagonal } \}
$$
which is an interesting nontrivial algebra.

The idea of preserving the information on the structure of the
equivalence relation in the description of quotient spaces has
analogs in Grothendieck's theory of stacks in algebraic geometry.

\subsection{Morita equivalence}

In noncommutative geometry, isomorphisms of $C^*$-algebras are too
restrictive to provide a good notion of isomorphisms of
noncommutative spaces. The correct notion is provided by Morita
equivalence of $C^*$-algebras.

We have equivalent ${\rm C}^*$-algebras ${\mathcal A}_1
\sim{\mathcal A}_2$ if there exists a bimodule ${\mathcal M}$,
which is a right Hilbert ${\mathcal A}_1$ module with an
${\mathcal A}_1$-valued inner product $\langle \cdot,\cdot
\rangle_{ {\mathcal A}_1}$, and a left Hilbert ${\mathcal
A}_2$-module with an ${\mathcal A}_2$-valued inner product
$\langle \cdot,\cdot \rangle_{ {\mathcal A}_2}$, such that we
have:
\begin{itemize}
\item We obtain all ${\mathcal A}_i$ as the closure of the span of
$$ \{ \langle \xi_1, \xi_2 \rangle_{{\mathcal A}_i}: \xi_1, \xi_2\in
{\mathcal M} \}. $$
\item $\forall \xi_1, \xi_2, \xi_3 \in {\mathcal M}$ we have
$$ \langle \xi_1, \xi_2 \rangle_{{\mathcal A}_1} \xi_3 = \xi_1 \langle
\xi_2, \xi_3 \rangle_{{\mathcal A}_2}. $$
\item ${\mathcal A}_1$ and ${\mathcal A}_2$ act on
${\mathcal M}$ by bounded operators,
$$ \langle a_2 \xi, a_2 \xi \rangle_{{\mathcal A}_1} \leq \| a_2 \|^2
\langle \xi, \xi \rangle_{{\mathcal A}_1} \ \ \ \ \langle a_1 \xi,
a_1 \xi \rangle_{{\mathcal A}_2} \leq \| a_1 \|^2 \langle \xi, \xi
\rangle_{{\mathcal A}_2} $$ for all $a_1 \in {\mathcal A}_1$, $a_2
\in {\mathcal A}_2$, $\xi \in {\mathcal M}$.
\end{itemize}

\medskip

This notion of equivalence roughly means that one can transfer
modules back and forth between the two algebras.

\subsection{The tools of noncommutative geometry}

Once one identifies in a specific problem a space that, by its
nature of quotient of the type described above, is best described
as a noncommutative space, there is a large set of well developed
techniques that one can use to compute invariants and extract
essential information from the geometry. The following is a list
of some such techniques, some of which will make their appearance
in the cases treated in these notes.

\begin{itemize}
\item Topological invariants: K-theory
\item Hochschild and cyclic cohomology
\item Homotopy quotients, assembly map (Baum-Connes)
\item Metric structure: Dirac operator, spectral triples
\item Characteristic classes, zeta functions
\end{itemize}

We will recall the necessary notions when needed. We now begin by taking
a closer look at the analog in the noncommutative world of
Riemannian geometry, which is provided by Connes' notion of
spectral triples.

\section{Spectral triples}

Spectral triples are a crucial notion in noncommutative geometry.
They provide a powerful and flexible generalization of the
classical structure of a Riemannian manifold. The two notions
agree on a commutative space. In the usual context of Riemannian
geometry, the definition of the infinitesimal element $ds$ on a
smooth spin manifold can be expressed in terms of the inverse of
the classical Dirac operator $D$. This is the key remark that
motivates the theory of spectral triples. In particular, the
geodesic distance between two points on the manifold is defined in
terms of $D^{-1}$ (\cf \cite{Co94} \S VI). The spectral triple
$(A,H,D)$ that describes a classical Riemannian spin manifold is
given by the algebra $A$ of complex valued smooth functions on the
manifold, the Hilbert space $H$ of square integrable spinor
sections, and the classical Dirac operator $D$. These data
determine completely and uniquely the Riemannian geometry on the
manifold. It turns out that, when expressed in this form, the
notion of spectral triple extends to more general non-commutative
spaces, where the data $(A,H,D)$ consist of a ${\rm C}^*$-algebra
$A$ (or more generally of some smooth subalgebra of a ${\rm
C}^*$-algebra) with a representation in the algebra of bounded
operators on a separable Hilbert space $H$, and an operator $D$ on
$H$ that verifies the main properties of a Dirac operator.

We recall the basic setting of Connes theory of spectral triples.
For a more complete treatment see \cite{Co94}, \cite{Connes}, 
\cite{ConnesMosc}.

\smallskip

\begin{defn}\label{specDef}
A spectral triple $({\mathcal A}, {\mathcal H}, D)$ consists of a
${\rm C}^*$-algebra ${\mathcal A}$ with a representation
$$ \rho : {\mathcal A} \to {\mathcal B}({\mathcal H}) $$
in the algebra of bounded operators on a separable Hilbert space
${\mathcal H}$, and an operator  $D$ (called the Dirac operator)
on ${\mathcal H}$, which satisfies the following properties:
\begin{enumerate}
\item $D$ is self--adjoint.
\item For all $\lambda\notin \R$, the resolvent $(D-\lambda)^{-1}$ is
a compact operator on ${\mathcal H}$.
\item The commutator $[D,\rho(a)]$ is a bounded
operator on ${\mathcal H}$, for all $a\in {\mathcal A}_0\subset
\cA$, a dense involutive subalgebra of $\cA$.
\end{enumerate}
\end{defn}

The property {\em 2.} of Definition \ref{specDef} can be regarded
as a generalization of the ellipticity property of the standard
Dirac operator on a compact manifold. In the case of ordinary
manifolds, we can consider as subalgebra $\cA_0$ the algebra of
smooth functions, as a subalgebra of the commutative ${\rm
C}^*$-algebra of continuous functions. In fact, in the classical
case of Riemannian manifolds, property {\em 3.} is equivalent the
Lipschitz condition, hence it is satisfied by a larger class than
that of smooth functions.

\smallskip

Thus, the basic geometric structure encoded by the theory of spectral
triples is Riemannian geometry, but in more refined cases, such as
K\"ahler geometry, the additional structure can be easily encoded as
additional symmetries. We will see, for instance, a case (\cf
\cite{CM} \cite{CMrev}) where
the algebra involves the action of the Lefschetz operator of a compact
K\"ahler manifold, hence it encodes the information (at the
cohomological level) on the K\"ahler form.

\smallskip

Since we are mostly interested in the relations between
noncommutative geometry and arithmetic geometry and number theory,
an especially interesting feature of spectral triples is that they
have an associated family of zeta functions and a theory of
volumes and integration, which is related to special values of
these zeta functions. (The following treatment is based on
\cite{Co94}, \cite{Connes}.)

\subsection{Volume form.}

A spectral triple $({\mathcal A}, {\mathcal H}, D)$ is said to be
of dimension $n$, or $n$--{\it summable} if the operator
$|D|^{-n}$ is an infinitesimal of order one, which means that the
eigenvalues $\lambda_k(|D|^{-n})$ satisfy the estimate
$\lambda_k(|D|^{-n})=O(k^{-1})$.

For a positive compact operator $T$ such that
$$ \sum_{j=0}^{k-1} \lambda_j(T) = O(\log k), $$
the Dixmier trace $\Tr_\omega(T)$ is the coefficient of this
logarithmic divergence, namely
\begin{equation}\label{DixTr}
\Tr_\omega(T) = \lim_\omega \frac{1}{\log k} \sum_{j=1}^k
\lambda_j(T). \end{equation} Here the notation $\lim_\omega$ takes
into account the fact that the sequence
$$ S(k,T):=\frac{1}{\log k} \sum_{j=1}^k \lambda_j(T) $$
is bounded though possibly non-convergent. For this reason, the
usual notion of limit is replaced by a choice of a linear form
$\lim_\omega$ on the set of bounded sequences satisfying suitable
conditions that extend analogous properties of the limit. When the
sequence $S(k,T)$ converges \eqref{DixTr} is just the ordinary
limit $\Tr_\omega(T)= \lim_{k\to \infty} S(k,T)$. So defined, the
Dixmier trace \eqref{DixTr} extends to any compact operator that
is an infinitesimal of order one, since any such operator is the
difference of two positive ones. The operators for which the
Dixmier trace does not depend on the choice of the linear form
$\lim_\omega$ are called {\em measurable operators}.

On a non-commutative space the operator $|D|^{-n}$ generalizes the
notion of a volume form. The volume is defined as
\begin{equation} \label{VolDix} V = \Tr_\omega
(|D|^{-n}). \end{equation} More generally, consider the algebra
$\tilde {\mathcal A}$ generated by ${\mathcal A}$ and
$[D,{\mathcal A}]$. Then, for $a\in \tilde {\mathcal A}$,
integration with respect to the volume form $|D|^{-n}$ is defined
as
\begin{equation} \label{intDix} \int a := \frac{1}{V} \Tr_\omega (a
|D|^{-n}). \end{equation}

The usual notion of integration on a Riemannian spin manifold $M$
can be recovered in this context (\cf \cite{Co94}, \cite{Landi})
through the formula ($n$ even):
$$ \int_M f dv = \left( 2^{n-[n/2]-1} \pi^{n/2} n \Gamma(n/2) \right)
\Tr_\omega (f |D|^{-n}). $$ Here $D$ is the classical Dirac
operator on $M$ associated to the metric that determines the
volume form $dv$, and $f$ in the right hand side is regarded as
the multiplication operator acting on the Hilbert space of square
integrable spinors on $M$.

\subsection{Zeta functions.}

An important function associated to the Dirac operator $D$ of a
spectral triple $({\mathcal A}, {\mathcal H}, D)$ is its zeta
function
\begin{equation}\label{zetaD}
 \zeta_D (z) := \Tr(|D|^{-z}) = \sum_{\lambda} \Tr(\Pi(\lambda,|D|))
\lambda^{-z}, \end{equation} where $\Pi(\lambda,|D|)$ denotes the
orthogonal projection on the eigenspace $E(\lambda,|D|)$.

\smallskip

An important result in the theory of spectral triples (\cite{Co94}
\S IV Proposition 4) relates the volume \eqref{VolDix} with the
residue of the zeta function \eqref{zetaD} at $s=1$ through the
formula
\begin{equation} \label{VolRes}
V = \lim_{s\to 1+} (s-1)\zeta_D(s) = Res_{s=1} \Tr (|D|^{-s}).
\end{equation}

\smallskip

There is a family of zeta functions associated to a spectral
triple $({\mathcal A}, {\mathcal H}, D)$, to which \eqref{zetaD}
belongs. For an operator $a\in \tilde {\mathcal A}$, we can define
the zeta functions
\begin{equation} \label{aDzeta}
\zeta_{a,D}(z) := \Tr (a |D|^{-z}) = \sum_{\lambda} \Tr(a\,
\Pi(\lambda,|D|)) \lambda^{-z}
\end{equation}
and
\begin{equation} \label{aDzetas}
\zeta_{a,D}(s,z) := \sum_{\lambda} \Tr(a\, \Pi(\lambda,|D|))
(s-\lambda)^{-z}.
\end{equation}
These zeta functions are related to the heat kernel $e^{-t|D|}$ by
Mellin transform
\begin{equation} \label{Mellin}
\zeta_{a,D}(z) = \frac{1}{\Gamma(z)} \int_0^\infty t^{z-1} \Tr(
a\, e^{-t|D|} )\, dt
\end{equation}
where
\begin{equation} \label{theta}
\Tr( a\, e^{-t|D|} ) = \sum_\lambda \Tr(a\, \Pi(\lambda,|D|))
e^{-t\lambda} =: \theta_{a,D}(t).
\end{equation}
Similarly,
\begin{equation} \label{Mellins}
\zeta_{a,D}(s,z) = \frac{1}{\Gamma(z)} \int_0^\infty
\theta_{a,D,s}(t) \, t^{z-1} \, dt
\end{equation}
with
\begin{equation} \label{thetas}
\theta_{a,D,s}(t) := \sum_\lambda \Tr(a\, \Pi(\lambda,|D|))
e^{(s-\lambda)t}.
\end{equation}
Under suitable hypothesis on the asymptotic expansion of
\eqref{thetas} (\cf Theorem 2.7-2.8 of \cite{Man4} \S 2), the
functions \eqref{aDzeta} and \eqref{aDzetas} admit a unique
analytic continuation (\cf \cite{ConnesMosc}) and there is an
associated regularized determinant in the sense of Ray--Singer
(\cf \cite{RS}):
\begin{equation} \label{zetadet}
 {\det_\infty}_{a,D}(s) := \exp \left( -\frac{d}{dz} \zeta_{a,D}(s,z)
|_{z=0} \right).
\end{equation}

\smallskip

The family of zeta functions \eqref{aDzeta} also provides a
refined notion of dimension for a spectral triple $({\mathcal A},
{\mathcal H}, D)$, called the {\it dimension spectrum}. This is a
subset $\Sigma=\Sigma({\mathcal A}, {\mathcal H}, D)$ in $\C$ with
the property that all the zeta functions \eqref{aDzeta}, as $a$
varies in $\tilde {\mathcal A}$, extend holomorphically to $\C
\setminus \Sigma$.

Examples of spectral triples with dimension spectrum not contained
in the real line can be constructed out of Cantor sets.

\subsection{Index map.} The data of a spectral triple determine
an index map. In fact, the self adjoint operator $D$ has a polar
decomposition, $D= F |D|$, where $|D|=\sqrt{D^2}$ is a positive
operator and $F$ is a sign operator, \ie $F^2=I$. Following
\cite{Co94} and \cite{Connes} one defines a cyclic cocycle
\begin{equation} \label{ch} \tau (a^0,a^1,\ldots,a^n) = \Tr (a^0
[F,a^1]\cdots [F,a^n]), \end{equation} where $n$ is the dimension
of the spectral triple. For $n$ even $\Tr$ should be replaced by a
super trace, as usual in index theory. This cocycle pairs with the
$K$-groups of the algebra ${\mathcal A}$ (with $K_0$ in the even
case and with $K_1$ in the odd case), and defines a class $\tau$
in the cyclic cohomology $HC^n({\mathcal A})$ of ${\mathcal A}$.
The class $\tau$ is called the Chern character of the spectral
triple $({\mathcal A}, {\mathcal H}, D)$.

\smallskip

In the case when the spectral triple $({\mathcal A}, {\mathcal H},
D)$ has {\em discrete} dimension spectrum $\Sigma$, there is a
local formula for the cyclic cohomology Chern character
\eqref{ch}, analogous to the local formula for the index in the
commutative case. This is obtained by producing Hochschild
representatives (\cf \cite{Connes}, \cite{ConnesMosc})
\begin{equation}\label{chHH}
\varphi(a^0,a^1,\ldots,a^n) = \Tr_\omega (a^0 [D,a^1]\cdots
[D,a^n]\, |D|^{-n}).
\end{equation}

\subsection{Infinite dimensional geometries}

The main difficulty in constructing specific examples of spectral
triples is to produce an operator $D$ that at the same time has
bounded commutators with the elements of ${\mathcal A}$ and
produces a non-trivial index map.

It sometimes happens that a noncommutative space does not admit 
a finitely summable spectral triple, namely one such that the
operator $|D|^{-p}$ is of trace class for some $p>0$. 
Obstructions to the existence of such spectral triples are
analyzed in \cite{Co-fredh}. It is then useful to consider a weaker
notion, namely that of $\theta$-summable spectral triples. These
have the property that the Dirac operator satisfies
\begin{equation}\label{thetasum}
\Tr(e^{-tD^2})<\infty \ \ \ \  \forall t>0. 
\end{equation}
Such spectral triples should be thought of as ``infinite dimensional
noncommutative geometries''.

We'll see examples of spectral triples that are $\theta$-summable but 
not finitely summable, because of the growth rate of the
multiplicities of eigenvalues.

\subsection{Spectral triples and Morita equivalences}\label{MeqS3}

If $({\mathcal A}_1,{\mathcal H},D)$ is a spectral triple, and we
have a Morita equivalence ${\mathcal A}_1 \sim {\mathcal A}_2$
implemented by a  bimodule ${\mathcal M}$ which is a finite
projective right Hilbert module over ${\mathcal A}_1$, then we can
transfer the spectral triple from ${\mathcal A}_1$ to ${\mathcal
A}_2$.

First consider the ${\mathcal A}_1$-bimodule $\Omega_D^1$
generated by
$$ \{ a_1 [D,b_1]: a_1,b_1\in {\mathcal A}_1 \}. $$

We define a connection
$$ \nabla : {\mathcal M} \to {\mathcal
M}\otimes_{{\mathcal A}_1} \Omega_D^1 $$ by requiring that
$$\nabla (\xi a_1) = (\nabla \xi) a_1 + \xi \otimes [D,a_1], $$
$\forall \xi\in {\mathcal M}$, $\forall a_1\in {\mathcal A}_1$ and
$\forall \xi_1,\xi_2\in {\mathcal M}$. We also require that
$$ \langle \xi_1, \nabla \xi_2 \rangle_{{\mathcal A}_1} - \langle
\nabla \xi_1, \xi_2\rangle_{{\mathcal A}_1} = [D, \langle \xi_1,
\xi_2 \rangle_{{\mathcal A}_1} ]. $$

This induces a spectral triple $({\mathcal A}_2,\tilde {\mathcal
H},\tilde D)$ obtained as follows.

The Hilbert space is given by $\tilde {\mathcal H} ={\mathcal
M}\otimes_{{\mathcal A}_1} {\mathcal H}$. The action takes the
form
$$
a_2\,\, (\xi \otimes_{{\mathcal A}_1} x) := (a_2\xi)
\otimes_{{\mathcal A}_1} x.
$$

The Dirac operator is given by
$$ \tilde D(\xi \otimes x)= \xi\otimes D(x) + ( \nabla \xi) x. $$

Notice that we need a Hermitian connection $\nabla$, because
commutators $[D,a]$ for $a\in {\mathcal A}_1$ are non-trivial,
hence $1\otimes D$ would not be well defined on $\tilde {\mathcal
H}$.

\subsection{K-theory of $C^*$-algebras}

The $K$-groups are important invariants of $C^*$-algebras that
capture information on the topology of non-commutative spaces,
much like cohomology (or more appropriately topological
$K$-theory) captures information on the topology of ordinary
spaces. For a $C^*$-algebra $\cA$, we have:

\subsubsection{$K_0(\cA)$.} This group is obtained by considering
idempotents ($p^2=p$) in matrix algebras over $\cA$. Notice that,
for a $C^*$-algebra, it is sufficient to consider projections
($P^2=P$, $P=P^*$). In fact, we can always replace an idempotent
$p$ by a projection $P=pp^* (1-(p-p^*)^2)^{-1}$, preserving the
von Neumann equivalence. Then we consider $P\simeq Q$ if and only
if $P=X^*X$ $Q=XX^*$, for $X=$ partial isometry ($X=XX^*X$), and
we impose the stable equivalence: $P\sim Q$, for $P\in \cM_n(\cA)$
and $Q\in \cM_m(\cA)$, if and only if there exists $R$ a
projection with $P\oplus R \simeq Q \oplus R$. We define
$K_0(\cA)^+$ to be the monoid of projections modulo these
equivalences and we let $K_0(\cA)$ be its Grothendieck group.

\subsubsection{$K_1(\cA)$.} Consider the group $\GL_n(\cA)$ of
invertible elements in the matrix algebra $\cM_n(\cA)$, and
$\GL_n^0(\cA)$ the identity component of $\GL_n(\cA)$. The
morphism
$$ \GL_n(\cA) \to \GL_{n+1}(\cA) \ \ \ \ a \mapsto
\left(\begin{array}{cc} a & 0 \\ 0 & 1 \end{array}\right) $$
induces
$$ \GL_n(\cA)/\GL_n^0(\cA) \to
\GL_{n+1}(\cA)/\GL_{n+1}^0(\cA).
$$
The group $K_1(\cA)$ is defined as the direct limit of these
morphisms. Notice that $K_1(\cA)$ is abelian even if the
$\GL_n(\cA)/\GL_n^0(\cA)$ are not.

\medskip

In general, $K$-theory is not easy to compute. This lead to 
two fundamental developments in noncommutative geometry. The
first is cyclic cohomology (\cf \cite{Co}, \cite{Co-transv},
\cite{Co94}), introduced by Connes in 1981, which provides cycles that
pair with $K$-theory, and a Chern character. 
The second development is
the geometrically defined $K$-theory of Baum--Connes
\cite{BaumConnes} and the assembly map 
\begin{equation}\label{assem}
\mu: K^*(X,G)\to K_*(C_0(X)\rtimes G)
\end{equation}
from the geometric to the analytic $K$-theory. Much work went,
in recent years, into exploring the range of validity of
the Baum--Connes conjecture, according to which the
assembly map gives an isomorphism between these two different
notions of $K$-theory.

\chapter{Noncommutative modular curves}

The results of this section are mostly based on the joint work of Yuri
Manin and the author \cite{ManMar},
with additional material from \cite{Mar-lyap} and \cite{Mar-cosm}.
We include necessary preliminary notions on modular curves and on
noncommutative tori, respectively based on the papers
\cite{Man-sym} and \cite{ConnesCR}, \cite{Rief1}.

We first recall some aspects of the classical theory of modular
curves. We will then show that these notions can be entirely
recovered and enriched through the analysis of a noncommutative
space associated to a compactification of modular curves, which
includes noncommutative tori as possible degenerations of elliptic
curves.

\section{Modular curves}

Let $G$ be a finite index subgroup of the modular group $\Gamma
=\PSL(2,\Z)$, and let $X_G$ denote the quotient
\begin{equation}\label{XG}
X_G := G \backslash \H^2
\end{equation}
where $\H^2$ is the 2-dimensional real hyperbolic plane, namely
the upper half plane $\{ z\in \C : \Im  z > 0 \}$ with the metric
$ds^2=|dz|^2/(\Im z)^2$. Equivalently, we identify $\H^2$ with the
Poincar\'e disk $\{ z : | z |<1 \}$ with the metric $ds^2 = 4
|dz|^2/(1-|z|^2)^2$.

We denote the quotient map by $\phi: \H^2 \to X_G$.

Let $\P$ denote the coset space $\P := \Gamma/G$. We can write the
quotient $X_G$ equivalently as
$$ X_G = \Gamma \backslash (\H^2 \times \P). $$

\begin{figure}
\begin{center}
\epsfig{file=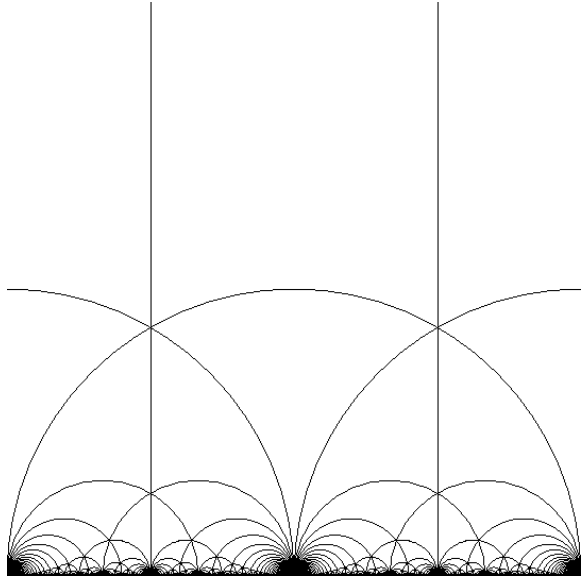}
\end{center}
\caption{Fundamental domains for $\PSL(2,\Z)$
\label{Fig-modcurve}}
\end{figure}

A tessellation of the hyperbolic plane $\H^2$ by fundamental domains
for the action of $\PSL(2,\Z)$ is illustrated in Figure
\ref{Fig-modcurve}. The action of $\PSL(2,\Z)$ by fractional linear
transformations $z \mapsto \frac{az+b}{cz+d}$ is usually written in
terms of the generators $S: z\mapsto -1/z$ and $T: z\mapsto z+1$
(inversion and translation). Equivalently, $\PSL(2,\Z)$ can be
identified with the free product $\PSL(2,\Z)\cong \Z/2 * \Z/3$, with
generators $\sigma$, $\tau$ of order two and three, respectively given
by $\sigma=S$ and $\tau=ST$.

An example of finite index subgroups is given by the congruence
subgroups $\Gamma_0(N)\subset \Gamma$, of matrices
$$ \left(\begin{array}{cc} a&b \\ c&d \end{array}\right) $$
with $c \equiv 0 \mod N$. A fundamental domain for $\Gamma_0(N)$ is
given by $F \cup_{k=0}^{N-1} ST^k (F)$, where $F$ is a fundamental
domain for $\Gamma$ as in Figure \ref{Fig-modcurve}\footnote{Figures
\ref{Fig-modcurve} and \ref{Fig-geods} are taken from
Curt McMullen's Gallery}.

The quotient space $X_G$ has the structure of a non-compact
Riemann surface. This has a natural algebro-geometric
compactification, which consists of adding the cusp points (points
at infinity). The cusp points are identified with the quotient
\begin{equation}\label{cusps}
G \backslash \P^1(\Q)\simeq \Gamma \backslash ( \P^1(\Q) \times \P
).
\end{equation}
Thus, we write the compactification as
\begin{equation}\label{alg-compact}
\overline{X_G}:= G \backslash ( \H^2 \cup \P^1(\Q)) \simeq \Gamma
\backslash \left( ( \H^2 \cup \P^1(\Q)) \times \P \right).
\end{equation}

The modular curve $X_\Gamma$, for $\Gamma =\SL(2,\Z)$ is the
moduli space of elliptic curves, with the point $\tau\in \H^2$
parameterizing the lattice $\Lambda=\Z \oplus \tau\Z$ in $\C$ and
the corresponding elliptic curve uniformized by
$$ E_\tau = \C / \Lambda. $$
The unique cusp point corresponds to the degeneration of the
elliptic curve to the cylinder $\C^*$, when $\tau\to \infty$ in
the upper half plane.

The other modular curves, obtained as quotients $X_G$ by a congruence
subgroup, can also be interpreted as moduli spaces: they are moduli
spaces of elliptic curves with level structure. Namely, for elliptic
curves $E=\C/\Lambda$, this is an additional information on the
torsion points $\frac{1}{N}\Lambda/\Lambda\subset \Q\Lambda/\Lambda$
of some level $N$.

For instance, in the case of the principal congruence subgroups
$\Gamma(N)$ of matrices
$$ M=\left(\begin{array}{cc} a&b \\ c&d \end{array}\right) \in \Gamma $$
such that $M \equiv Id \mod N$, points in the modular curve
$\Gamma(N)\backslash\H^2$ classify elliptic curves $E_tau=
\C/\Lambda$, with $\Lambda=\Z+\Z\tau$, together with a basis $\{
1/N, \tau/N \}$ for the torsion subgroup
$\frac{1}{N}\Lambda/\Lambda$. The projection $X_{\Gamma(N)}\to
X_\Gamma$ forgets the extra structure.

In the case of the groups $\Gamma_0(N)$, points in the quotient
$\Gamma_0(N)\backslash\H^2$ classify elliptic curves together with
a cyclic subgroup of $E_\tau$ of order $N$. This extra information
is equivalent to an isogeny $\phi: E_\tau \to E_{\tau'}$ where the
cyclic group is $Ker(\phi)$. Recall that an isogeny is a morphism
$\phi: E_\tau \to E_{\tau'}$ such that $\phi(0)=0$. These are
implemented by the action of $\GL_2^+(\Q)$ on $\H^2$, namely
$E_\tau$ and $E_{\tau'}$ are isogenous if and only if $\tau$ and
$\tau'$ in $\H^2$ are in the same orbit of $\GL_2^+(\Q)$.

\subsection{Modular symbols}

\begin{figure}
\begin{center}
\epsfig{file=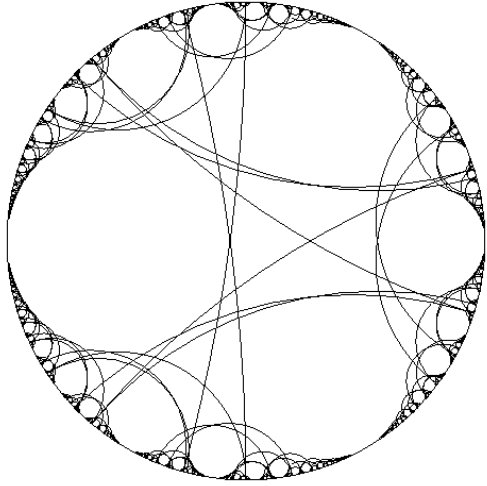}
\end{center}
\caption{Geodesics between cusps define modular symbols
\label{Fig-geods}}
\end{figure}

Given two points $\alpha,\beta \in \H^2 \cup \P^1(\Q)$, a real
homology class $\{ \alpha,\beta \}_G \in H_1(X_G,\R)$ is defined
as
$$ \{ \alpha,\beta \}_G : \omega \mapsto \int_{\alpha}^\beta
\phi^*(\omega), $$ where $\omega$ are holomorphic differentials on
$X_G$, and the integration of the pullback to $\H^2$ is along the
geodesic arc connecting two cusps $\alpha$ and  $\beta$ (\cf Figure
\ref{Fig-geods}).

The modular symbols $\{ \alpha,\beta \}_G$ satisfy the additivity
and invariance properties
$$ \{\alpha ,\beta\}_G + \{\beta ,\gamma\}_G = \{\alpha ,\gamma\}_G,
$$
and
$$
\{ g\alpha ,g\beta \}_G=\{\alpha ,\beta\}_G,
$$
for all $g\in G$.

Because of additivity, it is sufficient to consider modular
symbols of the form $\{ 0, \alpha \}$ with $\alpha\in \Q$,
$$ \{ 0,\alpha \}_G = -\sum_{k=1}^n \{ g_k(\alpha)\cdot 0,
g_k(\alpha)\cdot i\infty \}_G, $$ where $\alpha$ has continued
fraction expansion $\alpha = [a_0,\ldots,a_n]$, and
$$ g_k(\alpha) =\left(\begin{array}{cc} p_{k-1}(\alpha) & p_k(\alpha)
\\ q_{k-1}(\alpha) & q_k (\alpha) \end{array}\right), $$
with $p_k/q_k$ the successive approximations, and $p_n/q_n
=\alpha$.

\medskip

In the classical theory of the modular symbols of \cite{Man-sym},
\cite{Merel} (cf.~ also \cite{Ha}, \cite{Maz}) cohomology classes
obtained from cusp forms are evaluated against relative homology
classes given by modular symbols, namely, given a cusp form $\Phi$
on $\H$, obtained as pullback $\Phi=\varphi^*(\omega)/dz$ under
the quotient map $\varphi: \H \to X_G$, we denote by
$\Delta_\omega (s)$ the intersection numbers $\Delta_\omega (s) =
\int_{g_s (0)}^{g_s (i\infty)} \Phi(z) dz$, with $g_s G= s\in \P$.
These intersection numbers can be interpreted in terms of special
values of $L$--functions associated to the automorphic form which
determines the cohomology class.

We rephrase this in cohomological terms following \cite{Merel}. We
denote by $\tilde I$ and $\tilde R$ the elliptic points, namely
the orbits $\tilde I= \Gamma \cdot i$ and $\tilde R = \Gamma \cdot
\rho$. We denote by $I$ and $R$ the image in $X_G$ of the elliptic
points
$$ I= G\backslash \tilde I \ \ \ \  R= G\backslash \tilde R, $$
with $\rho=e^{\pi i/3}$. We use the notation
\begin{equation}\label{HBA}
H^B_A := H_1(\overline{X_G} \setminus A, B; \Z).
\end{equation}
These groups are related by the pairing
\begin{equation}\label{pairingAB}
 H^B_A \times H^A_B \to \Z.
\end{equation}

The modular symbols $\{ g(0), g(i\infty) \}$, for $gG\in \P$,
define classes in $H^{\text{cusps}}$. For $\sigma$ and $\tau$ the
generators of $\PSL(2,\Z)$ with $\sigma^2=1$ and $\tau^3=1$, we
set
\begin{equation}\label{PIPR}
\P_I = \langle \sigma \rangle \backslash \P \ \ \ \text{and} \ \ \
\P_R =\langle \tau \rangle \backslash \P.
\end{equation}
There is an isomorphism $\Z^{|\P |} \cong H^R_{\text{cusps} \cup
I}$. Given the exact sequences
$$ 0 \to H_{\text{cusps}} \stackrel{\iota'}{\to} H^R_{\text{cusps}}
\stackrel{\pi_R}{\to} \Z^{| \P_R |} \to \Z \to 0 $$ and
$$ 0\to \Z^{| \P_I |} \to H^R_{\text{cusps} \cup I}
\stackrel{\pi_I}{\to} H^R_{\text{cusps}} \to 0, $$ the image
$\pi_I(\tilde x)\in H^R_{\text{cusps}}$ of an element $\tilde x =
\sum_{s\in \P} \lambda_s s$ in $\Z^{|\P|}\cong
H^R_{\text{cusps}\cup I}$ represents an element $x \in
H_{\text{cusps}}$ iff the image $\pi_R(\pi_I(\tilde x)) =0$ in
$\Z^{| \P_R |}$. As proved in \cite{Merel}, for $s=gG \in \P$, the
intersection pairing $\bullet : H^{\text{cusps}} \times
H_{\text{cusps}} \to \Z$ gives
$$ \{ g(0), g(i\infty) \} \bullet x = \lambda_s - \lambda_{\sigma
s}. $$ Thus, we write the intersection number as a function
$\Delta_x : \P \to \R$ by
\begin{equation} \label{deltax} \Delta_x (s)= \lambda_s -
\lambda_{\sigma s}, \end{equation} where $x$ is given as above.

\subsection{The modular complex}

For $x$ and $y$ in $\H^2$ we denote by $\langle x, y \rangle$ the
oriented geodesic arc connecting them. Moreover, in the
decomposition $\PSL(2,\Z) = \Z/2 *\Z/3$ as a free product, we
denote by $\sigma$ the generator of $\Z/2$ and by $\tau$ the
generator of $\Z/3$.

\begin{defn} The {\em modular complex} is the cell complex defined as
follows.
\begin{itemize}
\item 0-cells: the cusps $G\backslash \P^1(\Q)$, and the elliptic
points $I$ and $R$.
\item 1-cells: the oriented geodesic arcs
$G \backslash (\Gamma \cdot \langle i\infty, i \rangle)$ and
$G\backslash (\Gamma \cdot \langle i, \rho \rangle)$, where by
$\Gamma \cdot$ we mean the orbit under the action of $\Gamma$.
\item 2-cells: $ G\backslash \{ \Gamma \cdot E \}$, where $E$ is the
polygon with vertices $\{ i,\rho,1+i, i\infty \}$ and sides the
corresponding geodesic arcs.
\item Boundary operator: $ \partial : C_2 \to C_1 $ is given by
$$ g E \mapsto g \langle i, \rho \rangle + g \langle \rho, 1+i
\rangle + g \langle 1+i, i\infty \rangle + g \langle i\infty, i
\rangle, $$ for $g\in \Gamma$, and the boundary $ \partial : C_1
\to C_0 $ is given by
$$ g \langle i\infty, i \rangle \mapsto g(i) - g(i\infty)  $$
$$ g \langle i, \rho \rangle \mapsto g(\rho)- g(i).  $$
\end{itemize}
\end{defn}

\begin{figure}
\begin{center}
\epsfig{file=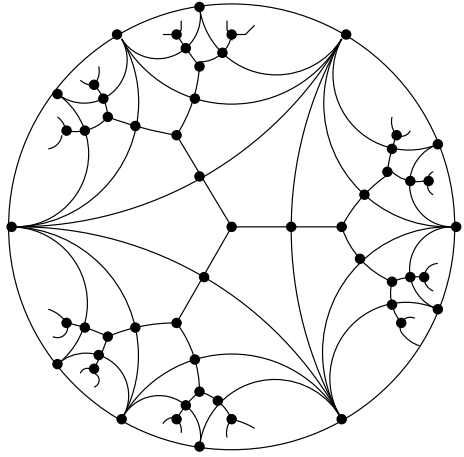}
\end{center}
\caption{The cell decomposition of the modular complex
\label{Fig-modcomplex}}
\end{figure}

This gives a cell decomposition of $\H^2$ adapted to the action of
$\PSL(2,\Z)$ and congruence subgroups, \cf Figure
\ref{Fig-modcomplex}.

We have the following result \cite{Man-sym}.

\begin{prop}
The modular complex computes the first homology of
$\overline{X_G}$:
\begin{equation}\label{mod-com-hom}
H_1(X_G)\cong \frac{Ker(\partial: C_1 \to C_0)}{Im(\partial:
C_2\to C_1)}.
\end{equation}
\end{prop}

We can derive versions of the modular complex that compute
relative homology.

Notice that we have $\Z[\text{cusps}]=C_0/ \Z[R\cup I]$, hence the
quotient complex
$$ 0 \to C_2 \stackrel{\partial}{\to} C_1 \stackrel{\tilde\partial}{\to}
\Z[\text{cusps}] \to 0,  $$ with $\tilde\partial$ the induced
boundary operator, computes the relative homology
$H_1(\overline{X_G},R\cup I)$. The cycles are given by $\Z[\P]$,
as combinations of elements $g \langle i, \rho \rangle$, $g$
ranging over representatives of $\P$, and by the elements $\oplus
a_g \langle g(i\infty), g(i) \rangle$ satisfying $\sum a_g
g(i\infty)=0$. In fact, these can be represented as relative
cycles in $(\overline{X_G},R\cup I)$. Similarly, the subcomplex
$$ 0 \to \Z[\P] \stackrel{\partial}{\to} \Z[R\cup I] \to 0, $$
with $\Z[\P]$ generated by the  elements $g \langle i, \rho
\rangle$, computes the homology $H_1(\overline{X_G} -
\text{cusps})$. The homology
$$ H_1(\overline{X_G} - \text{cusps}, R\cup I) \cong \Z[\P]  $$
is generated by the relative cycles $g \langle i, \rho \rangle$,
$g$ ranging over representatives of $\P$.

With the notation \eqref{HBA}, we consider the groups
$H^{\text{cusps}}_{R\cup I}$, $H_{\text{cusps}}^{R\cup I}$,
$H^{\text{cusps}}$, and $H_{\text{cusps}}$.

We have a long exact sequence of relative homology
\begin{equation}\label{rel-homol}
 0\to H_{\text{cusps}} \to H_{\text{cusps}}^{R\cup I}
\stackrel{(\tilde\beta_R,\tilde\beta_I)}{\to} H_0(R)\oplus H_0(I)
\to \Z \to 0,
\end{equation}
with $H_{\text{cusps}}$ and $H_{\text{cusps}}^{R\cup I}$ as above,
and with
$$ H_0(I)\cong \Z[ \P_I ], \,\,\,\,\,\,
\P_I= \langle \sigma \rangle\backslash \P = G\backslash \tilde I
$$
$$ H_0(R)\cong \Z[ \P_R ], \,\,\,\,\,\, \P_R=
\langle \tau \rangle\backslash \P = G\backslash \tilde R,  $$ that
is,
\begin{equation}\label{rel-hom2}
0\to H_{\text{cusps}} \to \Z[\P] \to \Z[\P_R] \oplus \Z[\P_I] \to
\Z \to 0.
\end{equation}

In the case of $H_{\text{cusps}}^{R\cup I}$ and $H_{R\cup
I}^{\text{cusps}}$ the pairing \eqref{pairingAB} gives the
identification of $\Z[\P]$ and $\Z^{|\P|}$, obtained by
identifying the elements of $\P$ with the corresponding delta
functions. Thus, we can rewrite the sequence \eqref{rel-hom2} as
\begin{equation}\label{rel-hom3}
0 \to H^{\text{cusps}} \to \Z^{ | \P | }
\stackrel{(\beta_R,\beta_I)}{\to} \Z^{ | \P_I | } \oplus  \Z^{ |
\P_R |} \to \Z \to 0,
\end{equation}
where $H_{R\cup I}^{\text{cusps}} \cong \Z^{ | \P | }$.

In order to understand more explicitly the map $(\beta_R,\beta_I)$
we give the following equivalent {\em algebraic} formulation of
the modular complex (\cf \cite{Man-sym} \cite{Merel}).

The homology group $H_{\text{cusps}}^{R\cup I}= \Z[\P]$ is
generated by the images in $X_G$ of the geodesic segments
$g\gamma_0:=g\langle i, \rho \rangle$, with $g$ ranging over a
chosenset of representatives of the coset space $\P$.

We can identify (\cf \cite{Merel}) the dual basis $\delta_s$ of
$H^{\text{cusps}}_{R\cup I}= \Z^{|\P|}$ with the images in $X_G$
of the paths $g\eta_0$, where for a chosen point $z_0$ with $0<
Re(z_0)< 1/2$ and $|z_0|>1$ the path $\eta_0$ is given by the
geodesic arcs connecting $\infty$ to $z_0$, $z_0$ to $\tau z_0$,
and $\tau z_0$ to $0$. These satisfy
$$ [g \gamma_0] \bullet [g \eta_0] =1 $$
$$ [g \gamma_0] \bullet [h \eta_0] =0, $$
for $g G \neq h G$, under the intersection pairing
\eqref{pairingAB}

Then, in the exact sequence \eqref{rel-hom3}, the identification
of $H^{\text{cusps}}$ with $Ker(\beta_R,\beta_I)$ is obtained by
the identification $\{ g(0), g(i\infty) \}_G \mapsto g\eta_0$, so
that the relations imposed on the generators $\delta_s$ by the
vanishing under $\beta_I$ correspond to the relations $\delta_s
\oplus \delta_{\sigma s}$ (or $\delta_s$ if $s=\sigma s$) and the
vanishing under $\beta_R$ gives another set of relations $\delta_s
\oplus \delta_{\tau s} \oplus \delta_{\tau^2 s}$ (or $\delta_s$ if
$s=\tau s$).

\section{The noncommutative boundary of modular curves}

The main idea that bridges between the algebro--geometric theory
of modular curves and noncommutative geometry consists of
replacing $\P^1(\Q)$ in the classical compactification, which
gives rise to a finite set of cusps, with $\P^1(\R)$. This
substitution cannot be done naively, since the quotient
$G\backslash \P^1(\R)$ is ill behaved topologically, as $G$ does
not act discretely on $\P^1(\R)$.

When we regard the quotient $\Gamma\backslash\P^1(\R)$, or more
generally $\Gamma\backslash(\P^1(\R)\times \P)$, itself as a
noncommutative space, we obtain a geometric object that is rich
enough to recover many important aspects of the classical theory
of modular curves. In particular, it makes sense to study in terms
of the geometry of such spaces the {\em limiting behavior} for
certain arithmetic invariants defined on modular curves when $\tau
\to \theta \in \R\setminus \Q$.

\section{Modular interpretation: noncommutative elliptic curves}

The boundary $\Gamma\backslash\P^1(\R)$ of the modular curve
$\Gamma\backslash\H^2$, viewed itself as a noncommutative space,
continues to have a modular interpretation, as observed originally
by Connes--Douglas--Schwarz (\cite{CDS}). In fact, we can think of
the quotients of $S^1$ by the action of rotations by an irrational
angle (that is, the noncommutative tori) as particular
degenerations of the classical elliptic curves, which are
``invisible'' to ordinary algebraic geometry. The quotient space
$\Gamma\backslash\P^1(\R)$ classifies these noncommutative tori up
to Morita equivalence (\cite{ConnesCR}, \cite{Rief1}) completes
the moduli space $\Gamma\backslash \H^2$ of the classical elliptic
curves. Thus, from a conceptual point of view it is reasonable to
think of $\Gamma\backslash\P^1(\R)$ as the boundary of
$\Gamma\backslash \H^2$, when we allow points in this classical
moduli space (that is, elliptic curves) to have non-classical
degenerations to noncommutative tori.

Noncommutative tori are, in a sense, a prototype example of
noncommutative spaces, in as one can see there displayed the full
range of techniques of noncommutative geometry (\cf
\cite{ConnesCR}, \cite{Co}). As ${\rm C}^*$--algebras,
noncommutative tori are {\em irrational rotation algebras}. We
recall some basic properties of noncommutative tori, which justify
the claim that these algebras behave like a noncommutative version
of elliptic curves. We follow mostly \cite{ConnesCR} \cite{Co94}
and \cite{Rief1} for this material.

\subsection{Irrational rotation and Kronecker
foliation}\label{Kron}

\begin{defn}\label{irr-rot}
The irrational rotation algebra $\cA_\theta$, for a given
$\theta\in \R$, is the universal ${\rm C}^*$-algebra ${\rm
C}^*(U,V)$, generated by two unitary operators $U$ and $V$,
subject to the commutation relation
\begin{equation}\label{commutUV}
UV=e^{2\pi i\theta} VU.
\end{equation}
\end{defn}

The algebra $A_\theta$ can be realized as a subalgebra of bounded
operators on the Hilbert space $\cH=L^2(S^1)$, with the circle
$S^1\cong \R/\Z$. For a given $\theta\in \R$, we consider two
operators, that act on a complete orthonormal basis $e_n$ of $\cH$
as
\begin{equation}\label{UV}
U e_n =  e_{n+1}, \ \ \ \ \ V e_n = e^{2\pi in\theta} e_n.
\end{equation}
It is easy to check that these operators satisfy the commutation
relation \eqref{commutUV}, since, for any $f\in \cH$ we have $VU
f\,(t)= Uf (t-\theta)=e^{2\pi i(t-\theta)} f(t-\theta)$, while $UV
f\, (t)= e^{2\pi i t} f(t-\theta)$.

The irrational rotation algebra be described in more geometric
terms by the foliation on the usual commutative torus by lines
with irrational slope. On $T^2=\R^2/\Z^2$ one considers the
foliation $dx =\theta dy$, for $x,y \in \R/\Z$. The space of
leaves is described as $X=\R/ (\Z + \theta \Z) \simeq S^1 / \theta
\Z$. This quotient is ill behaved as a classical topological
space, hence it cannot be well described by ordinary geometry.

A transversal to the foliation is given for instance by the choice
$T=\{ y=0 \}$, $T\cong S^1\cong\R/\Z$. Then the non-commutative
torus is obtained  (\cf \cite{ConnesCR} \cite{Co94}) as
\begin{equation}\label{NCtori-transv}
\cA_\theta = \{ \left( f_{ab} \right) \, \, \, \,
a,b\in T \, \, \text{ in the same leaf } \}
\end{equation}
where $(f_{ab})$ is a power series $b = \sum_{n \in \Z} b_n V^n$
and each $b_n$ is an element of the algebra $C(S^1)$. The
multiplication is given by
$$
V h V^{-1} = h \circ R_{\theta}^{-1},
$$
with
$$ R_{\theta} x = x + \theta  \mod 1. $$
The algebra $C(S^1)$ is generated by $U(t)=e^{2\pi i t}$, hence we
recover the generating system $(U,V)$ with the relation
$$
UV = e^{2\pi i \theta} VU.
$$

What we have obtained through this description is an
identification of the irrational rotation algebra of Definition
\ref{irr-rot} with the crossed product $C^*$-algebra
\begin{equation}
\label{Rthetacross}
 \cA_\theta =
C(S^1)\rtimes_{R_\theta} \Z
\end{equation}
representing the quotient $S^1/\theta\Z$ as a non-commutative
space.

\subsection{Degenerations of elliptic curves}\label{deg-ell}

An elliptic curve $E_\tau$ over $\C$ can be described as the
quotient $E_\tau = \C /(\Z + \tau \Z)$ of the complex plane by a
2-dimensional lattice $\Lambda=\Z + \tau \Z$, where we can take
$\Im(\tau)>0$. It is also possible to describe the elliptic curve
$E_q$, for $q\in \C^*$, $q=\exp(2\pi i \tau)$, $|q|<1$, in terms
of its Jacobi uniformization, namely as the quotient of $\C^*$ by
the action of the group generated by a single hyperbolic element
in $PSL(2,\C)$,
\begin{equation}\label{Jac-unif}
 E_q = \C^* / q^{\Z}.
\end{equation}
The fundamental domain for the action of $q^{\Z}$ is an annulus
$$ \{ z\in \C: \,\, |q|< z \leq 1 \} $$
of radii $1$ and $|q|$, and the identification of the two boundary
circles is obtained via the combination of scaling and rotation
given by multiplication by $q$.

\begin{figure}
\begin{center}
\epsfig{file=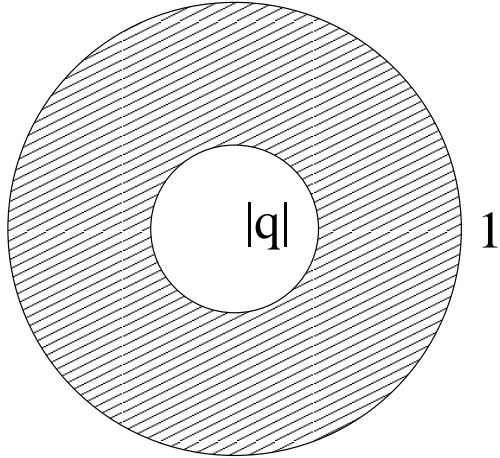}
\end{center}
\caption{The fundamental domain for the Jacobi uniformization of
the elliptic curve \label{Fig-tateU}}
\end{figure}

Now let us consider a degeneration where $q \to \exp(2\pi i\theta)
\in S^1$, with $\theta\in \R\setminus \Q$. We can say
heuristically that in this degeneration the elliptic curve becomes
a non-commutative torus
$$ E_q  \Longrightarrow  \cA_\theta, $$
in the sense that, as we let $q \to \exp(2\pi i\theta)$, the
annulus of the fundamental domain shrinks to a circle $S^1$ and we
are left with a quotient of $S^1$ by the infinite cyclic group
generated by the irrational rotation $\exp(2\pi i \theta)$. Since
this quotient is ill behaved as a classical quotient, such
degenerations do not admit a description within the context of
classical geometry. However, when we replace the quotient by the
corresponding crossed product algebra $C(S^1)\rtimes_\theta \Z$ we
find the irrational rotation algebra of definition \ref{irr-rot}.
Thus, we can consider such algebras as non-commutative
(degenerate) elliptic curves.

More precisely, when one considers degenerations of elliptic
curves $E_\tau = \C^*/q^\Z$ for $q=e^{2\pi i \tau}$, what one obtains
in the limit is the suspension of a noncommutative torus.
In fact, as the
parameter $q$ degenerates to a point on the unit circle, $q\to e^{2\pi
i \theta}$, the ``nice'' quotient $E_\tau = \C^*/q^\Z$ degenerates to
the ``bad'' quotient $E_\theta =\C^*/e^{2\pi i \tau\Z}$, whose
noncommutative algebra of coordinates is Morita equivalent $\cA_\theta
\otimes C_0(\R)$, with $\rho\in \R$ the radial coordinate $\C^* \ni
z=e^\rho e^{2\pi i s}$.

Because of the Thom isomorphism \cite{Co-Thom}, the K-theory of
the noncommutative space $E_\theta =C_0(\R^2)\rtimes (\Z\theta + \Z)$
satisfies
\begin{equation}\label{Kshift}
K_0(E_\theta)=K_1(\cA_\theta) \ \ \ \text{ and } \ \ \
K_1(E_\theta)=K_0(\cA_\theta),
\end{equation}
which is again compatible with the identification of the $E_\theta$
(rather than $\cA_\theta$) as degenerations of elliptic curves.
In fact, for instance, the Hodge filtration on the $H^1$ of an
elliptic curve and the equivalence between the elliptic curve and its
Jacobian, have analogs for the noncommutative torus $\cA_\theta$ in
terms of the filtration on $HC^0$ induced by the inclusion of
$K_0$ (\cf \cite{Co} p.~132--139, \cite{Co2} \S XIII), while by the Thom
isomorphism, these would again appear on the $HC^1$ in the case of
the ``noncommutative elliptic curve'' $E_\theta$.

The point of view of degenerations is sometimes a useful
guideline. For instance, one can study the limiting behavior of
arithmetic invariants defined on the parameter space of elliptic
curves (on modular curves), in the limit when $\tau \to \theta \in
\R\setminus \Q$. An instance of this type of result is the theory
of limiting modular symbols of \cite{ManMar}.

\subsection{Morita equivalent NC tori}

To extend the modular interpretation of the quotient
$\Gamma\backslash \H^2$ as moduli of elliptic curves to the
noncommutative boundary $\Gamma\backslash \P^1(\R)$, one needs to
check that points in the same orbit of the action of the modular
group $\PSL(2,\Z)$ by fractional linear transformations on
$\P^1(\R)$ define equivalent noncommutative tori, where
equivalence here is to be understood in the Morita sense.

Connes showed in \cite{ConnesCR} (\cf also \cite{Rief1}) that the
noncommutative tori $\cA_\theta$ and $\cA_{-1/\theta}$ are Morita
equivalent. Geometrically in terms of the Kronecker foliation and
the description \eqref{NCtori-transv} of the corresponding
algebras, the Morita equivalence $\cA_\theta\simeq
\cA_{-1/\theta}$ corresponds to changing the choice of the
transversal from $T=\{ y=0 \}$ to $T'=\{ x=0 \}$.

In fact, all Morita equivalences arise in this way, by changing
the choice of the transversal of the foliation, so that
$\cA_\theta$ and $\cA_{\theta'}$ are Morita equivalent if and only
if $\theta \sim \theta'$, under the action of $\PSL(2,\Z)$.

Connes constructed in \cite{ConnesCR} explicit bimodules realizing the
Morita equivalences between 
non-commutative tori $\cA_\theta$ and $\cA_{\theta'}$ with
$$ \theta'=\frac{a\theta + b}{c\theta + d}=g\theta, $$
for
$$ g=\left(\begin{array}{cc} a&b\\c&d \end{array}\right) \in
\Gamma,
$$
by taking
$\cM_{\theta,\theta'}$ to be the Schwartz space ${\mathcal S}(\R
\times \Z/c)$, with the right action of $\cA_\theta$
$$ Uf\, (x,u)=f\left( x-\frac{c\theta+d}{c}, u-1\right) $$
$$ Vf\, (x,u)=\exp(2\pi i(x-u d/c)) f(x,u) $$
and the left action of $\cA_{\theta'}$
$$ U'f\, (x,u)=f\left( x-\frac{1}{c}, u-a \right) $$
$$ V'f\, (x,u)=\exp\left( 2\pi i \left(\frac{x}{c\theta+d}
-\frac{u}{c}\right)\right) f(x,u). $$

\subsection{Other properties of NC elliptic curves}

There are other ways in which the irrational rotation algebra
behaves much like an elliptic curve, most notably the relation
between the elliptic curve and its Jacobian (\cf \cite{Co} and
\cite{Co2}) and some aspects of the theory of theta functions,
which we recall briefly.

The commutative torus $T^2=S^1\times S^1$ is connected, hence the
algebra ${\rm C}(T^2)$ does not contain interesting projections.
On the contrary, the noncommutative tori $\cA_\theta$ contain a
large family of nontrivial projections. Rieffel in \cite{Rief1}
showed that, for a given $\theta$ irrational and for all
$\alpha\in (\Z\oplus \Z\theta) \cap [0,1]$, there exists a
projection $P_\alpha$ in $\cA_\theta$, with
$\Tr(P_\alpha)=\alpha$. A different construction of projections in
$\cA_\theta$, given by Boca \cite{Boca}, has arithmetic relevance,
in as these projections correspond to the theta functions for
noncommutative tori defined by Manin in \cite{Man-theta}.

A method of constructing projections in $C^*$-algebras is based on
the following two steps (\cf \cite{Rief1} and \cite{Man5}):

\begin{enumerate}
\item Suppose given a bimodule ${_\cA\cM_\cB}$. If an element $\xi \in
{_\cA\cM_\cB}$ admits an invertible $*$-invariant square root
$\langle \xi, \xi \rangle_{\cB}^{1/2}$, then the element
$\mu:=\xi\langle \xi,\xi\rangle_{\cB}^{-1/2}$ satisfies $ \mu
\langle \mu ,\mu\rangle_{\cB} =\mu $.

\item Let $\mu\in {_\cA\cM_\cB}$ be a non-trivial element such that
$\mu \langle \mu ,\mu\rangle_{\cB} =\mu$. Then the element
$P:=_\cA\langle \mu, \mu \rangle$ is a projection.
\end{enumerate}

In Boca's construction, one obtains elements $\xi$ from Gaussian
elements in some Heisenberg modules, in such a way that the
corresponding $\langle \xi,\xi\rangle_{\cB}$ is a {\em quantum
theta function} in the sense of Manin \cite{Man-theta}. An
introduction to the relation between the Heisenberg groups and the
theory of theta functions is given in the third volume of
Mumford's Tata lectures on theta, \cite{Mum-Tata}.

\section{Limiting modular symbols}

We consider the action of the group $\Gamma=\PGL(2,\Z)$ on the
upper and lower half planes $\H^{\pm}$ and the modular curves
defined by the quotient $X_G=G\backslash \H^{\pm}$, for $G$ a
finite index subgroup of $\Gamma$. Then the noncommutative
compactification of the modular curves is obtained by extending
the action of $\Gamma$ on $\H^{\pm}$ to the action on the full
$$ \P^1(\C)=\H^\pm \cup \P^1(\R), $$
so that we have
\begin{equation}\label{P1CG}
 \overline{X}_G= G\backslash \P^1(\C) = \Gamma \backslash
(\P^1(\C)\times \P).
\end{equation}
Due to the fact that $\Gamma$ does not act discretely on
$\P^1(\R)$, the quotient \eqref{P1CG} makes sense as a
noncommutative space
\begin{equation}\label{ncP1CG}
C(\P^1(\C)\times \P) \rtimes \Gamma.
\end{equation}

Here $\P^1(\R)\subset \P^1(\C)$ is the {\em limit set} of the
group $\Gamma$, namely the set of accumulation points of orbits of
elements of $\Gamma$ on $\P^1(\C)$. We will see another instance
of noncommutative geometry arising from the action of a group of
M\"obius transformations of $\P^1(\C)$ on its limit set, in the
context of the geometry at the archimedean primes of arithmetic
varieties.

\subsection{Generalized Gauss shift and dynamics}

We have described the boundary of modular curves by the crossed
product $C^*$-algebra
\begin{equation}\label{cross-boundary-mc}
C(\P^1(\R)\times \P)\rtimes \Gamma.
\end{equation}
We can also describe the quotient space
$\Gamma\backslash(\P^1(\R)\times \P)$ in the following equivalent
way. If $\Gamma=\PGL(2,\Z)$, then $\Gamma$-orbits in $\P^1(\R)$
are the same as equivalence classes of points of $[0,1]$ under the
equivalence relation
$$ x \sim_T y \Leftrightarrow \exists n,m : T^n x = T^m y $$
where $T x = 1/x - [1/x]$ is the classical Gauss shift of the
continued fraction expansion. Namely, the equivalence relation is
that of having the same tail of the continued fraction expansion
(shift-tail equivalence).

A simple generalization of this classical result yields the
following.

\begin{lem}\label{lem-Tshift}
$\Gamma$ orbits in $\P^1(\R)\times \P$ are the same as equivalence
classes of points in $[0,1]\times \P$ under the equivalence
relation
\begin{equation}\label{equiv-shift}
 (x,s) \sim_T (y,t) \Leftrightarrow \exists n,m : T^n (x,s) = T^m
(y,t) \end{equation}
where $T: [0,1]\times \P \to [0,1] \times \P$
is the shift
\begin{equation}\label{gen-shift}
 T (x,s) = \left( \frac{1}{x} - \left[
\frac{1}{x} \right], \left(\begin{array}{cc} -[1/x] & 1 \\
1 & 0 \end{array}\right)\cdot s \right)
\end{equation}
generalizing the classical shift of the continued fraction
expansion.
\end{lem}

As a noncommutative space, the quotient by the equivalence
relation \eqref{equiv-shift} is described by the $C^*$-algebra of
the groupoid of the equivalence relation
$$
{\mathcal G}([0,1]\times \P, T)=\{ ((x,s), m-n, (y,t)):
T^m(x,s)=T^n(y,t) \}
$$
with objects ${\mathcal G}^0=\{ ((x,s), 0, (x,s)) \}$.

In fact, for any $T$-invariant subset $E\subset [0,1]\times \P$,
we can consider the equivalence relation \eqref{equiv-shift}. The
corresponding groupoid $C^*$-algebra $C^*({\mathcal G}(E,T))$
encodes the dynamical properties of the map $T$ on $E$.

Geometrically, the equivalence relation \eqref{equiv-shift} on
$[0,1]\times \P$ is related to the action of the geodesic flow on
the horocycle foliation on the modular curves.

\subsection{Arithmetic of modular curves and noncommutative boundary}

The result of Lemma \ref{lem-Tshift} shows that the properties of
the dynamical system $T$ or \eqref{equiv-shift} can be used to
describe the geometry of the noncommutative boundary of modular
curves. There are various types of results that can be obtained by
this method (\cite{ManMar} \cite{Mar-lyap}), which we will discuss
in the rest of this chapter.

\begin{enumerate}
\item Using the properties of this dynamical system it is possible to
recover and enrich the theory of modular symbols on $X_G$, by
extending the notion of modular symbols from geodesics connecting
cusps to images of geodesics in $\H^2$ connecting irrational
points on the boundary $\P^1(\R)$. In fact, the irrational points
of $\P^1(\R)$ define {\em limiting modular symbols}. In the case
of quadratic irrationalities, these can be expressed in terms of
the classical modular symbols and recover the generators of the
homology of the classical compactification by cusps
$\overline{X_G}$. In the remaining cases, the limiting modular
symbols vanishes almost everywhere.
\item It is possible to reinterpret Dirichlet
series related to modular forms of weight 2 in terms of integrals
on $[0,1]$ of certain intersection numbers obtained from homology
classes defined in terms of the dynamical system. In fact, even
when the limiting modular symbol vanishes, it is possible to
associate a non-trivial cohomology class in $X_G$ to irrational
points on the boundary, in such a way that an average of the
corresponding intersection numbers give Mellin transforms of
modular forms of weight 2 on $X_G$.
\item The Selberg zeta function of the modular curve can be
expressed as a Fredholm determinant of the Perron-Frobenius
operator associated to the dynamical system on the ``boundary''.
\item Using the first formulation of the boundary as
the noncommutative space \eqref{cross-boundary-mc} we can obtain a
canonical identification of the modular complex with a sequence of
K-groups of the $C^*$-algebra. The resulting exact sequence for
K-groups can be interpreted, using the description
\eqref{equiv-shift} of the quotient space, in terms of the
Baum--Connes assembly map and the Thom isomorphism.
\end{enumerate}

\medskip

All this shows that the noncommutative space
$C(\P^1(\R)\times\P)\rtimes \Gamma$, which we have so far
considered as a boundary stratum of $C(\H^2\times \P)\rtimes
\Gamma$, in fact contains a good part of the arithmetic
information on the classical modular curve itself. The fact that
information on the ``bulk space'' is stored in its boundary at
infinity can be seen as an instance of the physical principle of
holography (bulk/boundary correspondence) in string theory (\cf
\cite{ManMar2}). We will discuss the holography principle more in
details in relation to the geometry of the archimedean fibers of
arithmetic varieties.

\bigskip

\subsection{Limiting modular symbols}

Let $\gamma_\beta$ be an infinite geodesic in the hyperbolic plane
$\H$ with one end at $i\infty$ and the other end at $\beta \in \R
\smallsetminus \Q$. Let $x \in \gamma_\beta$ be a fixed base
point, $\tau$ be the geodesic arc length, and $y(\tau)$ be the
point along $\gamma_\beta$ at a distance $\tau$ from $x$, towards
the end $\beta$. Let $\{ x,y(\tau) \}_G$ denote the homology class
in $X_G$ determined by the image of the geodesic arc $\langle x,
y(\tau) \rangle$ in $\H$.

\begin{defn}\label{lim-MS}
The {\it limiting modular symbol} is defined as
\begin{equation} \label{limitmod}
\{\{* ,\beta\}\}_G :=\lim \frac{1}{\tau}\,\{ x,y(\tau) \}_G\in H_1
(X_G,\R),
\end{equation}
whenever such limit exists.
\end{defn}

The limit \eqref{limitmod} is independent of the choice of the
initial point $x$ as well as of the choice of the geodesic in $\H$
ending at $\beta$, as discussed in \cite{ManMar} (\cf Figure
\ref{Fig-modsymb}). We use the notation $\{\{* ,\beta\}\}_G$ as
introduced in \cite{ManMar}, where $*$ in the first argument
indicates the independence on the choice of $x$, and the double
brackets indicate the fact that the homology class is computed as
a limiting cycle.

\bigskip

\bigskip

\subsubsection{Dynamics of continued fractions}

As above, we consider on $[0,1]\times \P$ the dynamical system
$$ T: [0,1] \times \P \to [0,1] \times \P $$
\begin{equation} \label{shift} T(\beta,t) = \left( \frac{1}{\beta} - \left[
\frac{1}{\beta} \right], \left(\begin{array}{cc} -[1/\beta] & 1 \\
1 & 0 \end{array}\right)\cdot t \right). \end{equation}

This generalizes the classical shift map of the continued fraction
$$ T: [0,1]\to [0,1] \ \ \ \ \ \ T(x)=\frac{1}{x} - \left[
\frac{1}{x}\right]. $$

Recall the following basic notation regarding continued fraction
expansion. Let $k_1,\dots ,k_n$ be independent variables and, for
$n\ge 1$, let
$$
[k_1,\dots ,k_n]:= \frac{1}{k_1+\frac{1}{k_2+\dots
\frac{1}{k_n}}}= \frac{P_n(k_1,\dots ,k_n)}{Q_n(k_1,\dots ,k_n)}.
$$
The $P_n,Q_n$ are polynomials with integral coefficients, which
can be calculated inductively from the relations
$$
Q_{n+1}(k_1,\dots ,k_n,k_{n+1})=k_{n+1}Q_n(k_1,\dots ,k_n)
+Q_{n-1}(k_1,\dots ,k_{n-1}),
$$
$$
P_n(k_1,\dots ,k_n)=Q_{n-1}(k_2,\dots ,k_n),
$$
with $Q_{-1}=0, Q_{0}=1$. Thus, we obtain
$$
[k_1,\dots ,k_{n-1},k_n+x_n]
$$
$$
=\frac{P_{n-1}(k_1,\dots ,k_{n-1})\,x_n+P_{n}(k_1,\dots
,k_{n})}{Q_{n-1}(k_1,\dots ,k_{n-1})\,x_n+Q_{n}(k_1,\dots
,k_{n})}= \left( \begin{array}{cc}
P_{n-1} & P_{n}\\
Q_{n-1} & Q_{n}
\end{array} \right)\,(x_n),
$$
with the standard matrix notation for fractional linear
transformations,
$$
z\mapsto \frac{az+b}{cz+d}= \left( \begin{array}{cc}
a & b\\
c &  d
\end{array} \right)\,(z).
$$

\smallskip

If $\alpha \in (0,1)$ is an irrational number, there is a unique
sequence of integers $k_n(\alpha )\ge 1$ such that $\alpha$ is the
limit of $[\,k_1(\alpha ),\dots ,k_n(\alpha )\,]$ as $n\to\infty$.
Moreover, there is a unique sequence $x_n(\alpha )\in (0,1)$ such
that
$$
\alpha = [\,k_1(\alpha ),\dots ,k_{n-1}(\alpha ), k_n(\alpha
)+x_n(\alpha )\,]
$$
for each $n\ge 1$. We obtain
$$
\alpha = \left( \begin{array}{cc}
0 & 1\\
1 & k_1(\alpha )
\end{array} \right)\ \dots \
\left( \begin{array}{cc}
0 & 1\\
1 & k_n(\alpha )
\end{array} \right)\,(x_n(\alpha ) ).
$$
We set
$$
p_n(\alpha ):=P_n(k_1(\alpha ),\dots ,k_n(\alpha )),\ q_n(\alpha
):=Q_n(k_1(\alpha ),\dots ,k_n(\alpha ))
$$
so that $p_n(\alpha )/q_n(\alpha )$ is the sequence of convergents
to $\alpha$. We also set
$$
g_n(\alpha ):=\left( \begin{array}{cc}
p_{n-1}(\alpha ) & p_n(\alpha )\\
q_{n-1}(\alpha ) & q_n(\alpha ) \end{array} \right) \in \GL(2,\Z).
$$

Written in terms of the continued fraction expansion, the shift
$T$ is given by
$$ T: [k_0, k_1, k_2, \ldots] \mapsto [k_1,
k_2, k_3, \ldots]. $$

\medskip

The properties of the shift \eqref{shift} can be used to extend
the notion of modular symbols to geodesics with irrational ends
(\cite{ManMar}). Such geodesics correspond to infinite geodesics
on the modular curve $X_G$, which exhibit a variety of interesting
possible behaviors, from closed geodesics to geodesics that
approximate some limiting cycle, to geodesics that wind around
different homology class exhibiting a typically chaotic behavior.

\bigskip

\bigskip

\subsubsection{Lyapunov spectrum}

A measure of how chaotic a dynamical system is, or better of how
fast nearby orbits tend to diverge, is given by the Lyapunov
exponent.

\begin{defn}
the {\em Lyapunov exponent} of  $T: [0,1] \to [0,1]$ is defined as
\begin{equation} \label{LyapExp} \lambda(\beta) := \lim_{n\to \infty}
\frac{1}{n} \log | (T^n)^\prime (\beta) | = \lim_{n\to \infty}
\frac{1}{n} \log \prod_{k=0}^{n-1} | T^\prime (T^k \beta) |.
\end{equation}
\end{defn}

The function $\lambda(\beta)$ is $T$--invariant. Moreover, in the
case of the classical continued fraction shift $T\beta= 1/\beta
-[1/\beta]$ on $[0,1]$, the Lyapunov exponent is given by
\begin{equation} \label{LyapT}
\lambda(\beta) = 2 \lim_{n\to \infty}
\frac{1}{n} \log q_n(\beta),
\end{equation}
with $q_n(\beta)$ the successive denominators of the continued
fraction expansion.

In particular, the Khintchine--L\'evy theorem shows
that, for almost all $\beta$'s (with respect to the Lebesgue
measure on $[0,1]$) the limit \eqref{LyapT} is equal to
\begin{equation} \label{lambda0}
\lambda(\beta)=\pi^2 /(6\log 2)=:\lambda_0.
\end{equation}
There is, however, an exceptional set in $[0,1]$ of Hausdorff
dimension $\dim_H=1$ but with Lebesgue measure zero where limit
defining the Lyapunov exponent does not exist.

As we will see later, in ``good cases'' the value $\lambda(\beta)$
can be computed from the spectrum of Perron--Frobenius operator of
the shift $T$.

The {\em Lyapunov spectrum} is introduced (cf.~\cite{PoWei}) by
decomposing the unit interval in level sets of the Lyapunov
exponent $\lambda(\beta)$ of \eqref{LyapExp}. Let $L_c = \{ \beta
\in [0,1] \, | \lambda(\beta)= c \in \R \}$. These sets provide a
$T$--invariant decomposition of the unit interval,
$$ [0,1] = \bigcup_{c \in \R} L_c \cup \{ \beta \in [0,1] \, |
\lambda(\beta) \text{ does not exist} \}. $$ These level sets are
uncountable dense $T$--invariant subsets of $[0,1]$, of varying
Hausdorff dimension \cite{PoWei}. The Lyapunov spectrum measures
how the Hausdorff dimension varies, as a function $h(c) = \dim_H
(L_c)$.

\bigskip

\bigskip

\subsubsection{Limiting modular symbols and iterated shifts.}

We introduce a function $\varphi: \P \to H^{\text{cusps}}$ of the
form
\begin{equation}
\varphi(s) = \{ g(0), g(i\infty) \}_G , \label{varphi}
\end{equation}
where $g\in \PGL(2,\Z)$ (or $\PSL(2,\Z)$) is a representative of
the coset $s\in \P$.

Then we can compute the limit \eqref{limitmod} in the following
way.

\begin{thm}
Consider a fixed $c\in \R$ which corresponds to some level set $L_c$ of the
Lyapunov exponent \eqref{LyapT}. Then, for all $\beta \in L_c$,
the limiting modular symbol \eqref{limitmod} is computed by the
limit
\begin{equation}\label{limitLyap}
\lim_{n\to \infty} \frac{1}{c n} \sum_{k=1}^n \varphi \circ T^k
(t_0),
\end{equation}
where $T$ is the shift of \eqref{shift} and $t_0 \in \P$.
\label{thmLyap}
\end{thm}

Without loss of generality, one can consider the geodesic
$\gamma_\beta$ in $\H\times \P$ with one end at $(i\infty, t_0)$
and the other at $(\beta, t_0)$, for $\varphi(t_0)= \{ 0,i\infty
\}_G$.

\begin{figure}
\begin{center}
\epsfig{file=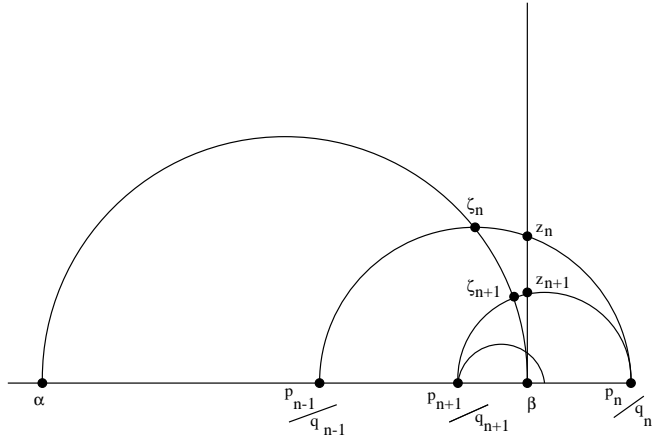}
\end{center}
\caption{Geodesics defining limiting modular symbols
\label{Fig-modsymb}}
\end{figure}

The argument given in \S 2.3 of \cite{ManMar} is based on the fact
that one can replace the homology class defined by the vertical
geodesic with one obtained by connecting the successive rational
approximations to $\beta$ in the continued fraction expansion
(Figure \ref{Fig-modsymb}). Namely, one can replace the path
$\langle x_0, y_n \rangle$ with the union of arcs
$$
\langle x_0, y_0 \rangle \cup \langle y_0, p_0/q_0 \rangle \cup
\bigcup_{k=1}^n \langle p_{k-1}/q_{k-1}, p_{k}/q_{k} \rangle \cup
\langle p_n/q_n, y_n \rangle
$$
representing the same homology class in $H_1(\overline{X}_G,\Z)$.

The result then follows by estimating the geodesic distance $\tau
\sim -\log \Im y + O(1)$, as $y(\tau) \to \beta$ and
$$ \frac{1}{2 q_n q_{n+1}} < \Im y_n < \frac{1}{2 q_n q_{n-1}}, $$
where $y_n$ is the intersection of $\gamma_\beta$ and the geodesic
with ends at $p_{n-1}(\beta)/q_{n-1}(\beta)$ and
$p_n(\beta)/q_n(\beta)$.

The matrix $g_k^{-1}(\beta)$, with
$$ g_k(\beta)= \left( \begin{array}{cc} p_{k-1}(\beta) & p_k(\beta) \\
q_{k-1}(\beta) & q_k(\beta) \end{array}\right), $$ acts on points
$(\beta,t)\in [0,1]\times \P$ as the $k$--th power of the shift
map $T$ of \eqref{shift}. Thus, we obtain
$$
\varphi(T^k t_0) = \{ g_k^{-1}(\beta)\, (0),g_k^{-1}(\beta)\,
(i\infty) \}_G = \left\{
\frac{p_{k-1}(\beta)}{q_{k-1}(\beta)},\frac{p_k(\beta)}{q_k(\beta)}
\right\}_G.
$$

\subsection{Ruelle and Perron--Frobenius operators}

A general principle in the theory of dynamical systems is that one
can often study the dynamical properties of a map $T$ (\eg
ergodicity) via the spectral theory of an associated operator.
This allows one to employ techniques of functional analysis and
derive conclusions on dynamics.

In our case, to the shift map $T$ of \eqref{shift}, we associate
the operator
\begin{equation}\label{PerFrob}
 (L_\sigma f) (x,t)= \sum_{k=1}^\infty
\frac{1}{(x+k)^{2\sigma}}\, f\left( \frac{1}{x+k},
\left(\begin{array}{cc} 0&1\\1&k
\end{array}\right) \cdot t \right)
\end{equation}
depending on a complex parameter $\sigma$.

More generally, the Ruelle transfer operator of a map $T$ is
defined as
\begin{equation}\label{Ruelle}
 (L_\sigma f) (x,t)= \sum_{(y,s)\in T^{-1}(x,t)} \exp(h(y,s))\,
f(y,s),
\end{equation}
where we take $h(x,t)=-2\sigma \log |T'(x,t)|$. Clearly this
operator is well suited for capturing the dynamical properties of
the map $T$ as it is defined as a weighted sum over preimages. On
the other hand, there is another operator that can be associated
to a dynamical system and which typically has better spectral
properties, but is less clearly related to the dynamics. The best
circumstances are when these two agree (for a particular value of
the parameter). The other operator is the Perron--Frobenius
operator ${\mathcal P}$. This is defined by the relation
\begin{equation}\label{Perron}
 \int_{[0,1]\times\P} f \, (g\circ T)\, d\mu_{Leb} =
\int_{[0,1]\times\P} ({\mathcal P} \, f) \, g\, d\mu_{Leb}.
\end{equation}
In the case of the shift $T$ of \eqref{shift} we have in fact that
$$ {\mathcal P} = L_\sigma |_{\sigma=1}. $$

\bigskip

\subsubsection{Spectral theory of $L_1$}

In the case of the modular group $G=\Gamma$, the spectral theory
of the Perron Frobenius operator of the Gauss shift was studied by
D.Mayer \cite{Mayer2}. More recently, Chang and Mayer
\cite{ChMayer} extended the results to the case of congruence
subgroups. A similar approach is used in \cite{ManMar} to study
the properties of the shift \eqref{shift}.

The Perron--Frobenius operator
$$ (L_1 f)(x,s)  = \sum_{k=1}^\infty
\frac{1}{(x+k)^2} f \left( \frac{1}{x +k}, \left(\begin{array}{cc}
0 & 1\\ 1 & k \end{array} \right) \cdot s \right) $$ for the shift
\eqref{shift} has the following properties.

\begin{thm}\label{spectralL1}
On a Banach space of holomorphic functions on $D\times\P$
continuous to boundary, with $D=\{z\in\C\,|\,|z-1|<3/2\}$,
under the condition (irreducibility)
\begin{equation}\label{irred-cond}
\P=\cup_{n=0}^{\infty} \left\{ \, \left( \begin{array}{cc}
0 & 1\\
1 & k_1
\end{array} \right)\ \dots \
\left( \begin{array}{cc}
0 & 1\\
1 & k_n
\end {array}\right)\,(t_0)\,|\,k_1,\dots , k_n\ge 1\right\},
\end{equation}
the Perron--Frobenius operator $L_1$ has the following properties:
\begin{itemize}
\item $L_1$ is a nuclear operator, of trace class.
\item $L_1$ has top eigenvalue $\lambda=1$. This eigenvalue is simple. The
corresponding eigenfunction is (up to normalization)
$$\frac{1}{(1+x)}.$$
\item The rest of the spectrum of $L_1$ is contained in a ball of radius
$r<1$.
\item There is a complete set of eigenfunctions.
\end{itemize}
\end{thm}

The irreducibility condition \eqref{irred-cond} is satisfied by
congruence subgroups.

\medskip

For other $T$-invariant subsets $E\subset[0,1]\times\P$, one can
also consider operators $L_{E,\sigma}$ and ${\mathcal P}_E$. When
the set has the property that ${\mathcal P}_E =L_{E,\delta_E}$,
for $\delta_E=\dim_H E$ the Hausdorff dimension, one can use the
spectral theory of the operator ${\mathcal P}_E$ to study the
dynamical properties of $T$.

The Lyapunov exponent can be read off the spectrum of the family
of operators $L_{E,\sigma}$.

\begin{lem}
Let $\lambda_\sigma$ denote the top eigenvalue of $L_{E,\sigma}$.
Then
$$ \lambda(\beta)= \frac{d}{d\sigma} \lambda_\sigma
|_{\sigma=\dim_H(E)} \ \ \   \mu_H \text{ a.e. in } E $$
\end{lem}

\bigskip

\subsubsection{The Gauss problem}

Let
\begin{equation}\label{mn-meas}
m_n(x):= \text{measure of}\ \{\alpha\in (0,1)\,|\,x_n(\alpha )\le
x\,\}
\end{equation}
with $\alpha = [\,a_1(\alpha ),\dots ,a_{n-1}(\alpha ), a_n(\alpha
)+x_n(\alpha )\,]$.

The asymptotic behavior of the measures $m_n$ is a famous problem
on the distribution of continued fractions formulated by Gauss,
who conjectured that
\begin{equation}\label{conjGauss}
m(x) = \lim_{n\to\infty}\,m_n(x) \stackrel{?}{=} \frac{1}{\log
2}\,\log(1+x).
\end{equation}
The convergence of \eqref{mn-meas} to \eqref{conjGauss} was only
proved by R.~Kuzmin in 1928. Other proofs were then given by
P.~L\'evy (1929), K.~Babenko (1978) and D.~Mayer (1991). The
arguments used by Babenko and Mayer use the spectral theory on the
Perron--Frobenius operator. Of these different arguments only the
latter extends nicely to the case of the generalized Gauss shift
\eqref{shift}.

The Gauss problem can be formulated in terms of a recursive
relation
\begin{equation}\label{Gauss-recursive}
m_{n+1}^{\prime}(x)= (L_1 \, m_n^{\prime})(x):=
\sum_{k=1}^{\infty} \frac{1}{(x+k)^2} m_n^{\prime}
\left(\frac{1}{x+k}\right).
\end{equation}
The right hand side of \eqref{Gauss-recursive} is the image of
$m_n'$ under the Gauss--Kuzmin operator. This is nothing but the
Perron--Frobenius operator for the shift $T$ in the case of the
group $\Gamma=\PGL(2,\Z)$.

As a consequence of Theorem \ref{spectralL1}, one obtains the
following result.

\begin{thm}\label{Gauss-gener}
Let
$$ m_n(x,t):= \mu_{Leb} \{ (y,s): x_n(y)\leq x,
g_n(y)^{-1}(s)=t \}. $$ Then we have
\begin{itemize}
\item $m_{n}^{\prime}(x,t)=L_1^n(1)$.
\item The limit $m(x,t) = \lim_{n\to\infty} m_n(x,t)$ exists
and equals
$$
m(x,t)=\frac{1}{|\P|\log 2} \log(1+x).
$$
\end{itemize}
\end{thm}

This shows that there exists a unique $T$-invariant measure on
$[0,1]\times \P$. This is uniform in the discrete set $\P$ (the
counting measure) and it is the Gauss measure of the shift of the
continued fraction expansion on $[0,1]$.

\bigskip

\subsection{Two theorems on limiting modular symbols}

The result on the $T$-invariant measure allows us to study the
general behavior of limiting modular symbols.

A special role is played by limiting modular symbols $\{\{ *,
\beta \}\}$ where $\beta$ is a quadratic irrationality in
$\R\smallsetminus \Q$.

\begin{thm}\label{limmodsym-per}
Let $g\in G$ be hyperbolic, with eigenvalue $\Lambda_g$
corresponding to the  attracting fixed point $\alpha_g^+$. Let
$\Lambda (g):=| \log \Lambda_g|$, and let $\ell$ be the period of
the continued fraction expansion of $\beta=\alpha_g^+$. Then
$$
\{\{*,\beta \}\}_G=\frac{\{0,g(0)\}_G}{\Lambda (g)} $$
$$=\frac{1}{\lambda(\beta)\ell} \sum_{k=1}^\ell \{
g_k^{-1}(\beta)\cdot 0,g_k^{-1}(\beta) \cdot i\infty \}_G .
$$
\end{thm}

This shows that, in this case, the limiting modular symbols are
linear combinations of classical modular symbols, with
coefficients in the field generated over $\Q$ by the Lyapunov
exponents $\lambda(\beta)$ of the quadratic irrationalities.

\begin{figure}
\begin{center}
\epsfig{file=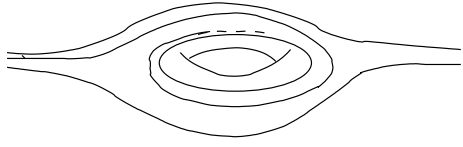}
\end{center}
\caption{Limiting modular symbols: limiting cycle
\label{Fig-curve1}}
\end{figure}

In terms of geodesics on the modular curve, this is the case where
the geodesic has a limiting cycle given by the closed geodesic
$\{0,g(0)\}_G$ (Figure \ref{Fig-curve1}).

\begin{figure}
\begin{center}
\epsfig{file=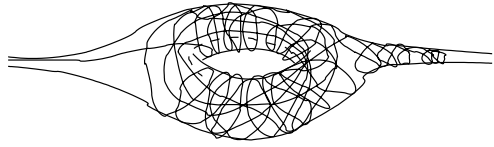}
\end{center}
\caption{Limiting modular symbols: chaotic tangling and untangling
\label{Fig-curve2}}
\end{figure}

There is then the ``generic case'', where, contrary to the
previous example, the geodesics wind around many different cycles
in such a way that the resulting homology class averages out to
zero over long distances (Figure \ref{Fig-curve2}).

\begin{thm}\label{limmodsym-gen}
For a $T$-invariant $E\subset [0,1]\times \P$, under the
irreducibility condition for $E$,
$$ R_\tau(\beta,s):=\frac{1}{\tau}\,\{x,y(\tau)\}_{G} $$
converges weakly to zero. Namely, for all $f\in L^1(E,d\mu_H)$,
$$ \lim_{\tau\to \infty}\int_E R_\tau(\beta,s)\, f(\beta,s)\,
d\mu_H(\beta,s) =  0. $$ This weak convergence can be improved to
strong convergence $\mu_H(E)$-almost everywhere. Thus, the
limiting modular symbol satisfies
$$ \{\{*,\beta\}\}_G \equiv 0 \
\ \ \text{ a.e. on } E. $$
\end{thm}

This results depends upon the properties of the operator $L_1$ and
the result on the $T$-invariant measure. In fact, to get the
result on the weak convergence (\cite{ManMar}) one notices that
the limit \eqref{limitLyap} computing limiting modular symbols can
be evaluated in terms of a limit of iterates of the
Perron--Frobenius operator, by
$$ \lim_n \frac{1}{\lambda_0\, n} \sum_{k=1}^n \int_{[0,1]\times \P}
f (\varphi \circ T^k) = \lim_n \frac{1}{\lambda_0 n} \sum_{k=1}^n
\int_{[0,1]\times \P} (L_1^k f) \varphi,
$$
where $\lambda(\beta)=\lambda_0$ a.e. in $[0,1]$.

By the convergence of $L_1^k 1$ to the density $h$ of the
$T$-invariant measure and of
$$ L_1^k f \to \left(\int f d\mu\right) \,\, h, $$
this yields
$$ \int_{[0,1]\times \P} \{\{ *, \beta \}\} \,\, f(\beta,t) \,\, d\mu(\beta,t) =
\left(\int_{[0,1]\times \P}  f\, d\mu\right)\,\, \left(
\int_{[0,1]\times \P} \varphi \,\, h \, d \mu \right) $$
$$ =\left(\int  f\, d\mu\right)\,\, \frac{1}{2\# \P} \sum_{s\in\P} \varphi(s).
$$
It is then easy to check that the sum $\sum_{s\in\P} \varphi(s)=0$
since each term in the sum changes sign under the action of the
inversion $\sigma\in \PGL(2,\Z)$ with $\sigma^2=1$, but the sum is
globally invariant under $\sigma$.

The argument extends to the case of other $T$-invariant subsets of
$[0,1]\times \P$ for which the corresponding Perron--Frobenius
operator $L_{E,\delta_E}$ has analogous properties (\cf
\cite{Mar-lyap}). The weak convergence proved by this type of
argument can be improved to strong convergence by applying the
strong law of large numbers to the ``random variables''
$\varphi_k= \varphi \circ T^k$ (\cf \cite{Mar-lyap} for details).
The result effectively plays the role of an ergodic theorem for
the shift $T$ on $E$.

\bigskip

\section{Hecke eigenforms}

A very important question is what happens to modular forms at the
noncommutative boundary of the modular curves. There is a variety
of phenomena in the theory of modular forms that hint to the fact
that a suitable class of ``modular forms'' survives on the
noncommutative boundary. Zagier introduced the term ``quantum
modular forms'' to denote this important and yet not sufficiently
understood class of examples. Some aspects of modular forms
``pushed'' to the noncommutative boundary were analyzed in
\cite{ManMar}, in the form of certain averages involving modular
symbols and Dirichlet series related to modular forms of weight 2.
We recall here the main results in this case.

We shall now consider the case of congruence subgroups $G=\Gamma_0(p)$,
for $p$ a prime.

Let $\omega$ be a holomorphic differential (also called a
differential of the first kind) on the modular curve
$X_{\Gamma_0(p)}=\Gamma_0(p)\backslash\H$. Let
$\Phi=\varphi^*(\omega) /dz$ be the pullback under the projection
$\varphi: \H \to X_{\Gamma_0(p)}$.

Let $\Phi$ be an eigenfunction for all the Hecke operators
$$ T_n = \sum_{d|n} \sum_{b=0}^{d-1} \left( \begin{array}{cc} n/d & b
\\ 0 & d \end{array}\right) $$
with $(p,n)=1$,
$$ T_n \Phi = c_n \Phi. $$
Then the L-function of $\omega$ is given by
$$ L_\omega(s) = -\frac{(2\pi)^s}{(2\pi i)\Gamma(s)}
\int_0^\infty \Phi(iy)\, y^{s-1} \, dy $$
$$ L_\omega(s) = \sum_{n=1}^\infty c_n \, n^{-s}. $$

There are many very interesting arithmetic properties of the
integrals of such Hecke eigenforms on modular symbols (Manin
\cite{Man-sym}).

In particular, one has (\cite{Man-sym}) the following relation
between $L_\omega$ and modular symbols:
\begin{equation}\label{Mansymbeq}
(\sigma(n)-c_n)\, L_\omega(1)= \sum_{d|n,\, b \text{ mod } d}
\int_{\{ 0, b/d\}} \omega,
\end{equation}
with $\sigma(n)=\sum_{d|n} d$, from which one obtains
\begin{equation}\label{Mansymbeq2}
\begin{array}{l}
\sum_{q: (q,p)=1} q^{-(2+t)} \sum_{q'\le q, (q,q')=1} \int_{\{ 0,
q'/q\}} \omega \\[2mm]
 = \left[ \frac{\zeta^{(p)} (1+t)}{\zeta^{(p)} (2+t)}
-\frac{L_{\omega}^{(p)}(2+t)}{\zeta^{(p)} (2+t)^2}
\right]\,\int_0^{i\infty}\phi^*(\omega)
\end{array}
\end{equation}

Here $L_{\omega}^{(p)}$ is the Mellin transform of $\Phi$ with
omitted Euler $p$--factor, and $\zeta(s)$ the Riemann zeta, with
corresponding $\zeta^{(p)}$, that is, $L_\omega^{(p)}(s) =\sum_{n:
(n,p)=1} c_n \, n^{-s}$ and $\zeta^{(p)}(s)= \sum_{n: (n,p)=1}
n^{-s}$.

\medskip

It is therefore interesting to study the properties of integrals
on limiting modular symbols (\cf \cite{ManMar}) and extensions of
\eqref{Mansymbeq} \eqref{Mansymbeq2}.

To this purpose, in \cite{ManMar} a suitable class of functions on
the boundary $\P^1(\R)\times \P$ (or just $[0,1]\times \P$) was
introduced.

\medskip

These are functions of the form
\begin{equation}\label{f-field}
\ell(f,\beta)=\sum_{k=1}^\infty f(q_k(\beta),q_{k-1}(\beta)).
\end{equation}
Here $f$ is a complex valued function defined on
pairs of coprime integers $(q,q')$ with $q\geq q'\geq 1$ and with
$f(q,q') =O(q^{-\epsilon})$ for some $\epsilon >0$, and
$q_k(\beta)$ are the successive denominators of the continued
fraction expansion of $\beta \in [0,1]$.

The reason for this choice of functions is given by the following
classical result of L\'evy (1929).

\begin{prop}
For a function $f$ as above, we have
\begin{equation} \label{Levy}
\int_0^1 \ell(f,\alpha )d\alpha = \sum{}^{\prime}\
\frac{f(q,q^{\prime})}{q(q+q^{\prime})}.
\end{equation}
Sums and integrals here converge absolutely and uniformly.
\end{prop}

We can interpret the summing over pairs of successive denominators
as a property that replaces modularity, when ``pushed to the
boundary''.

Through this class of functions it is possible recast certain
averages related to modular symbols on $X_G$, completely in terms
of function theory on the ``boundary'' space
$\Gamma\backslash(\P^1(\R)\times \P)$. A typical result of this
sort is the following (\cite{ManMar}).

Upon choosing
$$ f(q,q^{\prime})=\frac{q+q'}{q^{1+t}} \int_{\{ 0,q'/q\}_G} \omega, $$
one obtains
$$ \int_0^1 \ell(f,\beta )\,d\mu_{Leb}(\beta) =$$
$$  \left[\frac{\zeta(1+t)}{\zeta(2+t)}
-\frac{L_{\omega}^{(p)}(2+t)}{\zeta^{(p)} (2+t)^2}
\right]\,\int_0^{i\infty}\phi^*(\omega),
$$
which expresses the Mellin transform of a Hecke eigenform in terms
of a boundary average.

Analogous results exist for different choices of $f(q,q')$, for
Eisenstein series twisted by modular symbols, double logarithms at
roots of unity, etc.

In fact, this result also shows that, even when the limiting
modular symbol vanishes, one can still associate nontrivial
homology classes to the geodesics with irrational ends.

In fact, for the case of $G=\Gamma_0(p)$ let us consider
\begin{equation}\label{Cfbeta}
C(f,\beta):= \sum_{n=1}^\infty \frac{q_{n+1}(\beta)+
q_n(\beta)}{q_{n+1}(\beta)^{1+t}}\, \left\{
0,\frac{q_n(\beta)}{q_{n+1}(\beta)}\right\}_{\Gamma_0(p)}.
\end{equation}
This defines, for $\Re (t)>0$ and for almost all $\beta$, a
homology class in $H_1(X_G,\text{cusps},\R)$ such that the
integral average
$$
\ell(f,\beta)= \int_{C(f,\beta)} \omega
$$
on $[0,1]$ recovers Mellin transforms of cusp forms by
$$
\int_0^1\int_{C(f,\beta)} \omega \,d\mu_{Leb}(\beta) =
 \left[\frac{\zeta(1+t)}{\zeta(2+t)}
-\frac{L_{\omega}^{(p)}(2+t)}{\zeta^{(p)} (2+t)^2}
\right]\,\int_0^{i\infty}\phi^*(\omega).
$$

An estimate
$$ \frac{q_{n+1}(\beta)+
q_n(\beta)}{q_{n+1}(\beta)^{1+t}}\, \left\{
0,\frac{q_n(\beta)}{q_{n+1}(\beta)}\right\}_{G} \sim e^{-(5+2t)n
\lambda(\beta)} \sum_{k=1}^n \varphi\circ T^k (s) $$ shows that
the series involved are absolutely convergent.

\section{Selberg zeta function}

Infinite geodesics on $X_G$ are the image image of infinite
geodesics on $\H^2\times \P$ with ends on $\P^1(\R)\times\P$.
Thus, they can be coded by $(\omega^-,\omega^+,s)$, with
$(\omega^{\pm},s)=$ endpoints in $\P^1(\R)\times \{ s \}$, $s\in
\P$, $\omega^-\in (-\infty,-1]$, $\omega^+\in [0,1]$.

If $\omega^\pm$ are not cusp points, then we have infinite
continued fraction expansions of these endpoints,
$$ \begin{array}{ll}\omega^+ & = [k_0, k_1,
\ldots k_r, k_{r+1}, \ldots] \\
 \omega^- & =[k_{-1};k_{-2},\ldots, k_{-n}, k_{-n-1}, \ldots]
\end{array}, $$
hence the corresponding geodesics are coded by $(\omega,s)$, $s\in
\P$, where $\omega$ is a doubly infinite sequence
$$ \omega= \ldots k_{-(n+1)} k_{-n} \ldots k_{-1} k_0 k_1 \ldots
k_n \ldots $$

The equivalence relation of passing to the quotient by the group
action is implemented by the invertible (double sided) shift:
$$ T (\omega^+,\omega^-,s)=
  \left( \frac{1}{\omega^+} -\left[ \frac{1}{\omega^+}
\right],\frac{1}{\omega^- -[1/\omega^+]}, \left(\begin{array}{cc}
-[1/\omega^+]&1 \\ 1&0 \end{array}\right) \cdot s\right). $$

In particular, the closed geodesics in $X_G$ correspond to the
case where the endpoints $\omega^\pm$ are the attractive and
repelling fixed points of a hyperbolic element in the group. This
corresponds to the case where $(\omega,s)$ is a periodic point for
the shift $T$.

The Selberg zeta function is a suitable ``generating function''
for the closed geodesics in $X_G$. It is given by
$$ Z_G(s)= \prod_{m=0}^\infty \prod_\gamma \left( 1- e^{-(s+m)
\ell(\gamma)} \right), $$
where $\gamma$ is a primitive closed
geodesic of length $\ell(\gamma)$.

This can also be expressed in terms of the generalized Gauss shift
by the following (\cite{ChMayer}, \cite{ManMar})

\begin{thm}\label{SelbergZ-perrfro}
The Selberg zeta function for $G\subset \PGL(2,\Z)$ of finite
index is computed by the Fredholm determinant of the Ruelle
transfer operator \eqref{Ruelle},
$$ Z_G(s) = \det( 1-L_s ). $$
\end{thm}

An analogous result holds for finite index subgroups of
$\SL(2,\Z)$, where one gets $\det( 1-L_s^2)$. There are also
generalizations to other $T$-invariant sets $E$, where $\det
(1-L_{E,s})$ corresponds to a suitable ``Selberg zeta'' that only
counts certain classes of closed geodesics.

\section{The modular complex and $K$-theory of ${\rm C}^*$-algebras}

When we interpret the ``boundary'' $\P^1(\R)\times \P$ with the
action of $\Gamma$ as the noncommutative space ${\rm C}^*$-algebra
${\rm C}(\P^1(\R)\times \P)\rtimes \Gamma$, it is possible to
reinterpret some of the arithmetic properties of modular curves in
terms of these operator algebras. For instance, the modular
symbols determine certain elements in the $K$-theory of this ${\rm
C}^*$-algebra, and the modular complex and the exact sequence of
relative homology \eqref{rel-hom3} can be identified canonically
with the Pimsner six terms exact sequence for $K$--theory of this
${\rm C}^*$-algebra.

For $\Gamma=\PSL(2,\Z)$ and $G\subset \Gamma$ a finite index
subgroup, we use the notation $X:=\P^1(\R) \times \P$, and
$\Gamma_\sigma:=\langle \sigma \rangle=\Z/2$ and
$\Gamma_\tau:=\langle\tau\rangle=\Z/3$.

\begin{figure}
\begin{center}
\epsfig{file=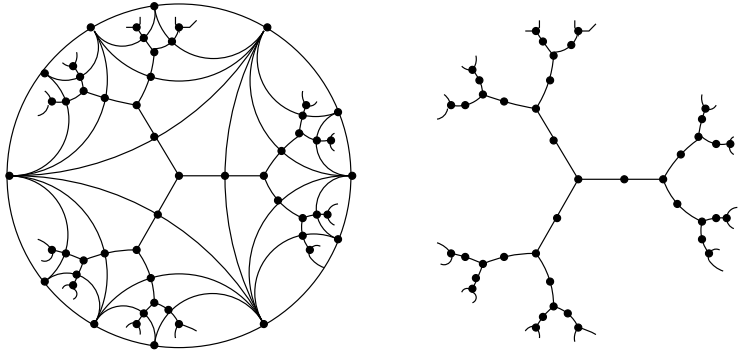}
\end{center}
\caption{The modular complex and the tree of $\PSL(2,\Z)$
\label{Fig-treemod}}
\end{figure}

The group $\Gamma$ acts on the tree with edges $T^1\simeq \Gamma$
and vertices $T^0\simeq \Gamma / \langle\sigma\rangle \cup \Gamma
/\langle\tau\rangle$. This tree is embedded in $\H^2$ with
vertices at the elliptic points and geodesic edges, as the dual
graph of the triangulation (Figure \ref{Fig-modcomplex}) of the
modular complex (Figure \ref{Fig-treemod}).

For a group acting on a tree there is a Pimsner six terms exact
sequence, computing the $K$-theory of the crossed product
$C^*$-algebra (for $\cA=C(X)$)
$$
\diagram K_0(\cA) \rto^{\alpha\qquad\qquad} & K_0(\cA\rtimes
\Gamma_\sigma) \oplus K_0(\cA\rtimes \Gamma_\tau)
\rto^{\qquad\qquad\tilde\alpha} &
K_0(\cA\rtimes\Gamma)\dto \\
K_1(\cA\rtimes\Gamma) \uto &  K_1(\cA\rtimes \Gamma_\sigma)
\oplus K_1(\cA\rtimes \Gamma_\tau) \lto^{\tilde\beta\qquad\qquad}
& K_1(\cA) \lto^{\qquad\qquad\beta}
\enddiagram
$$

A direct inspection of the maps in this exact sequence shows that
it contains a subsequence, which is canonically isomorphic to the
algebraic presentation of the modular complex, compatibly with the
covering maps between modular curves for different congruence
subgroups.

\begin{thm}\label{modC-PV}
The algebraic presentation of the modular complex
$$ 0 \to H_{cusps} \to \Z[\P] \to \Z[\P_R] \oplus \Z[\P_I] \to \Z
\to 0 $$ is canonically isomorphic to the sequence
$$ 0 \to Ker(\beta) \to \Z[\P] \to \Z[\P_I]\oplus \Z[\P_R] \to
Im(\tilde\beta) \to 0. $$
\end{thm}

A more geometric explanation for this exact sequence is the
following. The Lie group $\PSL(2,\R)$ can be identified with the
circle bundle (unit sphere bundle) over the hyperbolic plane
$\H=\PSL(2,\R)/PSO(2)$. Thus, the unit tangent bundle to the
modular curve $X_G$ is given by the quotient $T_1 X_G= G\backslash
\PSL(2,\R)$.

On the other hand, we can identify $\P^1(\R)=B\backslash
\PSL(2,\R)$, the quotient by the upper triangular matrices $B$, so
that we have Morita equivalent algebras
$$ C(\P^1(\R))\rtimes G \simeq C(T_1 X_G) \rtimes B. $$
The right hand side is the semidirect product of two $\R$-actions.
Thus, by the Thom isomorphism \cite{Co-Thom} we have
$$  K_*(C(T_1 X_G) \rtimes B)\cong K_*(C(T_1
 X_G)) $$
where $K^*(T_1 X_G)$ is related to $K^*(X_G)$ by the Gysin exact
sequence.

A construction of crucial importance in noncommutative geometry is
the classifying space for proper actions $\underline{B}G$ of
Baum--Connes \cite{BaumConnes}, the homotopy quotient $X\times_G 
\underline{E}G$ and the $\mu$-map
$$ \mu: K^*(X\times_G \underline{E}G)\to K_*(C(X)\rtimes G) $$
relating the topologically defined $K$-theory to the analytic
$K$-theory of $C^*$-algebras.

In a sense, for noncommutative spaces that are obtained as
quotients by ``bad'' equivalence relations, the homotopy quotient
is a ``commutative shadow'' from which much crucial information on
the topology can be read (\cf \cite{BaumConnes}, \cite{Co-transv}).

In our case we have $\underline{E}G =\H^2$ and the $\mu$-map is an
isomorphism (the Baum--Connes conjecture holds). By retracting
$\underline{E}G =\H^2$ to the tree of $\PSL(2,\Z)$ and applying
the Mayer--Vietoris sequence to vertices and edges one obtains the
Pimsner six terms exact sequence.

\section{Intermezzo: Chaotic Cosmology}

We digress momentarily on a topic from General
Relativity, which turns out to be closely related to the
``noncommutative compactification'' of the modular curve $X_G$
with $G$ the congruence subgroup $\Gamma_0(2)$ (\cite{ManMar}).

An important problem in cosmology is understanding how anisotropy
in the early universe affects the long time evolution of
space-time. This problem is relevant to the study of the beginning
of galaxy formation and in relating the anisotropy of the
background radiation to the appearance of the universe today.

We follow \cite{Barrow}  for a brief summary of
anisotropic and chaotic cosmology. The simplest significant
cosmological model that presents strong anisotropic properties is
given by the Kasner metric {\small
\begin{equation}\label{Kasner}
ds^2 = -dt^2 + t^{2p_1} dx^2 + t^{2p_2} dy^2 + t^{2p_3} dz^2,
\end{equation} }
where the exponents $p_i$ are constants satisfying $\sum
p_i=1=\sum_i p_i^2$. Notice that, for $p_i=d\log g_{ii}/d\log g$,
the first constraint $\sum_i p_i=1$ is just the condition that
$\log g_{ij}= 2\alpha \delta_{ij} + \beta_{ij}$ for a traceless
$\beta$, while the second constraint $\sum_i p_i^2=1$ amounts to
the condition that, in the Einstein equations written in terms of
$\alpha$ and $\beta_{ij}$, 
$$ \left(\frac{d\alpha}{dt}\right)^2= \frac{8\pi}{3} \left( T^{00}
+ \frac{1}{16\pi} \left( \frac{d\beta_{ij}}{dt}\right)^2 \right)
$$
$$ e^{-3\alpha} \frac{d}{dt} \left(
e^{3\alpha}\frac{d\beta_{ij}}{dt} \right) = 8\pi \left( T_{ij}
-\frac{1}{3} \delta_{ij} T_{kk} \right),
$$ 
the term $T^{00}$ is negligible with respect to the term
$(d\beta_{ij}/dt)^2/16\pi$, which is the ``effective energy
density'' of the anisotropic motion of empty space, contributing
together with a matter term to the Hubble constant.

Around 1970, Belinsky, Khalatnikov, and Lifshitz introduced a
cosmological model ({\em mixmaster universe}) where they allowed
the exponents $p_i$ of the Kasner metric to depend on a parameter
$u$,
\begin{equation}\label{u-param}
\begin{array}{l}
p_1=\frac{-u}{1+u+u^2} \\[3mm]
p_2=\frac{1+u}{1+u+u^2}\\[3mm]
p_3=\frac{u(1+u)}{1+u+u^2}
\end{array}
\end{equation}
Since for fixed $u$ the model is given by a Kasner space-time, the
behavior of this universe can be approximated for certain large
intervals of time by a Kasner metric. In fact, the evolution is
divided into Kasner eras and each era into cycles. During each era
the mixmaster universe goes through a volume compression. Instead
of resulting in a collapse, as with the Kasner metric, high
negative curvature develops resulting in a bounce (transition to a
new era) which starts again a behavior approximated by a Kasner
metric, but with a different value of the parameter $u$. Within
each era, most of the volume compression is due to the scale
factors along one of the space axes, while the other scale factors
alternate between phases of contraction and expansion. These
alternating phases separate cycles within each era.

More precisely, we are considering a metric
\begin{equation}\label{mixmaster}
ds^2 =-dt^2 + a(t) dx^2 + b(t) dy^2 + c(t) dz^2,
\end{equation}
generalizing the Kasner metric \eqref{Kasner}. We still require
that \eqref{mixmaster} admits $SO(3)$ symmetry on the space-like
hypersurfaces, and that it presents a singularity at $t\to 0$. In
terms of logarithmic time $d\Omega= -\frac{dt}{abc}$, the {\em
mixmaster universe} model of Belinsky, Khalatnikov, and Lifshitz
admits a discretization with the following properties:

\medskip

\noindent{\bf 1.} The time evolution is divided in Kasner eras
$[\Omega_n, \Omega_{n+1}]$, for $n\in \Z$. At the beginning of
each era we have a corresponding discrete value of the parameter
$u_n > 1$ in \eqref{u-param}.

\noindent{\bf 2.} Each era, where the parameter $u$ decreases with
growing $\Omega$, can be subdivided in cycles corresponding to the
discrete steps $u_n$, $u_n -1$, $u_n -2$, etc. A change $u\to u-1$
corresponds, after acting with the permutation $(12)(3)$ on the
space coordinates, to changing $u$ to $-u$, hence replacing
contraction with expansion and conversely. Within each cycle the
space-time metric is approximated by the Kasner metric
\eqref{Kasner} with the exponents $p_i$ in \eqref{u-param} with a
fixed value of $u=u_n -k >1$.

\noindent{\bf 3.} An era ends when, after a number of cycles, the
parameter $u_n$ falls in the range $0< u_n < 1$. Then the bouncing
is given by the transition $u\to 1/u$ which starts a new series of
cycles with new Kasner parameters and a permutation $(1)(23)$ of
the space axis, in order to have again $p_1< p_2< p_3$ as in
\eqref{u-param}.

\bigskip

Thus, the transition formula relating the values $u_n$ and
$u_{n+1}$ of two successive Kasner eras is
$$ u_{n+1} = \frac{1}{u_n - [u_n]}, $$
which is exactly the shift of the continued fraction expansion,
$Tx=1/x -[1/x]$, with $x_{n+1}=T x_n$ and $u_n = 1/x_n$.

The previous observation is the key to a geometric description of
solutions of the mixmaster universe in terms of geodesics on a
modular curve (\cf \cite{ManMar}).

\begin{thm}\label{thm-geod}
Consider the modular curve $X_{\Gamma_0(2)}$. Each infinite
geodesic $\gamma$ on $X_{\Gamma_0(2)}$ not ending at cusps
determines a mixmaster universe.
\end{thm}

In fact, an infinite geodesic on
$X_{\Gamma_0(2)}$ is the image under the quotient map
$$
\pi_\Gamma : \H^2 \times \P \to \Gamma\backslash (\H^2\times \P)
\cong X_G,
$$
where $\Gamma =\PGL(2,\Z)$, $G=\Gamma_0(2)$, and $\P=\Gamma/G\cong
\P^1(\F_2)=\{ 0,1,\infty \}$, of an infinite geodesic on
$\H^2\times \P$ with ends on $\P^1(\R)\times\P$. We consider the
elements of $\P^1(\F_2)$ as labels assigned to the three space
axes, according to the identification
\begin{equation}\label{axes}
\begin{array}{lll}
0=[0:1]& \mapsto & z \\
\infty=[1:0] & \mapsto & y \\
1=[1:1] & \mapsto & x.
\end{array}
\end{equation}

As we have seen, geodesics can be coded in terms of data
$(\omega^-,\omega^+, s)$ with the action of the shift $T$.

The data $(\omega,s)$ determine a mixmaster universe, with the
$k_n=[u_n]=[1/x_n]$ in the Kasner eras, and with the transition
between subsequent Kasner eras given by $x_{n+1}=Tx_n \in [0,1]$
and by the permutation of axes induced by the transformation
$$ \left( \begin{array}{cc} -k_n &1 \\ 1 & 0 \end{array}\right) $$
acting on $\P^1(\F_2)$. It is easy to verify that, in fact, this
acts as the permutation $0\mapsto \infty$, $1\mapsto 1$,
$\infty\mapsto 0$, if $k_n$ is even, and $0\mapsto \infty$,
$1\mapsto 0$, $\infty \mapsto 1$ if $k_n$ is odd, that is, after
the identification \eqref{axes}, as the permutation $(1)(23)$ of
the space axes $(x,y,z)$, if $k_n$ is even, or as the product of
the permutations $(12)(3)$ and $(1)(23)$ if $k_n$ is odd. This is
precisely what is obtained in the mixmaster universe model by the
repeated series of cycles within a Kasner era followed by the
transition to the next era.

Data $(\omega,s)$ and $T^m (\omega,s)$, $m\in \Z$, determine the
same solution up to a different choice of the initial time.

There is an additional time-symmetry in this model of the evolution of
mixmaster universes (\cf \cite{Barrow}). In fact, there is an
additional parameter $\delta_n$ in the system, which measures the initial
amplitude of each cycle. It is shown in \cite{Barrow} that this is
governed by the evolution of a parameter
$$ v_n = \frac{\delta_{n+1} (1+u_n)}{1-\delta_{n+1}} $$
which is subject to the transformation across cycles
$v_{n+1} = [u_n] + v_n^{-1}$.
By setting $y_n =v_n^{-1}$ we obtain
$$ y_{n+1} = \frac{1}{\left( y_n + \left[ 1/x_n \right] \right)},
$$
hence we see that we can interpret the
evolution determined by the data $(\omega^+,\omega^-,s)$ with the
shift $T$ either as giving the complete evolution of the
$u$-parameter towards and away from the cosmological singularity, or
as giving the simultaneous evolution of the two parameters $(u,v)$
while approaching the cosmological singularity.

This in turn determines the complete evolution of the parameters
$(u,\delta,\Omega)$, where $\Omega_n$ is the starting time of each
cycle. For the explicit recursion $\Omega_{n+1}=\Omega_{n+1}
(\Omega_n,x_n,y_n)$ see \cite{Barrow}.

The result of Theorem \ref{Gauss-gener} on the unique $T$-invariant 
measure on $[0,1]\times \P$
\begin{equation}\label{inv-meas}
 d\mu(x,s)= \frac{\delta(s)\, dx}{3\log(2)\, (1+x)},
\end{equation}
implies that the alternation of the space axes
is uniform over the time evolution, namely the three axes provide
the scale factor responsible for volume compression with equal
frequencies.

The Perron-Frobenius operator for the shift
\eqref{shift} yields the density of the invariant measure \eqref{inv-meas}
satisfying $L_1 f= f$. The top eigenvalue $\eta_\sigma$
of $L_{\sigma}$ is related to the topological pressure by
$\eta_\sigma=\exp( P(\sigma) )$. This can be estimated numerically,
using the technique of \cite{Babenko} and the integral kernel operator
representation of \S 1.3 of \cite{ManMar}.

The interpretation of solutions in terms of geodesics provides a
natural way to single out and study certain special classes of
solutions on the basis of their geometric properties. Physically, such
special classes of solutions exhibit different behaviors
approaching the cosmological singularity.

For instance, the data
$(\omega^+,s)$ corresponding to an eventually periodic sequence
$k_0 k_1 \ldots k_m \ldots$
of some period $\overline{a_0\ldots a_\ell}$ correspond to those
geodesics on $X_{\Gamma_0(2)}$ that asymptotically wind around the
closed geodesic identified with the doubly infinite sequence
$\ldots a_0\ldots a_\ell a_0\ldots a_\ell \ldots$.
Physically, these universes exhibit a pattern of cycles that
recurs periodically after a finite number of Kasner eras.

Another special class of solutions is given by the Hensley
Cantor sets (\cf \cite{Mar-cosm}). These are the mixmaster
universes for which there is a fixed upper bound $N$ to the number of
cycles in each Kasner era.

In terms of the continued fraction description, these solutions
correspond to data $(\omega^+,s)$ with $\omega^+$ in the Hensley
Cantor set $E_N\subset [0,1]$. The set $E_N$ is given by all points in
$[0,1]$ with all digits in the continued fraction expansion bounded by
$N$ (\cf \cite{Hen2}). In more geometric terms, these
correspond to geodesics on the modular curve $X_{\Gamma_0(2)}$ that
wander only a finite distance into a cusp.

On the set $E_N$ the Ruelle and Perron--Frobenius operators are given by 
\begin{equation}\label{RuelleEN}
(L_{\sigma,N} f)(x,s)  = \sum_{k=1}^N
\frac{1}{(x+k)^{2\sigma}} f \left( \frac{1}{x +k},
\left(\begin{array}{cc} 0 & 1\\ 1 & k \end{array} \right) \cdot s
\right).
\end{equation}

This operator still
has a unique invariant measure $\mu_N$, whose density satisfies
${\mathcal L}_{2\dim_H(E_N),N} f = f$, with
$$ \dim_H(E_N) = 1- \frac{6}{\pi^2 N} - \frac{72 \log N}{
\pi^4 N^2} + O(1/N^2) $$
the Hausdorff dimension of the Cantor set $E_N$.
Moreover, the top eigenvalue $\eta_\sigma$
of ${\mathcal L}_{\sigma,N}$ is related to the
Lyapunov exponent by
$$ \lambda(x) =  \frac{d}{d\sigma} \eta_\sigma|_{\sigma =\dim_H(E_N)}, $$
for $\mu_N$-almost all $x\in E_N$.

A consequence of this characterization of the time evolution in terms
of the dynamical system \eqref{shift} is that
we can study global properties of suitable {\em moduli spaces} of
mixmaster universes. For instance, the moduli space for time
evolutions of the $u$-parameter approaching the cosmological
singularity as $\Omega \to \infty$ is given by the quotient of
$[0,1]\times \P$ by the action of the shift $T$. 

Similarly when we restrict to special classes of solutions, \eg 
we can consider the ``moduli space'' $E_N \times \P$ modulo 
the action of $T$. In this example, the dynamical system $T$ acting on
$E_N \times \P$ is a {\em subshift of finite type}, and the 
resulting nocommutative space is a Cuntz--Krieger algebra \cite{CuKri}. 
This is an interesting class of $C^*$-algebras, generated by partial 
isometries. Another such algebra will play a fundamental role in
the geometry at arithmetic infinity, which is the topic of the 
next lecture. 

\begin{figure}
\begin{center}
\epsfig{file=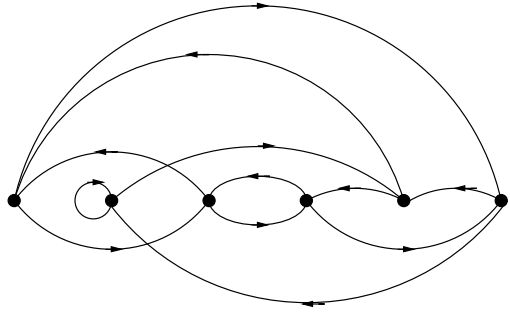}
\end{center}
\caption{The directed graph of the matrix $A_2$
\label{Fig-mixMatrix}}
\end{figure}

In the example of the mixmaster dynamics on the Hensley Cantor sets, the
shift $T$ is described by the Markov partition
$$ {\mathcal A}_N = \{ ((k,t),(\ell,s)) | U_{k,t} \subset
T(U_{\ell,s}) \},
$$
for $k,\ell \in \{ 1, \ldots, N \}$, and $s,t\in \P$, with sets
$U_{k,t} = U_k \times \{ t \}$, where $U_k \subset E_N$ are the clopen
subsets where the local inverses of $T$ are defined,
$$ U_k = \left[\frac{1}{k+1},\frac{1}{k}\right]\cap E_N. $$
This Markov partition determines a matrix $A_N$, with entries
$(A_N)_{kt,\ell s} =1$ if $U_{k,t} \subset T(U_{\ell,s})$ and zero
otherwise (\cf Figure \ref{Fig-mixMatrix}).

\begin{prop}\label{Matrix}
The $3\times 3$ submatrices $A_{k\ell}=( A_{(k,t),(\ell,s)} )_{s,t
\in\P}$ of the matrix $A_N$ are of the form
$$ A_{k\ell}=\left\{\begin{array}{lr} \left(\begin{array}{ccc}
0&0&1\\ 0&1&0\\ 1&0&0 \end{array}\right) & \ell =2m \\[6mm]
\left(\begin{array}{ccc}
0&0&1\\ 1&0&0\\ 0&1&0 \end{array}\right) & \ell =2m +1
\end{array}\right. $$
\end{prop}

This is just the condition $U_{k,t} \subset
T(U_{\ell,s})$ written as 
$$ \left(\begin{array}{cc} 0 & 1 \\ 1 & \ell \end{array} \right) \cdot s
=t , $$
together with the fact that the transformation
$$ \left(\begin{array}{cc} 0 & 1 \\ 1 & \ell \end{array} \right) $$
acts on $\P^1(\F_2)$ as the permutation $0\mapsto \infty$, $1\mapsto 1$,
$\infty\mapsto 0$, when $\ell$ is even, and $\infty \mapsto 0$,
$0\mapsto 1$, $1\mapsto \infty$ if $\ell$ is odd.

As a non-commutative space associated to the Markov partition we
consider the Cuntz--Krieger ${\rm C}^*$-algebra ${\mathcal O}_{A_N}$
(\cf \cite{CuKri}), which is the universal ${\rm C}^*$-algebra
generated by partial isometries $S_{kt}$ satisfying the relations
$$ \sum_{(k,t)}S_{kt}S_{kt}^* = 1, $$
$$ S_{\ell s}^* S_{\ell s} =\sum_{(k,t)} A_{(k,t),(\ell,s)}
S_{kt}S_{kt}^*. $$
Topological invariants of this ${\rm
C}^*$-algebra reflect dynamical properties of the shift $T$.

\chapter{Quantum statistical mechanics and Galois theory}

This chapter is dedicated to a brief review of 
the joint work of Alain Connes and the author on 
quantum statistical mechanics of $\Q$-lattices
\cite{CoMar}, \cite{cmln}. It also includes a discussion of the
relation to the explicit class field theory problem for imaginary
quadratic fields, following joint work of Alain Connes,
Niranjan Ramachandran, and the author \cite{CMR}, as well as a brief
review of some of Manin's idea on real quadratic fields and
noncommutative tori with real multiplication \cite{Man5}.

The main notion that we consider in this chapter is that of
commensurability classes of $\Q$-lattices, according to the
following definition.

\begin{defn}\label{defQlat}
A $\Q$-lattice in $\R^n$ consists of a pair $( \Lambda , \phi) $
of a lattice $\Lambda\subset \R^n$ (a cocompact free abelian
subgroup of $\R^n$ of rank $n$) together with a system of labels
of its torsion points given by a homomorphism of abelian groups
$$ \phi :  \Q^n/\Z^n \longrightarrow \Q\Lambda / \Lambda. $$
A $\Q$-lattice is {\em invertible}
if $\phi$ is an isomorphism. 

Two $\Q$-lattices are commensurable,
$$ (\Lambda_1, \phi_1) \sim (\Lambda_2, \phi_2), $$
iff $\Q\Lambda_1=\Q\Lambda_2$ and
$$ \phi_1 = \phi_2 \mod \Lambda_1 + \Lambda_2 $$
\end{defn}

One can check that, indeed, commensurability defines an equivalence
relation among $\Q$-lattices. The interesting aspect of this
equivalence relation is that the quotient provides another typical case
which is best described through noncommutative geometry. For this it
is crucial that we consider non-invertible $\Q$-lattices. In fact,
most $\Q$-lattices are not
commensurable to an invertible one, while
two invertible $\Q$-lattices are
commensurable if and only if they are equal. 

\begin{figure}
\begin{center}
\epsfig{file=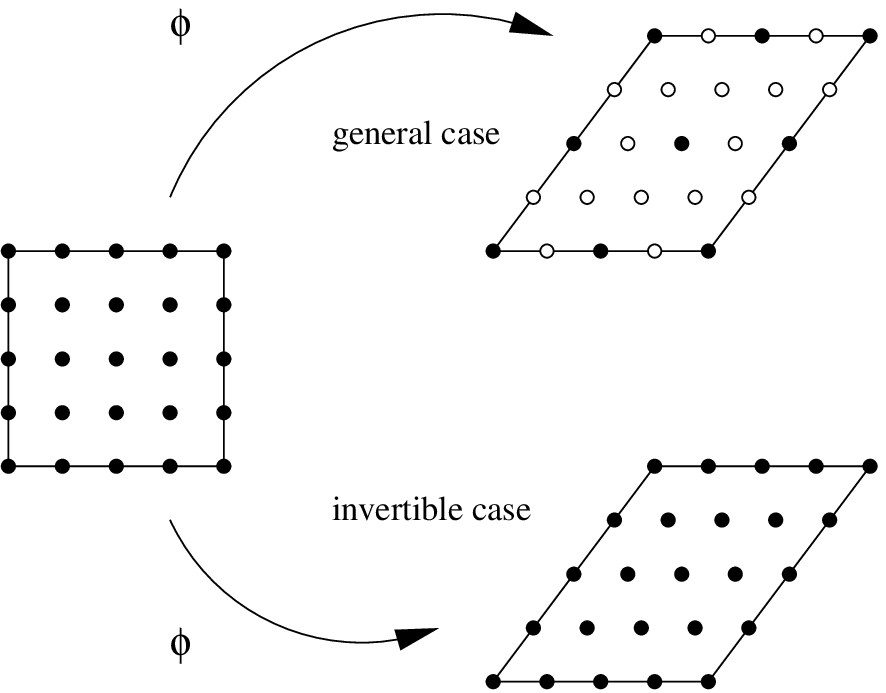}
\end{center}
\caption{Generic and invertible 2-dimensional $\Q$-lattices
\label{Fig-Qlat}}
\end{figure}

In the following, we denote by $\cL_n$ the set of commensurability classes of
$\Q$-lattices in $\R^n$. The topology on this space is encoded in a
noncommutative algebra of coordinates, namely a $C^*$-algebra $C^*(\cL_n)$.
We will discuss in some detail only the cases $n=1$ (Bost--Connes
system) and $n=2$ (the case considered in \cite{CoMar}). In these
cases, we will look at commensurability classes of
$\Q$-lattices up to scaling, namely the quotients $\cL_1/\R_+^*$ and
$\cL_2/\C^*$. The corresponding $C^*$-algebras 
$$ C^*(\cL_1/\R_+^*) \ \ \ \ \text{ and } \ \ \ \
C^*(\cL_2/\C^*) $$ 
are endowed with a natural time evolution (a 1-parameter 
family of automorphisms). Thus, one can consider them as 
quantum statistical mechanical systems and look for equilibrium
states, depending on a thermodynamic parameter $\beta$ (inverse
temperature). The interesting connection to arithmetic arises
from the fact that the action of symmetries of the system on
equilibrium states at zero temperature can be described in terms of
Galois theory. In the 1-dimensional case of Bost--Connes, this will
happen via the class field theory of $\Q$, namely the Galois theory of
the cyclotomic field $\Q^{cycl}$. In the two dimensional case, the
picture is more elaborate and involves the automorphisms of the
field of modular functions.

\section{Quantum Statistical Mechanics}

In classical statistical mechanics a state is a probability
measure $\mu$ on the phase space that assigns to each observable
$f$ an expectation value, in the form of an average
\begin{equation}\label{average-stat} \int f \, d\mu. \end{equation}
In particular, for a Hamiltonian system, the macroscopic thermodynamic
behavior is described via the Gibbs canonical ensemble. This is the
measure
\begin{equation}\label{Gibbs} d\mu_{G} = \frac{1}{Z} e^{-\beta H}
d\mu_{Liouville}, \end{equation} normalized by $Z= \int e^{-\beta
H} d\mu_{Liouville}$, with a thermodynamic
parameter $\beta =(k T)^{-1}$, for $T$ the temperature and $k$
the Boltzmann constant (which we can set equal to one in suitable
units).  

A quantum statistical mechanical system consists of the 
data of an algebra of observables (a $C^*$-algebra $\cA$), together
with a time evolution given as a 1-parameter family of automorphisms
$\sigma_t \in \Aut(\cA)$.
We refer to the pair $(\cA, \sigma_t)$ as a
$C^*$-dynamical system. These data describe the microscopic quantum
mechanical evolution of the system. 

The macroscopic thermodynamical properties are encoded in
the equilibrium states of the system, depending on the inverse
temperature $\beta$.  
While the Gibbs measure in the classical case is 
defined in terms of the Hamiltonian and the
symplectic structure on the phase space, the notion of
equilibrium state in the quantum statistical mechanical setting 
only depends on the algebra of
observables and its time evolution, and does
not involve any additional structure like the symplectic structure or
the approximation by regions of finite volume.

We first need to recall the notion of states. These can be thought of
as probability measures on noncommutative spaces.

\begin{defn}\label{defnstate}
Given a unital $C^*$-algebra $\cA$, a {\em state} on
$\cA$ is a linear functional
$\varphi: \cA \to \C$ satisfying normalization and positivity,
\begin{equation}\label{defstate}
 \varphi(1)=1, \ \ \ \varphi(a^*a)\geq 0. 
\end{equation}
When the $C^*$-algebra
 $\cA$ is non unital, the condition $\varphi (1) = 1 $
is replaced by $\| \varphi \|=1$ where
\begin{equation}
\label{eq1}
\|\varphi\|:= {\rm sup}_{x\in \cA,\|x\|\leq 1}|\varphi(x)| \,.
\end{equation}
Such states are restrictions of states on the unital
$C^*$-algebra $\tilde{\cA}$ obtained by adjoining a unit.
\end{defn}

Before giving the general definition of equilibrium states via the KMS
condition, we can see equilibrium states in the simple case of a
system with finitely many quantum degrees of freedom. In this case,
the algebra of observables is the algebra of operators in a
Hilbert space $\cH$ and the time evolution is given by
$\sigma_t(a)=e^{itH} \, a \, e^{-itH}$, where $H$ is a positive
self-adjoint operator such that $\exp(-\beta H)$ is trace class for
any $\beta >0$. 
For such a system, the analog of \eqref{Gibbs} is
\begin{equation}\label{GibbsQ}
\varphi(a) = \frac{1}{Z} \, \Tr \left(a\, e^{-\beta H} \right) \ \ \ \
\forall a\in A,
\end{equation}
with the normalization factor $Z=\Tr(\exp(-\beta H))$.

The Kubo--Martin--Schwinger
condition (KMS) (\cf \cite{BR}, \cite{Haag}, \cite{HHW}) 
describing equilibrium states of more general 
quantum statistical mechanical systems generalizes \eqref{GibbsQ}
beyond the range of temperatures where $\exp(-\beta H)$ is trace class.

\begin{defn}\label{KMSdefn}
Given a
$C^*$-dynamical system $(\cA,\sigma_t)$,  a state
$\varphi$ on $\cA$ satisfies the KMS condition at inverse
temperature $0<\beta<\infty$ iff, for all $a,b \in \cA$, there
exists a function $F_{a,b}(z)$ holomorphic on the strip $0< \Im(z)
<\beta$ continuous to the boundary and bounded, such that for all
$t\in \R$
\begin{equation}\label{KMScond}
  F_{a,b}(t)=\varphi(a\sigma_t(b)) \ \ \text { and } \ \
F_{a,b}(t+i\beta)=\varphi(\sigma_t(b)a). \end{equation}
KMS$_\infty$ states are  weak limits of KMS$_\beta$ states as
$\beta\to \infty$,
\begin{equation}\label{KMSinfty}
 \varphi_\infty (a) = \lim_{\beta\to \infty} \varphi_\beta(a) \ \
 \ \forall a \in \cA.
\end{equation}
\end{defn}

\begin{figure}
\begin{center}
\epsfig{file=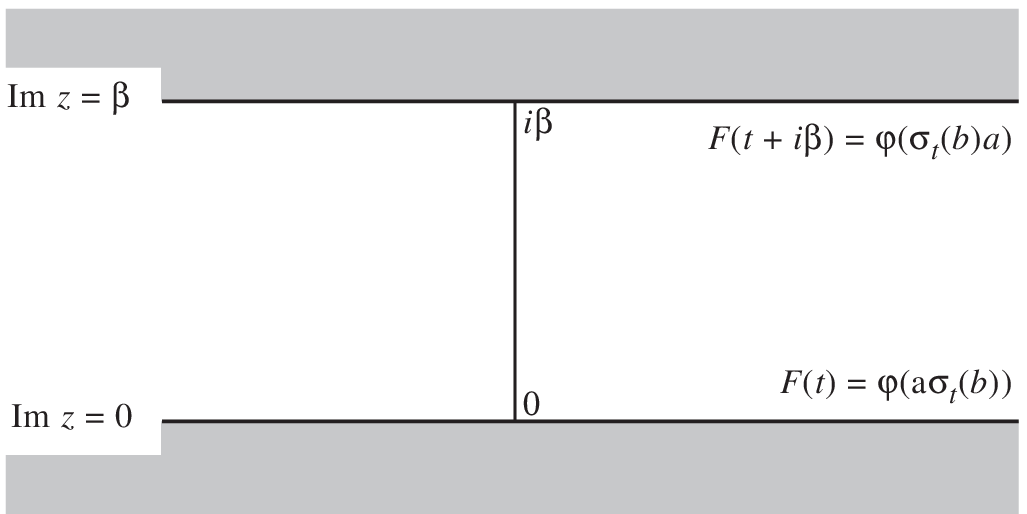}
\end{center}
\caption{The KMS condition \label{Fig-KMS}}
\end{figure}

Our definition of KMS$_\infty$ state is stronger than the one often 
adopted in the literature, which simply uses the existence, for each
$a,b\in \cA$, of a bounded 
holomorphic function $F_{a,b}(z)$ on the upper half plane such
that $F_{a,b}(t)=\varphi(a\sigma_t(b))$. It is easy to see that this
notions, which we simply call the {\em ground states} of the system,
is in fact weaker than the notion of KMS$_\infty$ states given in Definition
\ref{KMSdefn}. For example, in the simplest case of the trivial time
evolution,  all states are ground states, while only
tracial states are KMS$_\infty$ states.
Another advantage of our definition is that for all $0<\beta\leq
\infty$, the KMS$_\beta$ states form a Choquet simplex and we can
therefore consider the set $\cE_\beta$ of its extremal points.
These give a good notion of {\em points} for the underlying
noncommutative space.
This will be especially useful in connection to an arithmetic structure 
specified by an arithmetic subalgebra $\cA_\Q$ of $\cA$ on which
one evaluates the KMS$_\infty$ states. 
This will play a
key role in the relation between the symmetries of the system and
the action of the Galois group on states $\varphi \in \cE_\infty$
evaluated on $\cA_\Q$.

\subsection{Symmetries}

An important role in quantum statistical mechanics is played by
symmetries. Typically, symmetries of the algebra $\cA$ compatible
with the time evolution induce symmetries of the equilibrium
states $\cE_\beta$ at different temperatures. Especially important
are the phenomena of {\em symmetry breaking}. In such cases, there
is a global underlying group $G$ of symmetries of the algebra
$\cA$, but in certain ranges of temperature the choice of an
equilibrium state $\varphi$ breaks the symmetry to a smaller
subgroup $G_\varphi=\{ g \in G: \, g^* \varphi = \varphi \} $,
where $g^*$ denotes the induced action on states. Various systems
can exhibit one or more phase transitions, or none at all. A
typical situation in physical systems sees a unique KMS state for
all values of the parameter above a certain critical temperature
($\beta< \beta_c$). This corresponds to a chaotic phase such as
randomly distributed spins in a ferromagnet. When the system cools
down and reaches the critical temperature, the unique equilibrium
state branches off into a larger set KMS$_\beta$ and the symmetry
is broken by the choice of an extremal state in $\cE_\beta$. 

We will see that the case of $\cL_1/\R^*_+$ gives rise to a
system with a single phase transition (\cite{BC}), while in the
case of $\cL_2/\C^*$ the system has multiple phase
transitions.

\smallskip

A very important point is that we need to consider both symmetries
by automorphisms and by endomorphisms.

\smallskip

{\em Automorphisms:} A subgroup $G \subset \Aut(\cA)$ is
compatible with $\sigma_t$ if for all $g\in G$ and for all $t\in
\R$ we have $g\sigma_t = \sigma_t g$. There is then an induced
action of $G$ on KMS states and in particular on the set
$\cE_\beta$. If $u$ is a unitary, acting on $\cA$ by
$$ {\rm Ad} u: \,\,\,\, a\mapsto u a u^* $$
and satisfying $\sigma_t(u)=u$, then we say that ${\rm Ad} u$ is an inner
automorphism of $(\cA,\sigma_t)$. Inner automorphisms act
trivially on KMS states.

\smallskip

{\em Endomorphisms:} Let $\rho \sigma_t=\sigma_t \rho$ be a
$*$-homomorphism. Consider the idempotent $e=\rho(1)$. If
$\varphi\in \cE_\beta$ is a state such that $\varphi(e)\neq 0$,
then there is a well defined pullback $\rho^*\varphi$,
\begin{equation}\label{endoact}
 \rho^*(\varphi)= \frac{1}{\varphi(e)}\, \varphi \circ \rho.
\end{equation}
Let $u$ be an isometry compatible with the time evolution by
\begin{equation}\label{isomt}
\sigma_t(u)=\lambda^{it} u \ \ \ \lambda >0.
\end{equation}
One has $u^*u=1$ and $uu^*=e$.  We say that ${\rm Ad} u$ defined
by $a\mapsto u a u^*$ is an inner endomorphism of $(\cA,\sigma_t)$.
The condition \eqref{isomt} ensures that $({\rm Ad} u)^*\varphi$ is 
well defined according to \eqref{endoact} and the 
KMS condition shows that the induced
action of an inner endomorphism on KMS states is trivial.

\smallskip

One needs to be especially careful in defining the action of
endomorphisms by \eqref{endoact}. In fact, there are cases where
for KMS$_\infty$ states one finds only $\varphi(e)=0$, yet it is
still possible to define an interesting action of endomorphisms by
a procedure of ``warming up and cooling down''. For this to work
one needs sufficiently favorable conditions, namely that the
``warming up'' map
\begin{equation}\label{warming}
W_\beta(\varphi)(a)= \frac{ \Tr (\pi_\varphi(a)\,e^{-\beta
\,H})}{\Tr(\,e^{-\beta \,H})}
\end{equation}
gives a homeomorphism $W_\beta: \cE_\infty \to \cE_\beta$ for all
$\beta$ sufficiently large. One can then define the action by
\begin{equation}\label{warmcool}
 (\rho^* \varphi)(a) = \lim_{\beta\to \infty} \left( \rho^*
W_\beta(\varphi)\right)(a),
\end{equation}
for all $\varphi\in \cE_\infty$ and all $a\in\cA$.

\section{The Bost--Connes system}

We give the following geometric point of view on the BC system,
following \cite{CMR}. 

We use the notation $Sh(G,X):=G(\Q)\backslash (G(\A_f)\times X)$ for 
Shimura varieties (\cf \cite{Mi}). Here $\A_f=\hat\Z\otimes \Q$
denotes the finite ad\`eles of $\Q$, with $\hat\Z =\varprojlim_n
\Z/n\Z$. The simplest case is for $G=\GL_1$,
where we consider
\begin{equation}\label{ShGL1}
Sh(\GL_1,\{\pm 1\}) = \GL_1(\Q) \backslash (\GL_1(\A_f)\times \{
\pm 1 \}) = \Q^*_+\backslash \A_f^* .
\end{equation}

\begin{lem}\label{GL1Qlat}
The quotient \eqref{ShGL1} parameterizes the {\em invertible}
1-dimensional $\Q$-lattices up to scaling. 
\end{lem}

In fact, first notice that a 1-dimensional $\Q$-lattice can always be
written in the form
\begin{equation}\label{1dimQlat}
 ( \Lambda , \phi) \, = (\lambda\, \Z,\lambda\,\rho)
\end{equation}
for some $\lambda>0$ and some
\begin{equation}\label{rho}
\rho \in \Hom(\Q/\Z,\Q/\Z)=\varprojlim \Z/n\Z = \hat\Z.
\end{equation}
Thus, the set of 1-dimensional $\Q$-lattices up to scaling can be
identified with $\hat\Z$. The invertible lattices correspond to
the elements of $\hat \Z^*$, which in turn is identified with
$\GL_1^+(\Q) \backslash \GL_1(\A_f)$.

The quotient \eqref{ShGL1} can be thought of geometrically as
the Shimura variety associated to the {\em cyclotomic tower} (\cf
\cite{CMR}).  This is the tower $\cV$ of arithmetic varieties over
$V=\Sp(\Z)$, with $V_n= \Sp(\Z[\zeta_n])$, where $\zeta_n$ is a
primitive n-th root of unity. The group of deck transformations is
$\Aut_V(V_n)=\GL_1(\Z/n\Z)$, so that the group of deck transformations
of the tower is the projective limit
\begin{equation}\label{symmetries}
\Aut_V(\cV)= \varprojlim_n \Aut_V(V_n) = \GL_1(\hat\Z).
\end{equation}

If, instead of invertible $\Q$-lattices up to scale, we
consider the set of commensurability classes of all 1-dimensional
$\Q$-lattices up to scale, we find a noncommutative version of
the Shimura variety 
\eqref{ShGL1} and of the cyclotomic tower. This is described as the
quotient 
\begin{equation}\label{Sh1nc}
Sh^{(nc)}(\GL_1,\{ \pm 1 \}):= \GL_1(\Q)\backslash (\A_f \times
\{\pm 1 \})= \GL_1(\Q)\backslash \A^\cdot /\R^*_+ ,
\end{equation}
where $\A^\cdot:=\A_f\times \R^*$ are the ad\`eles of $\Q$ with
invertible archimedean component. 

The corresponding noncommutative
algebra of coordinates is given by the crossed product (\cf
\cite{laca-end}) 
\begin{equation}\label{adelic-cross}
C_0(\A_f) \rtimes \Q^*_+.
\end{equation}

One obtains an equivalent description of the noncommutative space
\eqref{Sh1nc}, starting from the description \eqref{1dimQlat} of
1-dimensional $\Q$-lattices. The commensurability relation is 
implemented by the action of $\N^\times =\Z_{>0}$ by
\begin{equation}\label{Nact}
 \alpha_n(f) (\rho)=\left\{ \begin{array}{lr} f(n^{-1} \rho) & \rho \in
 n\hat\Z \\ 0 & \text{otherwise.} \end{array}\right.
\end{equation}
Thus, the noncommutative algebra of coordinates of the space of
commensurability classes of 1-dimensional $\Q$-lattices up to scaling 
can be identified with the semigroup crossed product
\begin{equation}\label{semicross}
\cA_1:= C^*(\Q/\Z)\rtimes \N^\times ,
\end{equation}
where we have identified 
\begin{equation}\label{ChatZ}
C(\hat\Z)\simeq C^*(\Q/\Z),
\end{equation}
since $\hat\Z$ is the Pontrjagin dual of $\Q/\Z$. The algebra
\eqref{semicross} is Morita equivalent to the algebra
\eqref{adelic-cross} (\cf \cite{laca-end}). 

The algebra \eqref{semicross} is isomorphic to the algebra considered
in the work of Bost--Connes \cite{BC}, which was obtained there as a Hecke
algebra for the inclusion $\Gamma_0\subset \Gamma$ of the pair of groups
$(\Gamma_0, \Gamma)=(P_\Z, P_\Q)$, where $P$ is the $ax+b$ group.
These have the property that the left $\Gamma_0$ orbits of any
$\gamma\in \Gamma/\Gamma_0$ are finite (same for right orbits on the
left coset). The ratio of the lengths of left and right $\Gamma_0$
orbits determines a canonical time evolution $\sigma_t$ on the algebra (\cf
\cite{BC}). 

The algebra of the Bost--Connes system has an explicit
presentation in terms of two sets of operators. The first type
consists of {\em phase operators} $e(r)$, parameterized by elements
$r\in \Q/\Z$. 
These phase operators can be represented on the Fock space generated
by occupation numbers $|n\rangle$ as the operators
\begin{equation}\label{phaseOP}
e(r) |n\rangle =\alpha(\zeta_r^n) |n\rangle .
\end{equation}
Here we denote by $\zeta_{a/b}=\zeta_b^a$ the {\em abstract} roots
of unity generating $\Q^{cycl}$ and by $\alpha: \Q^{cycl}
\hookrightarrow \C$ an embedding that identifies $\Q^{cycl}$ with
the subfield of $\C$ generated by the {\em concrete} roots of
unity.

\begin{figure}
\begin{center}
\epsfig{file=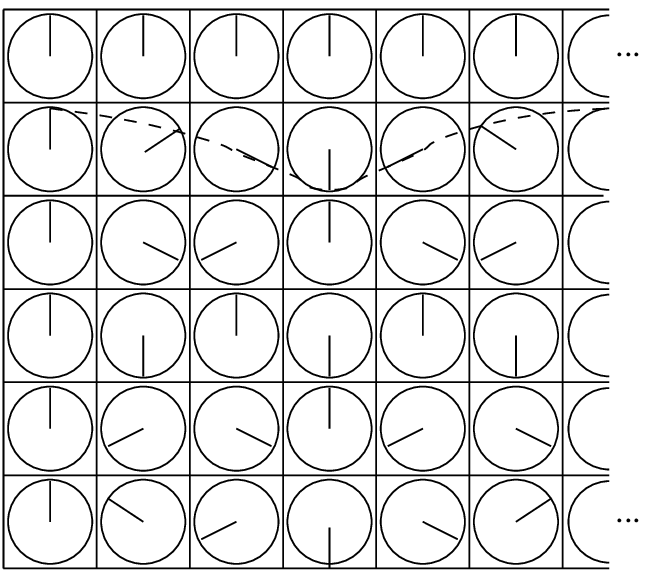}
\end{center}
\caption{Phasors with $\Z/6\Z$ discretization \label{Fig-phase}}
\end{figure}

These are the phase operators used in the theory of
quantum optics and optical coherence to model the phase
quantum-mechanically (\cf \cite{Lou} \cite{mandel}), as well as 
to model the {\em phasors} in quantum computing. They are based
on the choice of a certain scale $N$ at which the phase is discretized (Figure
\ref{Fig-phase}). Namely, the quantized
optical phase is defined as a state
$$ | \theta_{m,N}\rangle = e\left(
\frac{m}{N+1}\right) \cdot v_N, $$ where $v_N$ is a superposition
of occupation states
$$ v_N =\frac{1}{(N+1)^{1/2}}\sum_{n=0}^N | n
\rangle. $$ 
One needs then to ensure that the results are consistent over changes
of scale. 

The other operators that generate the Bost--Connes algebra can be
thought of as implementing the changes of scales in the optical
phases in a consistent way. These operators are isometries $\mu_n$
parameterized by positive integers $n\in \N^\times=\Z_{>0}$. The
changes of scale are described by the action of the $\mu_n$ on the
$e(r)$ by
\begin{equation}\label{muner}
 \mu_n e(r) \mu_n^* = \frac{1}{n} \sum_{ns=r} e(s).
\end{equation}
In addition to this compatibility condition, the operators $e(r)$
and $\mu_n$ satisfy other simple relations.

These give a presentation of the algebra of the BC system of
the form (\cite{BC}, \cite{Laca}):
\begin{itemize}
\item $\mu_n^*\mu_n =1$, for all $n\in \N^\times$,
\item $ \mu_k \mu_n=\mu_{kn} $, for all $k,n\in \N^\times$,
\item $e(0)=1$, $e(r)^*=e(-r)$, and $e(r)e(s)=e(r+s)$, for all $r,s\in
\Q/\Z$,
\item For all $n\in \N^\times$ and all $r\in \Q/\Z$, the relation
\eqref{muner} holds.
\end{itemize}
This shows that the Bost--Connes algebra is isomorphic to the algebra
$\cA_1$ of \eqref{semicross}.   

In terms of this explicit presentation, the time evolution is of the
form 
\begin{equation}\label{sigmatmue}
\sigma_t(\mu_n) = n^{it} \mu_n, \ \ \ \sigma_t(e(r))=e(r).
\end{equation}

The space \eqref{Sh1nc} can be compactified by replacing $\A^\cdot$ by
$\A$, as in \cite{Connes-Zeta}. This gives the quotient
\begin{equation}\label{compSh1nc}
 \overline{Sh^{(nc)}}(\GL_1,\{ \pm 1 \}) = \GL_1(\Q)\backslash \A
/\R^*_+.
\end{equation}
This compactification consists of adding the trivial lattice (with
a possibly nontrivial $\Q$-structure).

The dual space to \eqref{compSh1nc}, under the duality of type II and
type III factors introduced in Connes' thesis, is a principal $\R^*_+$-bundle
over \eqref{compSh1nc}, whose noncommutative algebra of coordinates is
obtained from the algebra of \eqref{compSh1nc} by
taking the crossed product by the time evolution $\sigma_t$. 
The space obtained this way is the space of ad\`ele classes
\begin{equation}\label{Lspace}
\cL_1=\GL_1(\Q)\backslash \A \to \GL_1(\Q)\backslash \A /\R^*_+
\end{equation}
that gives the spectral realization of zeros of the Riemann
$\zeta$ function in \cite{Connes-Zeta}. Passing to this dual space
corresponds to considering commensurability classes of
1-dimensional $\Q$-lattices (no longer up to scaling).

\subsection{Structure of KMS states}

The Bost--Connes algebra $\cA_1$ has irreducible representations on the
Hilbert space $\cH=\ell^2(\N^\times)$. These are parameterized by
elements $\alpha\in \hat\Z^*=\GL_1(\hat\Z)$. Any such element
defines an embedding $\alpha: \Q^{cycl}\hookrightarrow \C$ and the
corresponding representation is of the form
\begin{equation}\label{repGL1}
\begin{array}{rl}
\pi_\alpha (e(r))\, \epsilon_k = & \alpha(\zeta_r^k) \, \epsilon_k
\\[2mm]
 \pi_\alpha (\mu_n)\,  \epsilon_k = & \epsilon_{nk}
\end{array}
\end{equation}

The Hamiltonian implementing the time evolution \eqref{Nact} on
$\cH$ is
\begin{equation}\label{HamGL1}
H\,\epsilon_k=\log k \,\,\, \epsilon_k
\end{equation}
Thus, the partition function of the Bost--Connes system is the
Riemann zeta function
\begin{equation}\label{partzeta}
 Z(\beta)=\Tr\left( e^{-\beta H} \right) = \sum_{k=1}^\infty
k^{-\beta} = \zeta(\beta).
\end{equation}

Bost and Connes showed in \cite{BC} that KMS states have the
following structure, with a phase transition at $\beta=1$.
\begin{itemize}
\item In the range $\beta \leq 1$ there is a unique KMS$_\beta$
state. Its restriction to $\Q[\Q/\Z]$ is of the form
$$ \varphi_{\beta}(e(a/b))= b^{-\beta} \prod_{p \text{ prime},\,
p| b} \frac{1-p^{\beta-1}}{1-p^{-1}}. $$
\item For $1<\beta \leq \infty$ the set of extremal KMS
states $\cE_\beta$ can be identified with $\hat\Z^*$. It has a
free and transitive action of this group induced by an action on
$\cA$ by automorphisms. The extremal KMS$_\beta$ state
corresponding to $\alpha \in \hat\Z^*$ is of the form
\begin{equation}\label{statesalpha}
 \varphi_{\beta,\alpha}(x)= \frac{1}{\zeta(\beta)}\,
\Tr\left(\pi_\alpha(x)\, e^{-\beta H} \right).
\end{equation}
\item At $\beta=\infty$ the Galois group $\Gal(\Q^{cycl}/\Q)$ acts
on the values of states $\varphi\in \cE_\infty$ on an arithmetic
subalgebra $\cA_{1,\Q}\subset \cA_1$. These have the property that
$\varphi(\cA_\Q)\subset \Q^{cycl}$ and that the isomorphism (class
field theory isomorphism) $\theta: \Gal(\Q^{cycl}/\Q)
\stackrel{\cong}{\to} \hat\Z^*$ intertwines the Galois action on
values with the action of $\hat\Z^*$ by symmetries, namely,
\begin{equation}\label{GaloisEq1}
 \gamma^{-1} \, \varphi(x) = \varphi( \theta(\gamma)\, x),
\end{equation}
for all $\varphi\in \cE_\infty$, for all $\gamma\in
\Gal(\Q^{cycl}/\Q)$ and for all $x\in \cA_\Q$.
\end{itemize}

Here the arithmetic subalgebra $\cA_{1,\Q}$ can be taken to be the
algebra over 
$\Q$ generated by the $e(r)$ and the $\mu_n, \mu_n^*$. This $\Q$-algebra
can also be obtained (as shown in \cite{CoMar}) as the algebra 
generated by the
$\mu_n, \mu_n^*$ and by homogeneous functions of weight zero on 1-dim
$\Q$-lattices obtained as a normalization of the functions
\begin{equation}\label{epsilons}
\epsilon_{k,a} ( \Lambda , \phi)= \sum_{y\in  \Lambda +\phi(a)}
y^{-k}
\end{equation}
by covolume. Namely, one considers the functions $e_{k,a}:= c^k \,
\epsilon_{k,a}$, where $c(\Lambda)$ is proportional to the
covolume $|\Lambda|$ and satisfies
$$ (2\pi \sqrt{-1})\, c(\Z)=1. $$

The choice of an ``arithmetic subalgebra'' corresponds to endowing
the noncommutative space $\cA$ with an arithmetic structure. The
subalgebra corresponds to the rational functions and the values of
KMS$_\infty$ states at elements of this subalgebra should be
thought of as values of ``rational functions'' at the classical points
of the noncommutative space
(\cf \cite{CMR}).

\section{Noncommutative Geometry and Hilbert's 12th problem}

The most remarkable arithmetic feature of the result of Bost--Connes
recalled above is the Galois action on the ground states of the
system. First of all, the fact that the Galois action
on the values of states would preserve positivity (\ie would give
values of other states) is a very unusual property. Moreover,
the values of extremal ground states on elements of
the rational subalgebra generate the maximal abelian extension
$\Q^{ab}=\Q^{cycl}$ of $\Q$. The explicit action of the Galois
group $\Gal(\Q^{ab}/\Q)\simeq \GL_1(\hat\Z)$ is given by automorphism
of the system $(\cA_1,\sigma_t)$. Namely, the class field
theory isomorphism intertwines the two actions of the id\`eles
class group, as symmetry group of the system, and of the Galois
group, as permutations of the expectation values of the rational
observables.

In general, the  main theorem of class field theory provides a
classification of finite abelian extensions of a local or global
field $\K$ in terms of subgroups of a locally compact abelian group
canonically associated to the field. This is the multiplicative
group $\K^*=\GL_1(\K)$ in the local nonarchimedean case, while in the
global case
it is the quotient $C_\K/D_\K$ of the id\`ele class group $C_\K$ by the
connected component of the identity. 

Hilbert's 12th problem can be formulated as the question
of providing an explicit set of generators of the maximal abelian
extension $\K^{ab}$ of a number field $\K$, inside an algebraic
closure $\bar \K$, and of the action of the Galois group
$\Gal(\K^{ab}/\K)$. The maximal abelian extension $\K^{ab}$ of $\K$
has the property that
$$ \Gal(\K^{ab}/\K)=\Gal(\bar\K/\K)^{ab}. $$

The result that motivated Hilbert's formulation of the explicit class
field theory problem, is the Kronecker--Weber theorem, namely the
already mentioned case of the explicit class field theory of $\K=\Q$.
In this case, the maximal abelian extension
of $\Q$ can be identified with the cyclotomic field
$\Q^{cycl}$. Equivalently, one can say that the torsion points of the
multiplicative group $\C^*$ (\ie the roots of unity) generate
$\Q^{ab}\subset \C$.  

Remarkably, the only other case for number fields where
this program has been carried out completely is that of imaginary
quadratic fields, where the construction relies on the theory of
elliptic curves with complex multiplication, and on the Galois theory
of the field of modular functions (\cf \eg the survey
\cite{Steven}). 

Some generalizations of the original result of Bost--Connes to
other global fields (number fields and function fields) were obtained 
by Harari and Leichtnam \cite{HaLe},  P.~Cohen \cite{Cohen}, Arledge,
Laca and  Raeburn \cite{alr}. A more detailed account of various
results related to the BC system and generalizations is given in the
``further developments'' section of \cite{CoMar}. 
Because of this close relation to Hilbert's 12th problem, 
it is clear that finding generalizations of the Bost--Connes system to
other number fields is a very difficult problem, hence it is not
surprising that, so far, these constructions have not fully recovered
the Galois properties of the ground states of the BC system in the
generalized setting.

The strongest form of the result one may wish to obtain in this
direction can be formulated as follows. 

Given a number field $\K$,
we let $\A_\K$ denote the ad\`eles of $\K$ and we let
$\A_\K^*=\GL_1(\A_\K)$ be the group of id\`eles of $\K$. 
As above, we write $C_\K$ for
the group of id\`eles classes $C_\K=\A_\K^* /\K^*$ and $D_\K$ for the
connected component of the identity in $C_\K$.

One wishes to construct a $C^*$-dynamical system $(\cA,\sigma_t)$ and an
``arithmetic'' subalgebra $\cA_\K$, which satisfy the following properties:
\begin{enumerate}
\item The id\`eles class group $G=C_\K/D_\K$
acts by symmetries on $(\cA,\sigma_t)$ preserving the subalgebra
$\cA_\K$.
\item The states $\varphi \in \sE_\infty$, evaluated on elements of
$\cA_\K$, satisfy:
\begin{itemize}
\item $\varphi(a)\in \bar \K$, the algebraic
closure of $\K$ in $\C$;
\item the
elements of $\{ \varphi(a): a\in \cA_\K  \}$, for $\varphi\in
\sE_\infty$ generate $\K^{ab}$.
\end{itemize}
\item The class field theory isomorphism
$$
\theta:C_\K/D_\K \stackrel{\simeq}{\longrightarrow} \Gal (\K^{ab}/\K)
$$
intertwines the actions,
\begin{equation}\label{CFTiso}
 \alpha^{-1} \circ \varphi = \varphi \circ \theta^{-1}(\alpha),
\end{equation}
for all $\alpha \in \Gal (\K^{ab}/\K)$ and for all $\varphi \in
\sE_\infty$.
\end{enumerate}

This may provide a possible new approach, via noncommutative geometry,
to the explict class field theory problem.  
 
Given a number field $\K$ with $[\K:\Q]=n$, there is an embedding
$\K^*\hookrightarrow \GL_n(\Q)$ of its multiplicative group in
$\GL_n(\Q)$. Such embedding induces an embedding of
$\GL_1(\A_{\K,f})$ into $\GL_n(\A_f)$, where $\A_{\K,f}=\A_f\otimes
\K$ are the finite ad\`eles of $\K$.

This suggests that a possible strategy to approach the problem stated
above may be to first study quantum statistical mechanical systems
corresponding to $\GL_n$ analogs of the Bost--Connes system. 

The main result of the joint work of Alain Connes and the author in
\cite{CoMar} is the construction of such system in the $\GL_2$ case
and the analysis of the arithmetic properties of its KMS states.
In the case of $\GL_2$, one sees that
the geometry of modular curves and the algebra of modular forms
appear naturally. This leads to a formulation of a 
quantum statistical mechanical system related to the explicit
class field theory of imaginary quadratic fields \cite{CMR}.

The first case for which there is not yet a complete solution to the 
explicit class field theory problem is the case of real quadratic
fields, $\K=\Q(\sqrt{d})$, for some positive integer $d$.
It is natural to ask whether the approach outlined above, based on
noncommutative geometry, may provide any new information on this
case. Recent work of Manin \cite{Man3} \cite{Man5} suggests a close
relation between the real quadratic case and noncommutative
geometry. We will discuss the possible relation between his approach
and the $\GL_2$-system.

\section{The $\GL_2$ system}

In this section we will describe the main features of the $\GL_2$
analog of the Bost--Connes system, according to the results of
\cite{CoMar}.

We can start with the same geometric approach in terms of Shimura
varieties (\cf \cite{CMR}, \cite{cmln}) that we used above to introduce
the BC system. 

The Shimura variety associated to the tower of modular curves is
described by the ad\`elic quotient
\begin{equation}\label{ShiMod}
\begin{array}{rl}
Sh(\GL_2,\H^\pm) = & \GL_2(\Q)\backslash (\GL_2(\A_f)\times \H^\pm)
\\[3mm]
= & \GL_2^+(\Q)\backslash (\GL_2(\A_f)\times \H) =
\GL_2^+(\Q)\backslash \GL_2(\A) /\C^*.
\end{array}
\end{equation}

The tower of modular curves (\cf \cite{CMR}) has base $V=\P^1$ over
$\Q$ and $V_n=X(n)$ the modular curve corresponding to the principal
congruence subgroup $\Gamma(n)$. The automorphisms of the projection
$V_n\to V$ are given by $\GL_2(\Z/n\Z)/\{\pm 1\}$, so that the group
of deck transformations of the tower $\cV$ is in this case of the form
$$ \Aut_V (\cV)= \varprojlim_n  \GL_2(\Z/n\Z)/\{\pm 1\} =
\GL_2(\hat\Z)/\{\pm 1\}. $$
The inverse limit $\varprojlim \Gamma\backslash \H$ over
congruence subgroups $\Gamma\subset \SL(2,\Z)$ gives a connected
component of \eqref{ShiMod}, while by taking congruence subgroups in
$\SL(2,\Q)$ one obtains the ad\`elic version $Sh(\GL_2,\H^\pm)$.
The simple reason why it is
necessary to pass to the nonconnected case is the following. The
varieties in the tower are arithmetic varieties defined over
number fields. However, the number field typically changes along
the levels of the tower ($V_n$ is defined over the cyclotomic
field $\Q(\zeta_n)$). Passing to nonconnected Shimura varieties
allows for the definition of a canonical model where the
whole tower is defined over the same number field. 

The quotient \eqref{ShiMod} parameterizes invertible 2-dimensional
$\Q$-lattices up to scaling (\cf \cite{Mi}). When, instead of
restricting to the invertible case, we consider commensurability 
classes of 2-dimensional $\Q$-lattices up to scaling, we obtain a
noncommutative space whose classical points are the Shimura variety
\eqref{ShiMod}. More precisely, we have the following (\cf \cite{CMR}):

\begin{lem}\label{2dimQlatSh}
The space of commensurability classes of 2-dimensional
$\Q$-lattices up to scaling is described by the quotient
\begin{equation}\label{ShGL2NC}
Sh^{(nc)}(\GL_2,\H^{\pm}):=\GL_2(\Q) \backslash (M_2(\A_f)\times
\H^{\pm}).
\end{equation}
\end{lem}

Any 2-dimensional $\Q$-lattice can be written in the form
$$ (\Lambda,\phi)=(\lambda (\Z+\Z\tau),\lambda\rho), $$
for some $\lambda\in \C^*$, some $\tau\in \H$, and some $$\rho\in
M_2(\hat\Z)=\Hom(\Q^2/\Z^2,\Q^2/\Z^2).$$ Thus, the space of
2-dimensional $\Q$-lattices up to the scale factor $\lambda\in
\C^*$ and up to isomorphisms, is given by
\begin{equation}\label{2dQlat}
 M_2(\hat\Z)\times \H  \mod \Gamma=\SL(2,\Z).
\end{equation}
The commensurability relation is implemented by the partially
defined action of $\GL_2^+(\Q)$ on this space, given by
$$ g (\rho,z)=(g\rho, g(z)), $$
where $g(z)$ denotes action as fractional linear transformation.

Equivalently, the data $(\Lambda,\phi)$ of a $\Q$-lattice in $\C$ are
equivalent to data $(E,\eta)$ of an elliptic curve $E=\C/\Lambda$
and an $\A_f$-homomorphism
\begin{equation}\label{Afhom}
\eta : \Q^2 \otimes \A_f \to \Lambda \otimes \A_f,
\end{equation}
where $\Lambda\otimes \A_f = (\Lambda \otimes \hat\Z)\otimes \Q$,
with
\begin{equation}\label{torEhatZ}
\Lambda \otimes \hat\Z = \varprojlim_n \Lambda/n\Lambda,
\end{equation}
and $\Lambda/n\Lambda = E[n]$ the $n$th torsion of $E$. Since the
$\Q$-lattices need not be invertible, we do not require that $\eta$
be an $\A_f$-isomorphism (\cf \cite{Mi}).

The commensurability relation between $\Q$-lattices
corresponds to the equivalence $(E,\eta)\sim (E',\eta')$ given by an
isogeny $g: E \to E'$ and $\eta'=(g\otimes 1)\circ\eta$. Namely, the
equivalence classes can be identified with the quotient of
$M_2(\A_f)\times \H^{\pm}$ by the action of $\GL_2(\Q)$,
$(\rho,z)\mapsto (g\rho, g(z))$.

To associate a quantum statistical mechanical system to the
space of commensurability classes of 2-dimensional $\Q$-lattices up to
scaling, it is convenient to use the description \eqref{2dQlat}.
Then one can consider as noncommutative algebra of coordinates the
convolution algebra of continuous compactly supported
functions on the quotient of the space
\begin{equation}\label{Uspace}
 \cU:=\{ (g,\rho,z) \in \GL_2^+(\Q)\times M_2(\hat\Z)\times \H |
g\rho\in M_2(\hat\Z) \}
\end{equation}
by the action of $\Gamma\times \Gamma$, 
\begin{equation}\label{Gammaacts}
 (g,\rho,
z)\mapsto(\gamma_1g\gamma_2^{-1},\gamma_2 z). 
\end{equation}
One endows this algebra with the convolution product
\begin{equation}\label{Heckeprod2}
(f_1*f_2)(g,\rho,z)= \displaystyle{\sum_{s\in \Gamma\backslash
\GL_2^+(\Q): s\rho\in M_2(\hat\Z)}} f_1(gs^{-1},s\rho,s(z))
f_2(s,\rho,z)
\end{equation}
and the involution
$f^*(g,\rho,z)=\overline{f(g^{-1},g\rho,g(z))}$.

The time evolution is given by
\begin{equation}\label{evolution2}
\sigma_t(f) (g,\rho,z) = \det(g)^{it} \, f(g,\rho,z).
\end{equation}

For $\rho\in M_2(\hat\Z)$ let
\begin{equation}\label{Grho}
 G_\rho:= \{ g \in \GL_2^+(\Q)\, :\, g\rho\in M_2(\hat\Z) \}
\end{equation}
and consider the Hilbert space $\cH_\rho = \ell^2(\Gamma\backslash
G_\rho)$.

A 2-dimensional $\Q$-lattice $L=(\Lambda,\phi)=(\rho,z)$
determines a representation of the Hecke algebra by bounded
operators on $\cH_\rho$, setting
\begin{equation}\label{piLrep}
 (\pi_L(f)\xi)\, (g) =\displaystyle{\sum_{s\in \Gamma\backslash
G_\rho}} \, f(gs^{-1},s\rho,s(z))\, \xi(s).
\end{equation}
In particular, when the $\Q$-lattice $L=(\Lambda,\phi)$ is
invertible one obtains $$\cH_\rho\cong \ell^2(\Gamma\backslash
M_2^+(\Z)).$$ In this case, the Hamiltonian implementing the time
evolution \eqref{evolution2} is given by the operator
\begin{equation}\label{Hamiltonian2}
 H \, \epsilon_m = \log \det (m) \,\, \epsilon_m.
\end{equation}
Thus, in the special case of invertible $\Q$-lattices
\eqref{piLrep} yields a {\em positive energy} representation.
In general for $\Q$-lattices which are not commensurable to
an invertible one, the corresponding Hamiltonian $H $
is not bounded below.

The Hecke algebra \eqref{Heckeprod2} admits a $C^*$-algebra
completion $\cA_2$, where the norm is the sup over all
representations $\pi_L$.

The partition function for this $\GL_2$ system is given by
\begin{equation}\label{partition2}
Z(\beta)= \sum_{m\in \Gamma\backslash M_2^+(\Z)}\det(m)^{-\beta}
=\sum_{k=1}^\infty \sigma(k)\,
k^{-\beta}=\zeta(\beta)\zeta(\beta-1),
\end{equation}
where $\sigma(k)=\sum_{d|k} d$. This suggests the fact that
one expects two phase transitions to take place, at $\beta=1$ and
$\beta=2$, respectively. 

The set of components of $Sh(\GL_2,\H^\pm)$ is given by
\begin{equation}\label{components}
\pi_0(Sh(\GL_2,\H^\pm))=Sh(\GL_1,\{\pm 1\}).
\end{equation}
Similarly, on the level of the corresponding noncommutative spaces
\eqref{ShGL2NC} and \eqref{Sh1nc}, the 
identification \eqref{components} gives the compatibility
between the $\GL_1$ and the $\GL_2$ system (\cf \cite{CMR}). 
At the level of the
classical commutative spaces, this is given by the map
\begin{equation}\label{detsign}
 \det \times\, {\rm sign} : Sh(\GL_2,\H^{\pm}) \to Sh(\GL_1,\{\pm
1\}),
\end{equation}
which corresponds in fact to passing to the set $\pi_0$ of
connected components.

\subsection{KMS states of the $\GL_2$-system}

The main result of \cite{CoMar} on the structure of KMS states is the
following. 

\begin{thm}\label{GL2KMS}
The KMS$_\beta$ states of the $\GL_2$-system have the following
properties:
\begin{enumerate}
\item In the range $\beta\leq 1$ there are no KMS states.
\item In the range $\beta>2$ the set of extremal KMS states is
given by the classical Shimura variety
\begin{equation}\label{Ekmsbeta2}
\cE_\beta \cong \GL_2(\Q)\backslash \GL_2(\A) /\C^*.
\end{equation}
\end{enumerate}
\end{thm}

This shows that the extremal KMS states at sufficiently low
temperature are parameterized by the invertible $\Q$-lattices. The
explicit expression for these extremal KMS$_\beta$ states is
obtained as
\begin{equation}\label{KMSexpr}
 \varphi_{\beta,L}(f)=\frac{1}{Z(\beta)} \sum_{m\in
\Gamma\backslash M_2^+(\Z)} f(1,m\rho,m(z))\, \det(m)^{-\beta}
\end{equation}
where $L=(\rho,z)$ is an invertible $\Q$-lattice.

As the temperature rises, and we let $\beta\to 2$ from above,
all the different phases of the system merge, which is strong
evidence for the fact that 
in the intermediate range $1< \beta \leq 2$
the system should have only a single KMS state.

The symmetry group of $\cA_2$ (including both automorphisms and
endomorphisms) can be identified with the group
\begin{equation}\label{SymmGL2}
\GL_2(\A_f)=\GL_2^+(\Q) \GL_2(\hat\Z).
\end{equation}
Here the group $\GL_2(\hat\Z)$ acts by automorphisms,
\begin{equation}\label{GLhatZ}
 \theta_\gamma(f) (g,\rho,z)= f(g,\rho \gamma, z) ,
\end{equation}
as the group of deck transformations of
coverings of modular curves. 

The novelty with respect to the BC case is that we also have an action
of $\GL_2^+(\Q)$ by endomorphisms
\begin{equation}\label{endoaction}
 \theta_m(f) (g,\rho,z) = \left\{ \begin{array}{lr} f(g,\rho
\tilde
m^{-1}, z) & \rho\in m\,M_2(\hat\Z) \\[2mm]
0 & \text{ otherwise } \end{array} \right.
\end{equation}
where $\tilde m= \det(m) m^{-1}$.

The subgroup $\Q^*\hookrightarrow \GL_2(\A_f)$ acts by inner,
hence the group of symmetries of the set of extremal states
$\cE_\beta$ is of the form
\begin{equation}\label{Ssymmetries}
S=\Q^*\backslash \GL_2(\A_f).
\end{equation}

In the case of $\cE_\infty$ states (defined as weak limits) the
action of $\GL_2^+(\Q)$ is more subtle to define. In fact,
\eqref{endoaction} does not directly induce a nontrivial action on
$\cE_\infty$. However, there is a nontrivial action induced by the
action on $\cE_\beta$ states for sufficiently large $\beta$. The
action on the KMS$_\infty$ states is obtained by a ``warming up
and cooling down procedure'', as in \eqref{warming} and
\eqref{warmcool}.

Finally, the Galois action on the extremal KMS$_\infty$ states is
described by the following result \cite{CoMar}.

\begin{thm}\label{GalGL2infty}
There exists an arithmetic subalgebra $\cA_{2,\Q}$ of unbounded
multiplier of the $C^*$-algebra $\cA_2$, such that the following
holds. For $\varphi_{\infty,L}\in \cE_\infty$ with $L=(\rho,\tau)$
generic, the values on arithmetic elements satisfy
\begin{equation}\label{values}
\varphi_{\infty,L}(\cA_{2,\Q})\subset F_\tau,
\end{equation}
where $F_\tau$ is the embedding of the modular field in $\C$
given by evaluation at $\tau\in \H$. There is an isomorphism
\begin{equation}\label{thetaphi1}
\theta_\varphi : \Gal(F_\tau/\Q)
\stackrel{\simeq}{\longrightarrow} \Q^* \backslash \GL_2(\A_f),
\end{equation}
that intertwines the Galois
action on the values of the state with the action of symmetries,
\begin{equation}\label{intertwineGL2}
\gamma^{-1}\, \varphi(f) = \varphi( \theta_\varphi(\gamma) f), \ \ \ \
\forall f\in \cA_{2,\Q}, \ \ \forall\gamma\in \Gal(F_\tau/\Q).
\end{equation}
\end{thm}

Here recall that the modular field $F$ is the field of modular
functions over $\Q^{ab}$, namely the union of the fields $F_N$ of modular
functions of level $N$ rational over the cyclotomic field
$\Q(\zeta_n)$, that is, such that the $q$-expansion in powers of
$q^{1/N}=\exp(2\pi i \tau/N)$ has all coefficients in $\Q(e^{2\pi
i/N})$.

It has explicit generators given by the
Fricke functions (\cite{Sh}, \cite{Lang}). If $\wp$ is the
Weierstrass $\wp$-function, which gives the parameterization $$z
\mapsto (1,\wp(z;\tau,1),\wp'(z;\tau,1))$$ of the elliptic curve
$$y^2 = 4x^3 - g_2(\tau) x - g_3(\tau)$$ by the quotient
$\C/(\Z+\Z \tau)$, then the Fricke functions are homogeneous
functions of 1-dimensional lattices of weight zero, parameterized
by $v\in \Q^2/\Z^2$, of the form
\begin{equation}\label{Fricke}
 f_v(\tau)=- 2^7 3^5 \frac{g_2(\tau)g_3(\tau)}{\Delta(\tau)}\,\, \wp
(\lambda_\tau(v); \tau, 1),
\end{equation}
where $\Delta=g_2^3 - 27 g_3^2$ is the discriminant and
$\lambda_\tau(v):=v_1 \tau + v_2$.

For generic $\tau$ (such that $j(\tau)$ is transcendental), evaluation
of the modular functions at the point $\tau\in \H$ gives an embedding
$F_\tau \hookrightarrow \C$. There is a corresponding isomorphism
\begin{equation}\label{GalFtau}
 \theta_\tau: \Gal(F_\tau/\Q)
\stackrel{\simeq}{\to} \Aut(F).
\end{equation}
The isomorphism \eqref{thetaphi1} is given by
\begin{equation}\label{thetaphi2}
 \theta_\varphi(\gamma)=\rho^{-1}\, \theta_\tau(\gamma) \, \rho,
\end{equation}
for $\theta_\tau$ as in \eqref{GalFtau}.
In fact, the group of automorphisms $\Aut(F)$ has a completely explicit
description, due to a result of Shimura \cite{Sh}, which identifies 
it with the quotient
$$ \Aut(F)\cong \Q^* \backslash \GL_2(\A_f). $$

We still need to explain what is the arithmetic subalgebra that
appears in Theorem \ref{GalGL2infty}.

We consider continuous functions on the quotient of
\eqref{Uspace} by the action \eqref{Gammaacts}, which have finite
support in the variable $g\in \Gamma\backslash\GL_2^+(\Q)$. For
convenience we adopt the notation
$$ f_{(g,\rho)}(z) = f(g,\rho,z) $$
so that $f_{(g,\rho)}\in C(\H)$.  Let $p_N: M_2(\hat\Z)\to
M_2(\Z/N\Z)$ be the canonical projection. We say that $f$ is of
level $N$ if
$$ f_{(g,\rho)} = f_{(g,p_N(\rho))} \ \ \ \ \forall (g,\rho). $$
Then $f$ is completely determined by the functions
$$ f_{(g,m)} \in C(\H), \ \ \ \ \text{ for } m\in M_2(\Z/N\Z). $$
Notice that the invariance
$$ f(g\gamma,\rho,z)=f(g,\gamma\rho,\gamma(z)), $$ for all
$\gamma \in \Gamma$ and for all $(g,\rho,z)\in \cU$, implies that
we have
\begin{equation}\label{modular}
f_{(g,m)| \gamma} = f_{(g,m)}, \ \ \ \ \forall \gamma \in
\Gamma(N)\cap g^{-1}\Gamma g.
\end{equation}
so that $f$ is invariant under a congruence subgroup.

Elements $f$ of $\cA_{2,\Q}$ are characterized by the
following properties.
\begin{itemize}
\item The support of $f$ in $\Gamma\backslash\GL_2^+(\Q)$ is
finite.
\item The function $f$ is of finite level with
$$ f_{(g,m)} \in F \ \ \ \ \forall (g,m). $$
\item The function $f$ satisfies the {\em cyclotomic condition}:
$$ f_{(g,\alpha(u)m)} = {\rm cycl}(u) \, f_{(g,m)}, $$
for all $g\in \GL_2^+(\Q)$ diagonal and all $u\in \hat\Z^*$, with
$$ \alpha(u)=\begin{pmatrix} u& 0 \\ 0 & 1 \end{pmatrix}. $$
\end{itemize}
Here ${\rm cycl}: \hat\Z^* \to \Aut(F)$
denotes the action of the Galois group $\hat\Z^*
\simeq \Gal(\Q^{ab}/\Q)$ 
on the coefficients of the $q$-expansion in powers of
$q^{1/N}=\exp(2\pi i \tau/N)$.

If we took only the first two conditions, this would allow the
algebra $\cA_{2,\Q}$ to contain the cyclotomic field
$\Q^{ab}\subset \C$, but this would prevent the existence of
states satisfying the desired Galois property. In fact,
if the subalgebra contains scalar elements in $\Q^{cycl}$,
the sought for Galois property 
would not be compatible with the $\C$-linearity of states. The
cyclotomic condition then 
forces the spectrum of the elements of $\cA_{2,\Q}$ to
contain all Galois conjugates of any root of unity that appears in the
coefficients of the $q$-series, so that 
elements of $\cA_{2,\Q}$ cannot be scalars. 

The algebra $\cA_{2,\Q}$ defined by the properties above is a
subalgebra of unbounded multipliers of $\cA_2$, which is globally
invariant under the group of symmetries $\Q^*\backslash \GL_2(\A_f)$.

\section{Quadratic fields}

The first essential step in the direction of generalizing the
BC system to other number fields and exploring a possible approach to
the explicit class field theory problem via noncommutative geometry is
the construction of a system that recovers the explicit class field
theory for imaginary quadratic fields $\K=\Q(\sqrt{-d})$, for $d$ a
positive integer. This case is being investigated as part of ongoing
joint work of Alain Connes, Niranjan Ramachandran, and the author
\cite{CMR}. 

A system for imaginary quadratic fields can be constructed by
considering 1-dimensional $\K$-lattices. These are lattices in $\K$
with a homomorphism $\phi: \K/\O  \to \Q\Lambda/\Lambda$, where
$\O$ is the ring of integers of $\K$. A choice of
a generator $\tau\in \H$ such that $\K=\Q(\tau)$ realizes
1-dimensional $\K$-lattices as a particular set of 2-dimensional
$\Q$-lattices. The commensurability relation for 1-dimensional
$\K$-lattices can be defined as a commensurability of the underlying
2-dimensional $\Q$-lattices where the isomorphism is realized by
multiplication by an element in $\K^*$, viewed as embedded in
$\GL_2^+(\Q)$. 

The set of commensurability classes of 1-dimensional $\K$-lattices up
to scale is then described by the quotient
\begin{equation}\label{1Klat}
\A_{\K,f}/\K^*,
\end{equation}
where $\A_{\K,f}=\A_f\otimes \K$ are the finite ad\`eles of $\K$,
while the set of invertible 1-dimensional $\K$-lattices up
to scale can be identified with the id\`ele class group 
\begin{equation}\label{CKDK2}
C_\K/D_\K = \A_{\K,f}^*/\K^*.
\end{equation}

The noncommutative algebra of coordinates of the quotient
\eqref{1Klat} is obtained as in the case of the $\GL_2$-system as the
convolution algebra of functions on 
$$ \cU_\K =\{ (x,\rho)\in \K^* \times \hat\O  :\, x\rho\in
\hat\O \}, $$
with $\hat\O  =\hat\Z \otimes \O $, endowed 
with the time evolution 
\begin{equation}\label{timeK}
\sigma_t(f)(x,\rho,\lambda)= N(x)^{-it} f (x,\rho,\lambda),
\end{equation}
where $N: K^* \to \Q^*$ is the norm map.

This quantum statistical mechanical system has properties that are, in
a sense, intermediate between the BC system and the $\GL_2$-system. 

The symmetry group is the group of id\`eles classes
$\A_{\K,f}^*$, with the subgroup $\K^*$ acting by inner, so that
the induced action on KMS states is by the id\`ele class group
\eqref{CKDK2}. As in the case of the $\GL_2$-system, only the subgroup
$\hat\O /\O ^*$, with $\O ^*$ the group of units, acts by
automorphisms, while the full action of 
$\A_{\K,f}^*/\K^*$ also involves endomorphisms. 

This, in particular, shows the appearance of the class number of
$\K$, as one can see from the commutative diagram
\begin{eqnarray}
\diagram 1\rto & \hat\cO^* /\cO^*  \rto\dto^{\simeq} &   \A_{\K,f}^*/\K^*
\dto^{\simeq} \rto & {\rm Cl}(\cO )\dto^{\simeq} \rto & 1 \\
1\rto & \Gal(\K^{ab}/H) \rto & \Gal(\K^{ab}/\K) \rto & \Gal(H/\K) \rto & 1,
\enddiagram
\label{Hfielddiagr}
\end{eqnarray}
where $H$ is the Hilbert class field of $\K$, \ie, its maximal abelian
unramified extension. The ideal class group ${\rm Cl}(\cO )$ is naturally
isomorphic to the Galois group $\Gal(H/\K)$.
The case of class number one is analogous to the BC system, as
was already observed (\cf \cite{HaLe}). 

The same definition of the arithmetic algebra $\cA_{2,\Q}$ for the
$\GL_2$-system provides a subalgebra of the system for the imaginary
quadratic field $\K$. In fact, one considers 
functions in $\cA_{2,\Q}$ whose
finite support in the $g$ variable lies in $\K^*$ embedded in
$\GL_2^+(\Q)$, and then takes the restriction to the set of 1-dimensional
$\K$-lattices, by evaluation at $\tau\in \H$ and restriction to
$\rho\in \hat\O $, embedded in $M_2(\hat\Z)$.  

Notice that, because of the fact that the variable $z\in \H$ is now
set to be $z=\tau$, we obtain a subalgebra as in the BC case
and not an algebra of unbounded multipliers as in the $\GL_2$-case. 

With this choice of the arithmetic subalgebra, one obtains a result
analogous to Theorem \ref{GalGL2infty}. Namely, since now the
2-dimensional $\Q$-lattices $L=(\rho,\tau)$ we are considering have
a fixed $\tau\in \H\cap \K$, we are no longer in the generic case of 
Theorem \ref{GalGL2infty}. For a CM point $\tau\in \H\cap \K$
evaluation of elements $f$ in the modular field at $\tau$ no longer
gives an embedding. The image $F_\tau \subset \C$ is in this case a
copy of the maximal abelian extension of $\K$ (\cf \cite{Sh})
\begin{equation}\label{FtauCM}
F \to F_\tau \simeq \K^{ab}\subset \C.
\end{equation}
The explicit action of the Galois group $\Gal(\K^{ab}/\K)$ is obtained
through the action of automorphisms of the modular field via Shimura
reciprocity \cite{Sh}, as described in the following diagram with
exact rows: 
\begin{eqnarray}
\diagram 1\rto & K^* \rto & \GL_1(\A_{K,f})
\dto^{q_\tau} \rto^{\theta} &
\Gal(K^{ab}/K) \rto & 1 \\
1\rto & \Q^* \rto & \GL_2(\A_f) \rto^{\sigma} & \Aut(F) \rto & 1.
\enddiagram
\label{Shirecdiagr}
\end{eqnarray}
Here $\theta$ is the class field theory isomorphism and $q_\tau$ is
the embedding determined by the choice of $\tau\in H$ with
$\K=\Q(\tau)$. Thus, the explicit Galois action is given by
$$ \theta(\gamma)( f(\tau) ) = f^{\sigma(q_\tau(\gamma))}(\tau), $$
where $f\mapsto f^\alpha$ denotes the action of $\alpha\in \Aut(F)$.

This provides the intertwining between the symmetries of the
quantum statistical mechanical system and the Galois action on the
image of states on the elements of the arithmetic subalgebra.

The structure of KMS states at positive temperatures is similar to the
Bost--Connes case (\cf \cite{CMR}).

The next fundamental question in the direction of generalizations
of the BC system to other number fields is how to approach the
more case of real quadratic fields. We give a brief outline of
Manin's ideas on real multiplication \cite{Man3} \cite{Man5}, and
suggest how they may combine with the $\GL_2$-system described above.

\subsection{Real multiplication}

Recently, Manin developed a theory of real multiplication for
noncommutative tori \cite{Man5},
aimed at providing a setting, within noncommutative
geometry, where to treat the problem of abelian extensions of real
quadratic fields on a somewhat similar footing as the known case of
imaginary quadratic fields, for which the theory of elliptic curves
with complex multiplication provides the right geometric setup.

The first entry in the dictionary developed in \cite{Man5} between 
elliptic curves with complex multiplication and noncommutative tori
with real multiplication consists of a parallel between lattices 
and pseudolattices in $\C$, 
$$ (\C, \text{lattice}\,\, \Lambda) \ \ \leftrightsquigarrow \ \  (\R,
\text{pseudo-lattice}\,\, L). $$ 
By a pseudo-lattice one means the data of a free abelian group of rank
two, with an injective homomorphism to a 1-dimensional complex vector
space, such that the image lies in an oriented real line.

This aims at generalizing the 
well known equivalence between the category
of elliptic curves and the category of lattices, realized by the
period functor, to a setting that includes noncommutative tori. 

As we have seen previously, any 2-dimensional lattice is, up to
isomorphism, of the form 
$\Lambda_\tau = \Z + \Z\tau$, for $\tau \in \C\smallsetminus \R$ and
non-trivial morphisms between such lattices are given by the action of
matrices $M_2(\Z)$ by fractional linear transformations. Thus, the moduli
space of lattices up to isomorphism is given by the quotient of
$\P^1(\C)\smallsetminus \P^1(\R)$ 
by $\PGL_2(\Z)$. 

A pseudo-lattice is, up to
isomorphism, of the form $L_\theta = \Z + \Z\theta$, for $\theta \in
\R\smallsetminus \Q$, and non-trivial morphisms of pseudo-lattices are
again given by matrices in $M_2(\Z)$ acting by fractional linear
transformations on $\P^1(\R)$. The moduli space of
pseudo-lattices is given by the quotient of $\P^1(\R)$ by the action
of $\PGL(2,\Z)$. Since the action does not give rise to a nice
classical quotient, this moduli space should be treated as a
noncommutative space. 

As described at lenght in the previous chapter, the resulting space
represents a component in the ``boundary'' of 
the classical moduli space of elliptic curves, which parameterizes
those degenerations from lattices to pseudo-lattices, that are
invisible to the usual algebro-geometric setting. 

The cusp, \ie the
orbit of $\P^1(\Q)$, corresponds to the degenerate case where the
image in $\C$ of the rank two free abelian group has rank one. 

This gives another entry in the dictionary, regarding the 
moduli spaces:
$$ \H/\PSL(2,\Z)  \ \ \leftrightsquigarrow \ \   C(\P^1(\R))\rtimes
\PGL(2,\Z). $$    

In the correspondence of pseudo-lattices and noncommutative tori, the
group of invertible morphisms of pseudo-lattices corresponds to
isomorphisms of noncommutative tori realized by strong Morita
equivalences. In this context a ``morphism'' is not
given as a morphism of algebras but as a map of the category of
modules, obtained by tensoring with a bimodule. 
The notion of Morita equivalences as morphisms fits into
the more general context of correspondences for operator algebras as
in \cite{Co94} \S V.B, as well as in the algebraic approach to
noncommutative spaces of \cite{Rosenberg98}. 

The category of noncommutative tori is defined by considering as
morphisms the isomorphism classes of ($C^*$-)bimodules that are range
of projections. The functor from noncommutative tori to
pseudo-lattices (\cf \cite{Man3} \S 3.3 and \cite{Man5} \S
1.4) is then given on objects by 
\begin{equation}\label{functorNCT-PL}
\T^2_\theta \mapsto \left( K_0(\cA_\theta), HC_0(\cA_\theta) 
, \tau: K_0(\cA_\theta) \to
HC_0(\cA_\theta) \right). 
\end{equation}
Here $\cA_\theta$ is the algebra of the noncommutative torus,
$HC_0(\cA_\theta)= \cA_\theta/[\cA_\theta,\cA_\theta]$, with $\tau$
the universal trace, and the orientation is determined by the cone of
positive  elements in $K_0$. On morphisms the functor is given by
\begin{equation}\label{functorNCT-PL2}
\cM_{\theta,\theta'} \mapsto \left( [ \cE ]
\mapsto [ \cE \otimes_{\cA_\theta} \cM_{\theta,\theta'} ] \right), 
\end{equation}
where $\cM_{\theta,\theta'}$ are the bimodules constructed by Connes
in \cite{ConnesCR}. 
A crucial point in this definition is the fact that, for noncommutative
tori, finite projective modules are classified by the value of a unique
normalized trace (\cf \cite{ConnesCR} \cite{Rief1}). 

The functor of \eqref{functorNCT-PL} \eqref{functorNCT-PL2} is weaker
than an equivalence of categories. For instance, it maps trivially all
ring homomorphisms that act trivially on $K_0$. 
However, this correspondence is sufficient to develop a theory of
real multiplication for noncommutative tori parallel to the theory of
complex multiplication for elliptic curves (\cf \cite{Man3}
\cite{Man5}).  

For lattices/elliptic curves, the typical situation is that
$\End(\Lambda)=\Z$, but there are exceptional lattices for which
$\End(\Lambda)\supsetneq \Z$. In this case, there exists a complex
quadratic field $\K$ such that $\Lambda$ is isomorphic to a lattice in
$\K$. More precisely, the endomorphism group is given by
$\End(\Lambda) = \Z + f \O $, where $\O $ is the ring of integers of
$\K$ and the integer $f\geq 1$ is called the conductor. Such lattices
are said to have complex multiplication. The elliptic curve
$E_\K$ with $E_\K(\C)=\C/ \O $ is endowed with a complex
multiplication map, given on the universal cover by
$x\mapsto ax$, $a\in \O $. 

Similarly, there is a parallel situation for pseudo-lattices:
$\End(L)\supsetneq \Z$ happens when there exists a real quadratic
field $\K$ such that $L$ is isomorphic to a pseudolattice contained in
$\K$. In this case, one also has $\End(L)=\Z + f \O $. 
Such pseudo-lattices are said to have {\em real multiplication}. 
The pseudo-lattices $L_\theta$ that have real multiplication
correspond to values of $\theta \in \R \smallsetminus \Q$ that are
quadratic irrationalities. These are characterized by the 
having eventually periodic continued fraction expansion. The ``real
multiplication map'' is given by tensoring with a bimodule: in the
case of $\theta$ with periodic continued fraction expansion, there is
an element $g\in \PGL(2,\Z)$ such that $g\theta =\theta$, to which we
can associate an $\cA_\theta$--$\cA_\theta$ bimodule $\cE$. 

An analog of isogenies for noncommutative tori is obtained by
considering morphisms of pseudolattices $L_\theta \mapsto
L_{n\theta}$ which correspond to a Morita morphism given by
$\cA_\theta$ viewed as an $\cA_{n\theta}$--$\cA_\theta$ bimodule,
where $\cA_{n\theta} \hookrightarrow \cA_\theta$ by $U \mapsto U^n$,
$V\mapsto V$.

By considering isogenies, one can enrich the dictionary between moduli
space of elliptic 
curves and moduli space of Morita equivalent noncommutative tori. This
leads one to consider the whole tower of modular curves parameterizing
elliptic curves with level structure and the corresponding tower of
noncommutative modular curves described in the previous chapter (\cf
\cite{ManMar}), 
\begin{equation}\label{NC-modcurv}
 \H/G \ \ \leftrightsquigarrow \ \ C(\P^1(\R)\times \P)\rtimes
\Gamma, 
\end{equation}
for $\Gamma =\PGL(2,\Z)$ (or $\PSL(2,\Z)$) and $G\subset \Gamma$ a
finite index subgroup, with $\P=G\backslash \Gamma$ the coset
space. 
 
In the problem of constructing the maximal abelian extension of
complex quadratic fields, one method is based on evaluating at the
torsion points of the elliptic curve $E_\K$ a power of the Weierstrass
function. This means considering the corresponding values of a
coordinate on the projective line $E_\K / \O ^*$. The analogous
object in the noncommutative setting, replacing this projective line,
should be (\cf \cite{Man5}) a crossed product of functions on 
$\K$ by the $ax +b$ group with $a\in \O ^*$ and $b\in \O $, for a
real quadratic field $\K$.
 
\smallskip

A different method to construct abelian extensions of a complex
quadratic field $\K$ is via Stark numbers. Following the notation of
\cite{Man5}, if $j: \Lambda \to V$
is an injective homomorphism of a free abelian group of rank two to a
one-dimensional complex vector space, and $\lambda_0 \in
\Lambda\otimes \Q$, then one can consider the zeta function
\begin{equation}\label{zeta-im}
 \zeta(\Lambda,\lambda_0,s)= \sum_{\lambda \in \Lambda}
\frac{1}{|j(\lambda_0 + \lambda)|^{2s} }. 
\end{equation}
Stark proved (\cf \cite{Stark}) that the numbers
\begin{equation}\label{St-im}
 S(\Lambda,\lambda_0) = \exp \, \zeta^\prime (\Lambda,\lambda_0,0) 
\end{equation}
are algebraic units generating certain abelian extensions of $\K$. The
argument in this case is based on a direct computational tool (the
Kronecker second limit formula) and upon reducing the problem to the
theory of complex multiplication. There is no known independent
argument for the Stark conjectures, while the analogous question is
open for the case of real quadratic fields. 

For a real quadratic field $\K$, instead of zeta
functions of the form \eqref{zeta-im}, the conjectural Stark units are
obtained from zeta functions of the form (in the notation of
\cite{Man5}) 
\begin{equation}\label{zeta-re}
\zeta(L,l_0,s)= {\rm sgn}\, l_0' N(\ma)^s \sum \frac{{\rm sgn}
(l_0+l)'}{|N(l_0+l)|^s}, 
\end{equation}
where $l\mapsto l'$ is the action of the nontrivial element in
$Gal(\K/\Q)$ and $N(l)=ll'$. The element $\l_0\in O $ is chosen so
that the ideal $\ma =(L,l_0)$ and $(l_0)\ma^{-1}$ are coprime wih $L
\ma^{-1}$, and the summation over $l\in L$ in \eqref{zeta-re} is
restricted by taking only one representative from each coset class
$(l_0+l)\epsilon$, for units $\epsilon$ satisfying $(l_0+L)\epsilon =
(l_0+L)$, \ie $\epsilon \equiv 1 \mod L \ma^{-1}$. Then the Stark
numbers are given as in \eqref{St-im} by (\cf \cite{Man5}) 
\begin{equation}\label{St-re}
 S(L,l_0) = \exp \, \zeta^\prime (L,l_0,0).  
\end{equation}

In \cite{Man5}, Manin develops an approach to the computation of
sums \eqref{zeta-re} based on a version of theta functions for
pseudo-lattices, which are obtained by averaging theta constants of
complex lattices along geodesics with ends at a pair of conjugate
quadratic irrationalities $\theta,\theta'$ in $\R\smallsetminus \Q$. 

This procedure fits into a general philosophy, according to which one
can recast part of the arithmetic theory of modular curves in terms of
the noncommutative boundary \eqref{NC-modcurv} by studying the
limiting behavior when $\tau \to \theta \in \R\setminus \Q$ along
geodesics, or some averaging along such geodesics. In general, some
nontrivial result is obtained when approaching $\theta$ along a path
that corresponds to a geodesic in the modular curve that spans a
limiting cycle, which is the case precisely when the endpoint $\theta$
is a quadratic irrationality. An example of this type of behavior is
the theory of ``limiting modular symbols'' developed in \cite{ManMar},
\cite{Mar-lyap}.

\subsection{Pseudolattices and the $\GL_2$-system}

The noncommutative space \eqref{ShGL2NC} of the $\GL_2$-system
also admits a compactification, now given by adding the
boundary $\P^1(\R)$ to $\H^\pm$, as in the noncommutative
compactification of modular curves of \cite{ManMar},
\begin{equation}\label{compSh2nc}
\begin{array}{rl}
\overline{Sh^{(nc)}}(\GL_2,\H^{\pm}):= & \GL_2(\Q) \backslash
(M_2(\A_f)\times \P^1(\C))\\[2mm] = &\GL_2(\Q) \backslash
M_2(\A)/\C^*,
\end{array}
\end{equation}
where $\P^1(\C)=\H^\pm \cup \P^1(\R)$. 

This corresponds to adding
to the space of commensurability classes of 2-dimensional
$\Q$-lattices the pseudolattices in the sense of \cite{Man5},
here considered together with a $\Q$-structure.
It seems then that Manin's real multiplication program may fit in with
the boundary of the noncommutative space of the $\GL_2$-system.
The crucial question in this respect becomes the construction of
an arithmetic algebra associated to the noncommutative modular
curves. The results illustrated in the previous chapter, regarding
identities involving modular forms at the boundary of the classical
modular curves and limiting modular symbols, as well as the still
mysterious phenomenon of ``quantum modular forms'' identified by
Zagier, point to the fact that there should exist a rich class
of objects replacing modular forms, when ``pushed to the boundary''.

Regarding the role of modular forms, notice that, in the case of
\eqref{compSh2nc}, we can again consider the dual system. This is a
$\C^*$-bundle
\begin{equation}\label{L2}
\cL_2= \GL_2(\Q) \backslash M_2(\A).
\end{equation}
On this dual space modular forms appear naturally instead of
modular functions and the algebra of coordinates contains the
modular Hecke algebra of Connes--Moscovici (\cite{CoMo1},
\cite{CoMo2}) as arithmetic subalgebra.

Thus, the noncommutative geometry of the space of $\Q$-lattices modulo
the equivalence relation of commensurability provides a setting
that unifies several phenomena involving the interaction of
noncommutative geometry and number theory. These include
the Bost--Connes system, the noncommutative space underlying the
construction of the spectral realization of the zeros of the
Riemann zeta function in \cite{Connes-Zeta}, the modular Hecke
algebras of \cite{CoMo1} \cite{CoMo2}, and the noncommutative
modular curves of \cite{ManMar}.

\chapter{Noncommutative geometry at arithmetic infinity}

This chapter is based on joint work of Katia Consani and the
author (\cite{CM}, \cite{CM1}, \cite{CM2}, \cite{CM3},
\cite{CMrev}) that proposes a model for the dual graph of the
maximally degenerate fibers at the archimedean places of an
arithmetic surface in terms of a noncommutative space (spectral
triple) related to the action of a Schottky group on its limit
set. This description of $\infty$-adic geometry provides a
compatible setting that combines Manin's result \cite{Man-hyp} on
the Arakelov Green function for arithmetic surfaces in terms of
hyperbolic geometry and for a cohomological construction of
Consani \cite{KC} associated to the archimedean fibers of
arithmetic varieties, related to Deninger's calculation of the
local $L$-factors as regularized determinants (\cf \cite{Den1},
\cite{Den}).  

\section{Schottky uniformization}

Topologically a compact Riemann surface $X$ of genus $g$ is
obtained by gluing the sides of a $4g$-gon. Correspondingly, the
fundamental group has a presentation
$$ \pi_1(X)=\langle a_1,\ldots, a_g, b_1, \ldots, b_g\,  | \,\,\, \prod_i
[a_i, b_i] =1 \rangle, $$ where the generators $a_i$ and $b_i$
label the sides of the polygon.

In the genus $g=1$ case, the parallelogram is the fundamental
domain of the $\pi_1(X)\simeq \Z^2$ action on the plane $\C$, so
that $X=\C /(\Z+\Z\tau)$ is an elliptic curve.

For genus at least $g\geq 2$, the hyperbolic plane $\H^2$ admits a
tessellation by regular $4g$-gons, and the action of the
fundamental group by deck transformation is realized by the action
of a subgroup $\pi_1(X)\simeq G \subset \PSL(2,\R)$ by isometries
of $\H^2$. This endows the compact Riemann surface $X$ with a
hyperbolic metric and a Fuchsian uniformization
$$ X = G\backslash \H^2. $$

Another, less well known, type of uniformization of compact
Riemann surfaces is Schottky uniformization. We recall briefly
some general facts on Schottky groups.

\subsection{Schottky groups}

A {\em Schottky group} of rank $g$ is a discrete subgroup
$\Gamma\subset \PSL (2,\C)$, which is purely loxodromic and
isomorphic to a free group of rank $g$. The group $\PSL(2,\C)$
acts on $\P^1(\C)$ by fractional linear transformations
$$ \gamma : z \mapsto \frac{(az+b)}{(cz+d)}. $$
Thus, $\Gamma$ also acts on $\P^1(\C)$.

\begin{figure}
\begin{center}
\epsfig{file=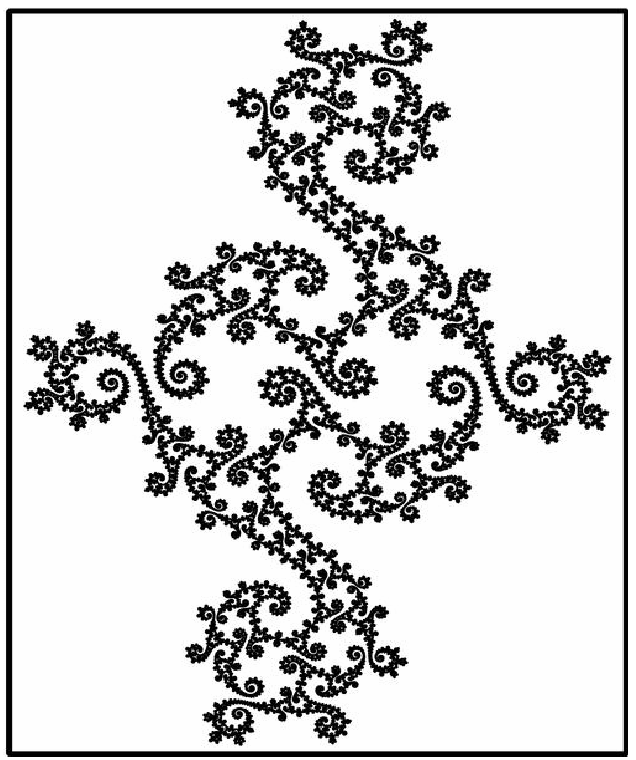}
\end{center}
\caption{Limit set of a Schottky group
\label{Fig-LimSchottky}}
\end{figure}

We denote by $\Lambda_\Gamma$ the {\em limit set} of the action
of $\Gamma$. This is the smallest non--empty closed
$\Gamma$--invariant subset of $\P^1(\C)$. This set can also be
described as the closure of the set of the attractive and
repelling fixed points $z^{\pm}(\gamma)$ of the loxodromic
elements $\gamma\in \Gamma$. In the case $g=1$ the limit set
consists of two points, but for $g\geq 2$ the limit set is usually
a fractal of some Hausdorff dimension $0\leq \delta =\dim_H
(\Lambda_\Gamma) < 2$ (\cf \eg Figure
\ref{Fig-LimSchottky}\footnote{Figures \ref{Fig-LimSchottky} and
\ref{Fig-FS} are taken from ``Indra's pearls'' by
Mumford, Series, and Wright, \cite{MSW}}).

We denote by $\Omega_\Gamma$ the {\em domain of discontinuity} of
$\Gamma$, that is, the complement of $\Lambda_\Gamma$ in
$\P^1(\C)$. The quotient
\begin{equation}
X(\C) = \Gamma \backslash \Omega_\Gamma \label{RS}
\end{equation}
is a Riemann surface of genus $g$ and the covering $\Omega_\Gamma
\to X(\C)$ is called a {\em Schottky uniformization} of $X(\C)$.

Every compact Riemann surface $X(\C)$ admits a Schottky
uniformization.

\begin{figure}
\begin{center}
\epsfig{file=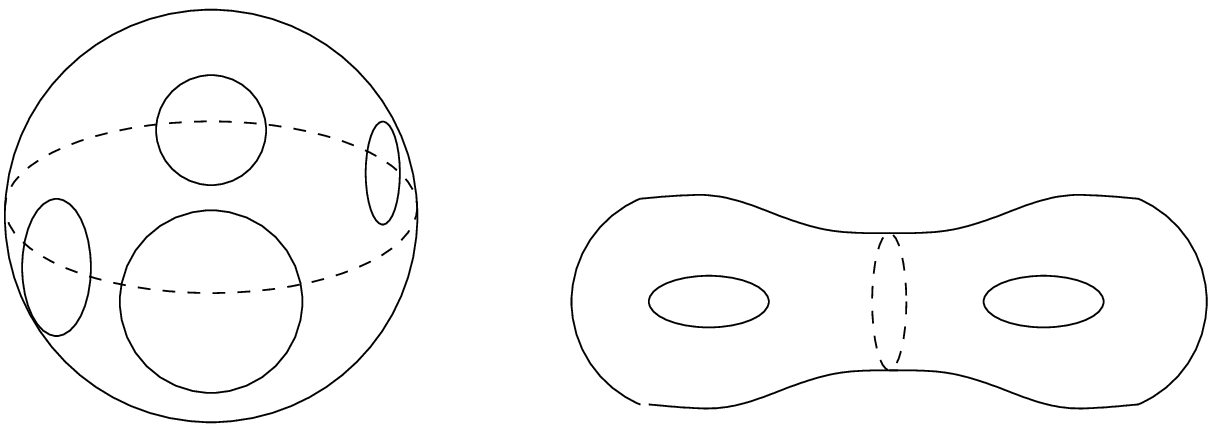}
\end{center}
\caption{Schottky uniformization for $g=2$
\label{Fig-Schottky2}}
\end{figure}

Let $\{ \gamma_i \}_{i=1}^g$ be a set of generators of the
Schottky group $\Gamma$. We write $\gamma_{i+g}=\gamma_i^{-1}$.
There are $2g$ Jordan curves $\gamma_k$ on the sphere at infinity
$\P^1(\C)$, with pairwise disjoint interiors $D_k$, such that the
elements $\gamma_k$ are given by fractional linear transformations
that map the interior of $C_k$ to the exterior of $C_{j}$ with
$|k-j|=g$. The curves $C_k$ give a {\em marking} of the Schottky
group. The markings are circles in the case of {\em classical}
Schottky groups. A fundamental domain for the action of a
classical Schottky group $\Gamma$ on $\P^1(\C)$ is the region
exterior to $2g$-circles (\cf Figure \ref{Fig-Schottky2}).

\subsection{Schottky and Fuchsian}

Notice that, unlike Fuchsian uniformization, where the covering
$\H^2$ is the universal cover, in the case of Schottky
uniformization $\Omega_\Gamma$ is very far from being simply
connected, in fact it is the complement of a Cantor set.

The relation between Fuchsian and Schottky uniformization is given
by passing to the covering that corresponds to the normal subgroup
$N\langle a_1,\ldots,a_g \rangle$ of $\pi_1(X)$ generated by half
the generators $\{ a_1, \ldots, a_g \}$,
$$ \Gamma \simeq \pi_1(X)/ N\langle a_1,\ldots,a_g \rangle, $$
with a corresponding covering map
$$
\diagram \H^2 \rrto^J \drto_{\pi_G} & &
\Omega_\Gamma\dlto^{\pi_\Gamma}   \\ & X & &
\enddiagram
$$
At the level of moduli, there is a corresponding map between
Teichm\"uller space $\T_g$ and Schottky space $\cS_g$, which
depends on $3g-3$ complex moduli.

\subsection{Surface with boundary: simultaneous uniformization}

To better visualize geometrically the Schottky uniformization of a
compact Riemann surface, we can relate it to a simultaneous
uniformization of the upper and lower half planes that yields two
Riemann surfaces with boundary, joined at the boundary.

\begin{figure}
\begin{center}
\epsfig{file=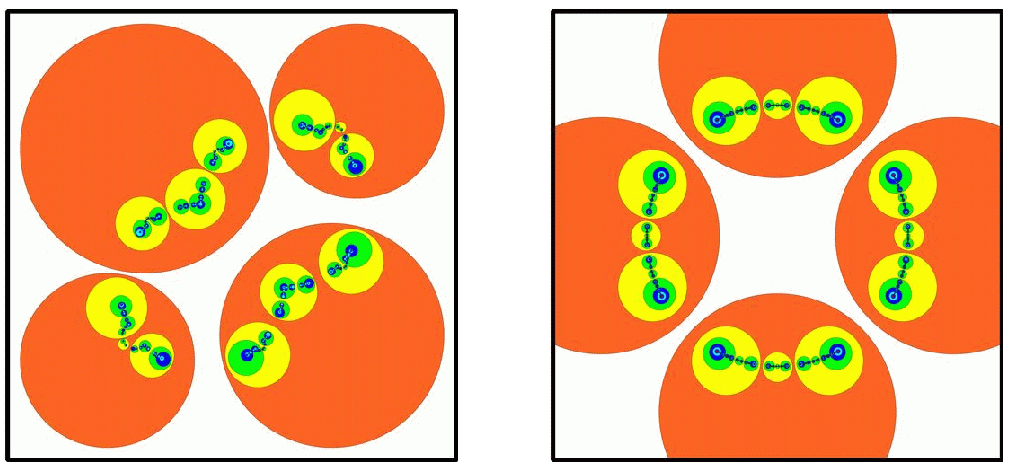}
\end{center}
\caption{Classical and Fuchsian Schottky groups
\label{Fig-FS}}
\end{figure}

A Schottky group that is specified by real parameters so that it
lies in $\PSL(2,\R)$ is called a {\em Fuchsian Schottky group} (\cf
Figure \ref{Fig-FS}).
Viewed as a group of isometries of the hyperbolic plane $\H^2$, or
equivalently of the Poincar\'e disk, a Fuchsian Schottky group $G$
produces a quotient $G\backslash \H^2$, which is topologically a
Riemann surface with boundary.

A {\em quasi-circle} for $\Gamma$ is a Jordan curve $C$ in
$\P^1(\C)$ which is invariant under the action of $\Gamma$. In
particular, such curve contains the limit set $\Lambda_\Gamma$. It
was proved by Bowen that, if $X(\C)$ is a Riemann surface of genus
$g\geq 2$ with Schottky uniformization, then there exists always a
quasi-circle for $\Gamma$.

We have then $\P^1(\C)\setminus C = \Omega_1 \cup \Omega_2$ and,
for $\pi_\Gamma : \Omega_\Gamma \to X(\C)$ the covering map,
$$ \hat C = \pi_\Gamma (C \cap \Omega_\Gamma) \subset X(\C) $$
is a set of curves on $X(\C)$ that disconnect the Riemann surface
in the union of two surfaces with boundary, uniformized
respectively by $\Omega_1$ and $\Omega_2$.

\begin{figure}
\begin{center}
\epsfig{file=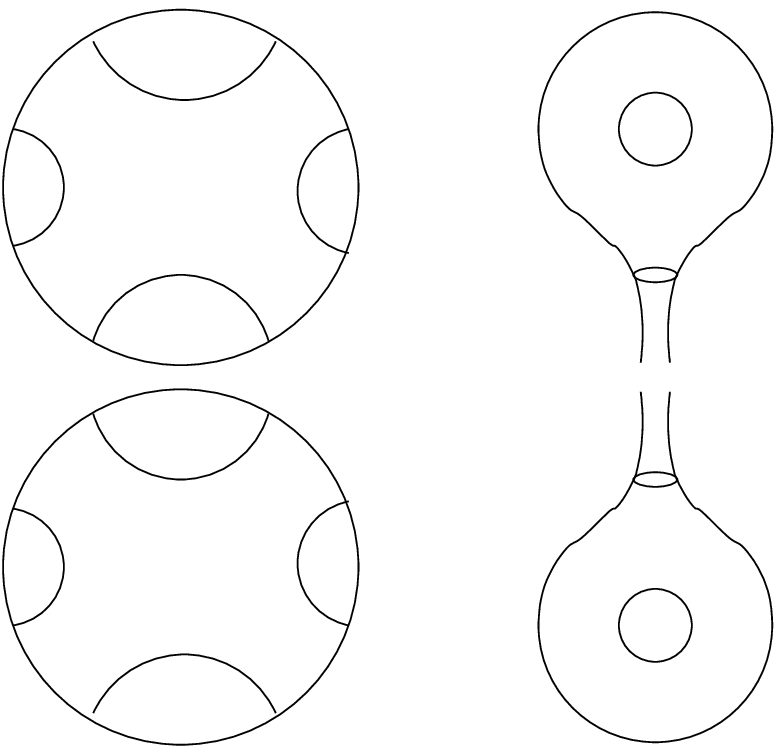}
\end{center}
\caption{Fuchsian Schottky groups: Riemann surfaces with boundary
\label{Fig-Schottky2b}}
\end{figure}

There exist conformal maps
\[ \alpha_i: \Omega_i \stackrel{\simeq}{\to} U_i,\quad U_1 \cup
U_2 = \P^1(\C) \setminus \P^1(\R)
\]
with $U_i \simeq \H^2 =$ upper and lower half planes in
$\P^1(\C)$, and with
$$ G_i :=\{\alpha_i\gamma{\alpha_i}^{-1}~:~\gamma\in\tilde\Gamma\} \simeq
\Gamma $$ Fuchsian Schottky groups $G_i \subset\PSL(2,\R)$. Here
$\tilde\Gamma \subset \SL(2,\R)$ is the $\Gamma$-stabilizer of
each of the two connected components in $\P^1(\C)\setminus C$.

The compact Riemann surface $X(\C)$ is thus obtained as
$$ X(\C) = X_1 \cup_{\partial X_1 = \hat C = \partial X_2}X_2, $$
with $X_i = U_i/G_i$  Riemann surfaces with boundary $\hat C$ (\cf
Figure \ref{Fig-Schottky2b}).

In the case where $X(\C)$ has a real structure $\iota: X\to X$,
and the fixed point set $Fix(\iota)=X(\R)$ of the involution
$\iota$ is nonempty, we have in fact $\hat C = X(\R)$, and the
quasi-circle is given by $\P^1(\R)$.

Notice that, in the case of a Fuchsian Schottky group, the
Hausdorff dimension $\dim_H \Lambda_\Gamma$ of the limit set is in
fact bounded above by $1$, since $\Lambda_\Gamma$ is contained in
the rectifiable quasi-circle $\P^1(\R)$.

\subsection{Hyperbolic Handlebodies}

\begin{figure}
\begin{center}
\epsfig{file=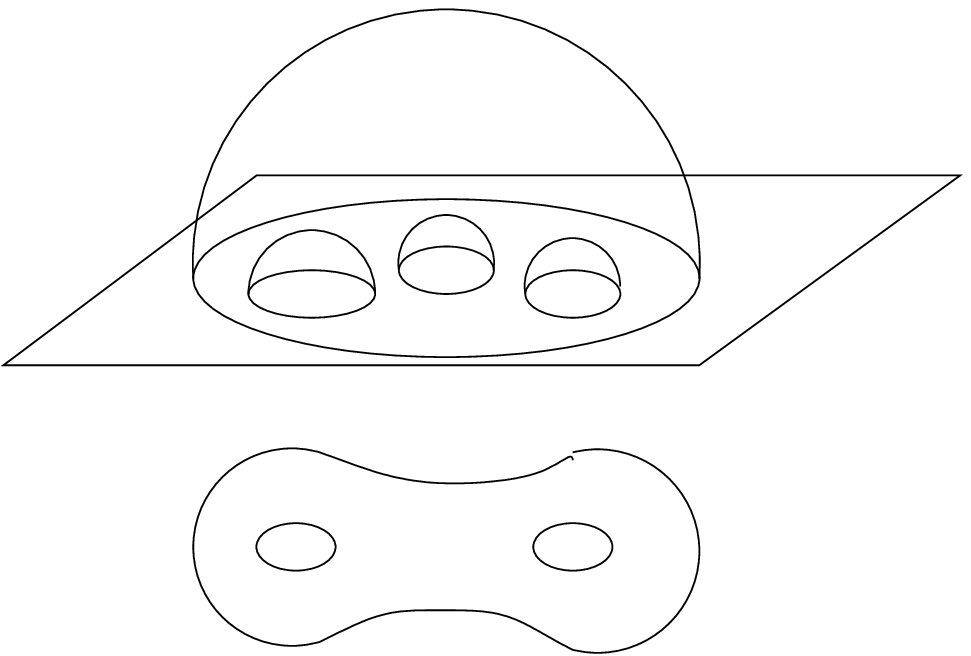}
\end{center}
\caption{Genus two: fundamental domain in $\H^3$
\label{Fig-3dSchottky}}
\end{figure}

The action of a rank $g$ Schottky group $\Gamma\subset \PSL(2,\C)$
on $\P^1(\C)$, by fractional linear transformations, extends to an
action by isometries on real hyperbolic 3-space $\H^3$. For a
classical Schottky group, a fundamental domain in $\H^3$ is given by
the region external to $2g$ half spheres over the circles $C_k\subset
\P^1(\C)$ (\cf Figure \ref{Fig-3dSchottky}).

The quotient
\begin{equation}\label{handlebodies}
\mX_\Gamma = \H^3 /\Gamma
\end{equation}
is topologically a handlebody of genus $g$ filling the Riemann
surface $X(\C)$.

Metrically, $\mX_\Gamma$ is a real hyperbolic 3-manifold of
infinite volume, having $X(\C)$ as its conformal boundary at
infinity $X(\C)=\partial \mX_\Gamma$.

We denote by $\overline{\mX}_\Gamma$ the compactification obtained
by adding the conformal boundary at infinity,
\begin{equation}\label{handle-compact}
\overline{\mX}_\Gamma = (\H^3 \cup \Omega_\Gamma)/ \Gamma.
\end{equation}

In the genus zero case, we just have the sphere $\P^1(\C)$ as the
conformal boundary at infinity of $\H^3$, thought of as the unit
ball in the Poincar\'e model.

In the genus one case we have a solid torus $\H^3/q^\Z$, for $q\in
\C^*$ acting as $$ q (z,y) =(q z, |q| y) $$ in the upper half
space model, with conformal boundary at infinity the Jacobi
uniformized elliptic curve $\C^* /q^\Z$.

In this case, the limit set consists of the point $\{0, \infty\}$,
the domain of discontinuity is $\C^*$ and a fundamental domain is
the annulus $\{ |q|< |z| \leq 1 \}$ (exterior of two circles).

The relation of Schottky uniformization to the usual Euclidean
uniformization of complex tori $X=\C/(\Z+\tau\Z)$ is given by $q=
\exp(2\pi i \tau)$.

We shall return to the genus one case where we discuss a physical
interpretation of Manin's result on the Green function. To our
purposes, however, the most interesting case is when the genus is
$g\geq 2$. In this case, the limit set $\Lambda_\Gamma$ is a
Cantor set with an interesting dynamics of the action of $\Gamma$.
It is the dynamics of the Schottky group on its limit set that
generates an interesting noncommutative space.

\begin{figure}
\begin{center}
\epsfig{file=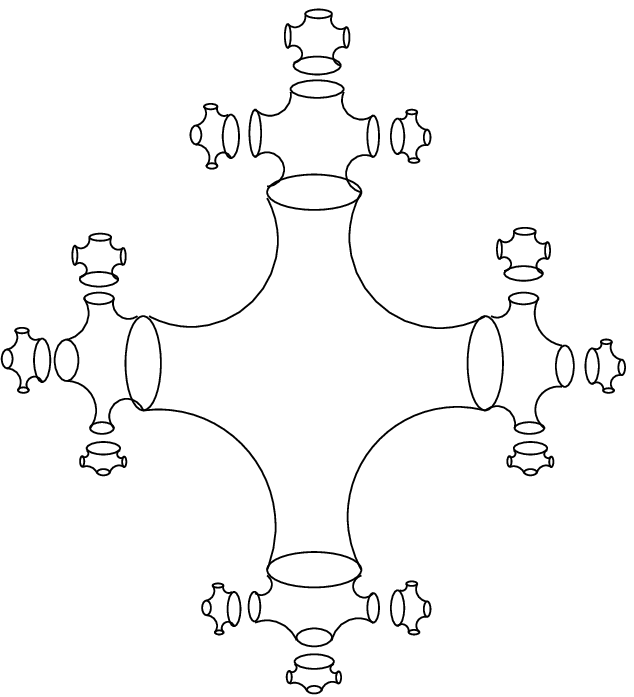}
\end{center}
\caption{Handlebody of genus 2: fundamental domains in $\H^3$
\label{Fig-fundSchottkyH3}}
\end{figure}

\subsection{Geodesics in $\mX_\Gamma$}

The hyperbolic handlebody $\mX_\Gamma$ has infinite volume, but it
contains a region of finite volume, which is a deformation retract
of $\mX_\Gamma$. This is called the ``convex core'' of
$\mX_\Gamma$ and is obtained by taking the geodesic hull of the
limit set $\Lambda_\Gamma$ in $\H^3$ and then the quotient by
$\Gamma$.

We identify different classes of infinite geodesics in
$\mX_\Gamma$.

\begin{itemize}

\item {\em Closed geodesics}: since $\Gamma$ is purely loxodromic,
for all $\gamma \in \Gamma$ there exist two fixed points
$\{z^\pm(\gamma)\} \in \P^1(\C)$. The geodesics in $\H^3 \cup
\P^1(\C)$ with ends at two such points $\{z^\pm(\gamma)\}$, for
some $\gamma \in \Gamma$, correspond to closed geodesics in the
quotient $\mX_\Gamma$.

\medskip

\item {\em Bounded geodesics}: The images in $\mathfrak
X_\Gamma$ of geodesics in $\H^3 \cup \P^1(\C)$ having both ends on
the limit set $\Lambda_\Gamma$ are geodesics that remain confined
within the convex core of $\mX_\Gamma$.

\medskip

\item {\em Unbounded geodesics}: these are the geodesics in
$\mX_\Gamma$ that eventually wander off the convex core towards
the conformal boundary $X(\C)$ at infinity. They correspond to
geodesics in $\H^3\cup \P^1(\C)$ with at least one end at a point
of $\Omega_\Gamma$.

\end{itemize}

In the genus one case, there is a unique primitive closed
geodesic, namely the image in the quotient of the geodesic in
$\H^3$ connecting $0$ and $\infty$. The bounded geodesics are
those corresponding to geodesics in $\H^3$ originating at $0$ or
$\infty$.

The most interesting case is that of genus $g\geq 2$, where the
bounded geodesics form a complicated tangle inside $\mX_\Gamma$.
Topologically, this is a generalized solenoid, namely it is
locally the product of a line and a Cantor set.

\section{Dynamics and noncommutative geometry}

Since the uniformizing group $\Gamma$ is a free group, there is a
simple way of obtaining a coding of the bounded geodesics in
$\mX_\Gamma$. The set of such geodesics can be identified with
$\Lambda_\Gamma \times_\Gamma \Lambda_\Gamma$, by specifying the
endpoints in $\H^3\cup \P^1(\C)$ modulo the action of $\Gamma$.

If $\{ \gamma_i \}_{i=1}^g$ is a choice of generators of $\Gamma$
and $\gamma_{i+g}=\gamma_i^{-1}$, $i=1,\ldots, g$, we can
introduce a {\em subshift of finite type} $({\mathcal S},T)$ where
\begin{equation}\label{Sshift}
{\mathcal S}= \{\ldots a_{-m} \ldots a_{-1} a_0 a_1 \ldots a_\ell
\ldots  \, \, |  a_i \in \{ \gamma_i \}_{i=1}^{2g}, \, \, a_{i+1} \neq
a_i^{-1}, \forall i\in \Z \}
\end{equation}
is the set of doubly infinite words in the generators and their
inverses, with the admissibility condition that no cancellations
occur. The map $T$ is the invertible shift
\begin{equation}\label{Tshift}
 T(\ldots a_{-m}
\ldots a_{-1} a_0
 a_1 \ldots  a_{\ell} \ldots ) =
 \ldots a_{-m+1}  \ldots  a_{0} a_1  a_2 \ldots  a_{\ell+1}
\ldots  \end{equation}

Then we can pass from the discrete dynamical system $({\mathcal
S},T)$ to its suspension flow and obtain the mapping torus
\begin{equation}\label{mapTorus}
 {\mathcal S}_T := {\mathcal S} \times [0,1] / (x,0)\sim (Tx,1).
\end{equation}
Topologically this space is a solenoid, that is, a bundle over
$S^1$ with fiber a Cantor set.

\subsection{Homotopy quotient}

The space $\cS_T$ of \eqref{mapTorus} has a natural interpretation
in noncommutative geometry as the homotopy quotient in the sense
of Baum--Connes \cite{BaumConnes} of the noncommutative space
given by the $C^*$-algebra
\begin{equation}\label{CSTprod}
C(\cS)\rtimes_T \Z
\end{equation}
describing the action of the shift \eqref{Tshift} on the totally
disconnected space \eqref{Sshift}. The noncommutative space
\eqref{CSTprod} parameterizes bounded geodesics in the handlebody
$\mX_\Gamma$.

The homotopy quotient is given by
$$ \cS\times_\Z \R = \cS_T. $$

The $K$-theory of the $C^*$-algebra \eqref{CSTprod} can be
computed via the Pimsner--Voiculescu six terms exact sequence,
$$
\diagram K_1(A) \rto & K_0(C(\cS))\rto^{\delta=1-T} & K_0(C(\cS))
\dto \\ K_1(C(\cS)) \uto & K_1(C(\cS))\lto & K_0(A)\lto
\enddiagram
$$
where $A=C(\cS)\rtimes_T \Z$. Here, since the space $\cS$ is
totally disconnected, we have $K_1(C(\cS))=0$ and
$K_0(C(\cS))=C(\cS,\Z)$ (locally constant integer valued functions).
Thus the exact sequence becomes
\begin{equation}\label{PVseqT}
0 \to K_1(A) \to C(\cS,\Z) \stackrel{\delta=1-T}{\to} C(\cS,\Z)
\to K_0(A) \to 0,
\end{equation}
and we obtain $K_1(A)=\Z=\Ker(\delta)$ and
$K_0(A)=\Coker(\delta)$. In dynamical systems language, these are
respectively the {\em invariants} and {\em coinvariants} of the
invertible shift $T$ (\cf \cite{BoHa} \cite{PaTu}).

In terms of the homotopy quotient, one can describe this exact
sequence more geometrically in terms of the Thom isomorphism and
the $\mu$-map
$$ \mu: K^{*+1} ({\mathcal S}_T)\cong H^{*+1}({\mathcal S}_T,\Z) \to
K_*({\rm C}({\mathcal S})\rtimes_T \Z). $$

Thus, we obtain:
$$
\begin{array}{l}
K_1(A)\simeq H^0({\mathcal S}_T)=\Z  \\[2mm]
K_0(A)\simeq H^1({\mathcal S}_T)  \end{array}
$$
The $H^1(\cS_T)$ can be identified with the \v{C}ech cohomology
group given by the homotopy classes $[{\mathcal S}_T,U(1)]$, by
mapping $f\mapsto [\exp(2\pi i t f(x))]$, for $f\in
C(\cS,\Z)/\delta$.

\subsection{Filtration}

The identification
$$ H^1({\mathcal S}_T,\Z) \cong K_0({\rm
C}({\mathcal S}) \rtimes_T \Z)
$$
of the cohomology of ${\mathcal S}_T$ with the $K_0$-group of the
crossed product ${\rm C}^*$-algebra for the action of $T$ on
${\mathcal S}$ endows $H^1({\mathcal S}_T,\Z)$ with a filtration.

\begin{thm}\label{filtrH1}
The first cohomology of $\cS_T$ is a direct limit
$$ H^1(\cS_T) = \varinjlim_n F_n, $$
of free abelian groups $F_0 \hookrightarrow F_1 \hookrightarrow
\cdots F_n\hookrightarrow \cdots$ of ranks ${\rm rank}\, F_0=2g$
and ${\rm rank}\, F_n = 2g(2g-1)^{n-1}(2g-2) +1$, for $n\geq 1$.
\end{thm}

In fact, by the Pimsner-Voiculescu six terms exact sequence, the
group $K_0({\rm C}({\mathcal S}) \rtimes_T \Z)$ can be identified
with the cokernel of the map $1-T$ acting as $f \mapsto f-f\circ
T$ on the $\Z$ module ${\rm C}({\mathcal S},\Z)\simeq K_0({\rm
C}({\mathcal S}))$.

Then the filtration is given by the submodules of functions
depending only on the $a_0 \ldots a_n$ coordinates in the doubly
infinite words describing points in ${\mathcal S}$. Namely, one
first shows that
$$ H^1(\cS_T,\Z) = C(\cS,\Z)/ \delta (C(\cS,\Z)) = \cP/\delta \cP $$
where $\cP$ denotes the $\Z$-module of locally constant
$\Z$-valued functions that depend only on ``future coordinates''.
This can be identified with functions on the limit set
$\Lambda_\Gamma$, since each point in $\Lambda_\Gamma$ is
described by an infinite (to the right) admissible sequence in the
generators $\gamma_i$ and their inverses. We have $\cP \simeq
C(\Lambda_\Gamma,\Z)$.

The module $\cP$ clearly has a filtration by the submodules
$\cP_n$ of functions of the first $n+1$ coordinates. These have
${\rm rank}\, \cP_n = 2g(2g-1)^n$.

We set $F_n := \cP_n/\delta \cP_{n-1}$. On can show that there are
induced injections $F_n \hookrightarrow F_{n+1}$ and that
$$ H^1(\cS_T) = \varinjlim_n F_n. $$

Moreover, we have ${\rm rank} F_n=\theta_n - \theta_{n-1} +1$,
where $\theta_n$ is the number of admissible words of length
$n+1$. All the $\Z$ modules $F_n$ and the quotients $F_n/ F_{n-1}$
are torsion free (\cf \cite{PaTu}).

\subsection{Hilbert space and grading}

It is convenient to consider the complex vector space
$$ {\mathcal P}_\C=C(\Lambda_\Gamma,\Z) \otimes \C $$
and the corresponding exact sequence computing the cohomology with
$\C$ coefficients:
\begin{equation}\label{PVC}
 0 \to \C \to {\mathcal P}_\C \stackrel{\delta}{\longrightarrow}
{\mathcal P}_\C \to H^1({\mathcal S}_T,\C) \to 0
\end{equation}

The complex vector space $\cP_\C$ sits in the Hilbert space
$$ \cP_\C \subset \cL=L^2(\Lambda_\Gamma, d\mu), $$
where $\mu$ is the Patterson--Sullivan measure on the limit set,
satisfying
$$ \gamma^* d\mu = |\gamma'|^{\dim_H(\Lambda_\Gamma)}
d\mu, $$ with $\dim_H(\Lambda_\Gamma)$ the Hausdorff dimension.

This gives a Hilbert space $\cL$, together with a filtration
$\cP_n$ by finite dimensional subspaces. In such setting, it is
natural to consider a corresponding {\em grading operator},
\begin{equation}\label{Dgrading}
 D= \sum_n n \hat\Pi_n,
\end{equation}
where $\Pi_n$ denotes the orthogonal projection onto $\cP_n$ and
$\hat\Pi_n = \Pi_n -\Pi_{n-1}$.

\subsection{Cuntz--Krieger algebra}

There is a noncommutative space that encodes nicely the dynamics
of the Schottky group $\Gamma$ on its limit set $\Lambda_\Gamma$.
Consider the $2g\times 2g$ matrix $A$ that gives the admissibility
condition for sequences in ${\mathcal S}$: this is the matrix with
$\{ 0,1 \}$ entries satisfying $A_{ij}=1$ for $|i-j|\neq g$ and
$A_{ij}=0$ otherwise.

The Cuntz--Krieger algebra ${\mathcal O}_A$ associated to this
matrix is the universal $C^*$-algebra generated by partial
isometries $S_i$, $i=1,\ldots, 2g$, satisfying the relations
$$ \sum_j S_j S^*_j =1 \ \ \ \  S_i^* S_i = \sum_j A_{ij} S_j S^*_j. $$
Recall that a partial isometry is an operator $S$ satisfying
$S=SS^*S$.

This algebra is related to the Schottky group by the following
result.

\begin{prop}\label{crossOA}
There is an isomorphism
\begin{equation}\label{OA-cross1}
{\mathcal O}_A \cong {\rm C}(\Lambda_\Gamma)\rtimes \Gamma.
\end{equation}
\end{prop}

Up to a stabilization (tensoring with compact operators), the
algebra has another crossed product description as
\begin{equation}\label{OA-cross2}
{\mathcal O}_A \simeq {\mathcal F}_A \rtimes_T \Z,
\end{equation}
with ${\mathcal F}_A$ an AF-algebra (a direct limit of finite
dimensional $C^*$-algebras).

Consider the cochain complex of Hilbert spaces
$$ 0 \to \C \to \cL \stackrel{\delta}{\to}\cL \to \cH \to 0 $$
determined by the Pimsner--Voiculescu sequence \eqref{PVseqT}.

\begin{prop}\label{represOA}
The $C^*$-algebra ${\mathcal O}_A$ admits a faithful
representation on the Hilbert space ${\mathcal
L}=L^2(\Lambda_\Gamma, d\mu)$.
\end{prop}

This is obtained as follows.

For $d_H=\dim_H(\Lambda_\Gamma)$ the Hausdorff dimension, consider
the operators
\begin{equation}\label{operPi}
 P_i\, f = \chi_{\gamma_i}\, f \ \ \ T_i\, f = | (\gamma_i^{-1})' |^{d_H/2}
\, f\circ \gamma_i^{-1},
\end{equation}
where $\{ \gamma_i \}$ are the generators of $\Gamma$ and their
inverses and
\begin{equation}\label{Tgamma}
T_\gamma\, f = |\gamma'|^{d_H/2} f\circ \gamma.
\end{equation}

Then the operators
\begin{equation}\label{Si-rep}
 S_i = \sum_j A_{ij} T_i^* P_j
\end{equation}
are partial isometries on $\cL$ satisfying the Cuntz--Krieger
relations for the matrix $A$ of the subshift of finite type
\eqref{Sshift}. This gives the representation of $\cO_A$.

\subsection{Spectral triple for Schottky groups}

On the Hilbert space $\cH=\cL\oplus\cL$ consider the diagonal
action of the algebra $\cO_A$ and the Dirac operator $\sD$ defined
as
\begin{equation}\label{DiracDH}
 \begin{array}{ll} \sD|_{\cL\oplus 0} & = \sum_n\, (n+1)\,\,
(\hat\Pi_n
\oplus 0) \\[3mm]
 \sD|_{0\oplus \cL} & = -\sum_n\, n\,\, (0\oplus  \hat\Pi_n). \end{array}
\end{equation}

The choice of the sign in \eqref{DiracDH} is not optimal from the
point of view of the $K$-homology class. A better choice would be
$$ F= \left(\begin{array}{cc} 0 & 1 \\ 1 & 0 \end{array}\right). $$
This would require in turn a modification of $|D|$. A construction
along these lines is being considered in joint work of the author
with Alina Vdovina and Gunther Cornelissen. In our setting here,
the reason for the choice \eqref{DiracDH} in \cite{CM} is the
formula \eqref{DiracFrobeniusCone} relating $\sD$ to the
``logarithm of Frobenius'' at arithmetic infinity.

\begin{thm}\label{triple}
For a Schottky group $\Gamma$ with $\dim_H(\Lambda_\Gamma)<1$, the
data $(\cO_A,\cH,\sD)$, for $\cH=\cL\oplus\cL$ with the diagonal
action of $\cO_A$ through the representation \eqref{Si-rep} and
the Dirac operator \eqref{DiracDH}, define a non-finitely
summable, $\theta$-summable spectral triple.
\end{thm}

The key point of this result is the compatibility relation between
the algebra and the Dirac operator, namely the fact that the
commutators $[\sD,a]$ are bounded operators, for all $a\in
\cO_A^{alg}$, the involutive subalgebra generated algebraically by
the $S_i$ subject to the Cuntz--Krieger relations.

This follows by an estimate on the norm of the commutators $\| [
\sD, S_i ] \|$ and $\| [ \sD, S_i^* ] \|$, in terms of the
Poincar\'e series of the Schottky group ($d_H<1$)
$$ \sum_{\gamma\in \Gamma} |\gamma\,'|^s, \ \ \  s=1> d_H, $$
where the Hausdorff dimension $d_H$ is the exponent of convergence
of the Poincar\'e series.

The dimension of the $n$-th eigenspace of $\sD$ is
$2g(2g-1)^{n-1}(2g-2)$ for $n\geq 1$, $2g$ for $n=0$ and
$2g(2g-1)^{-n-1}(2g-2)$ for $n\leq -1$, so the spectral triple is
not finitely summable, since $|\sD|^z$ is not of trace class. It
is $\theta$ summable, since the operator $\exp(-t\sD^2)$ is of
trace class, for all $t>0$.

Using the description \eqref{OA-cross2} of the noncommutative
space as crossed product of an AF algebra by the action of the
shift, $\cF_A \rtimes_T \Z$, one may be able to find a 1-summable
spectral triple. Here the dense subalgebra should not contain any
of the group elements. In fact, by the result of \cite{Co-fredh},
the group $\Gamma$ is a free group, hence its growth is too fast
to support finitely summable spectral triples on its group ring.

\section{Arithmetic infinity: archimedean primes}

An algebraic number field $\K$, which is an extension of $\Q$ with
 $[\K:\Q]=n$ admits $n=r_1+2r_2$ embeddings
\begin{equation}\label{embeddings}
 \alpha: \K \hookrightarrow \C .
\end{equation}
These can be subdivided into $r_1$ real embeddings $\K \hookrightarrow
\R$ and $r_2$ pairs of complex conjugate embeddings.
The embeddings \eqref{embeddings} are called the {\em archimedean
primes} of the number field. The set of archimedean primes is often
referred to as ``arithmetic infinity'', a terminology borrowed from
the case of the unique embedding of $\Q \hookrightarrow \R$, which is
called the ``infinite prime''.

A general strategy in arithmetic geometry is to adapt the tools of
classical algebraic geometry to the arithmetic setting. In particular,
over $\Q$ the set of primes $\Sp(\Z)$ is the analog in arithmetic
geometry of the affine line. It becomes clear then
that some compactification is necessary, at least in order to have
a well behaved form of intersection theory in arithmetic
geometry. Namely, we need to
pass from the ``affine'' $\Sp(\Z)$ to the projective case.
The compactification is obtained by adding the infinite prime to the
set of finite primes. A goal of arithmetic geometry then becomes
developing a setting that treats the infinite prime and the finite
primes of equal footing.

More generally, for a number field $\K$ with $O_\K$ its ring of
integers, the set of primes $\Sp(O_\K)$ is compactified by adding
the set of ``archimedean primes''
\begin{equation}\label{SpOK-comp}
\overline{\Sp(O_\K)} = \Sp(O_\K) \cup \{ \alpha: \K \hookrightarrow \C
\}.
\end{equation}

\subsection{Arithmetic surfaces}

Let $X$ be a smooth projective algebraic curve defined over
$\Q$. Then, by clearing denominators one obtains an equation with $\Z$
coefficients. This determines a scheme $X_\Z$ over $\Sp(\Z)$,
$$ X_\Z \otimes_\Z \Sp(\Q)= X, $$
where the closed fiber of $X_{\Z}$ at a prime $p\in \Sp(\Z)$ is the
reduction $X \mod p$.
Thus, viewed as an arithmetic varities, an algebraic curve becomes a
2-dimensional fibration over the affine line $\Sp(\Z)$.

One can also consider reductions of $X$ (defined over $\Z$)
modulo $p^n$ for some prime $p\in \Sp(\Z)$. The limit as $n\to\infty$
defines a $p$--adic completion of $X_{\Z}$. This can be thought of as
an ``infinitesimal neighborhood'' of the fiber at $p$.

The picture is more complicated at arithmetic infinity, since one does
not have a suitable notion of ``reduction mod $\infty$'' available to
define the closed fiber. On the other hand, one does have the analog
of the $p$--adic completion at hand. This is given by the Riemann
surface (smooth projective algebraic curve over $\C$) determined by
the equation of the algebraic curve $X$ over $\Q$, under the
embedding of $\Q\subset \C$,
$$ X(\C) = X \otimes_\Q \Sp(\C) $$
with the absolute value $|\cdot |$ at the infinite
prime replacing the $p$-adic valuations.

Similarly, for $\K$ a number field and $O_\K$ its ring of integers, a
smooth proper algebraic curve $X$ over $\K$ determines a smooth
minimal model $X_{O_\K}$, which defines an arithmetic surface ${\mathcal
X}_{O_\K}$ over $\Sp(O_\K)$. The closed fiber $X_{\wp}$ of
${\mathcal X}_{O_\K}$ over a prime $\wp\in O_\K$ is given by the
reduction mod $\wp$.

When $\Sp(O_\K)$ is compactified by adding the
archimedean primes, one also obtains $n$ Riemann surfaces
$X_\alpha(\C)$, obtained from the equation defining $X$ over $\K$
under the embeddings $\alpha: \K \hookrightarrow \C$.
Of these Riemann surfaces, $r_1$ carry a real involution.

\begin{figure}
\begin{center}
\epsfig{file=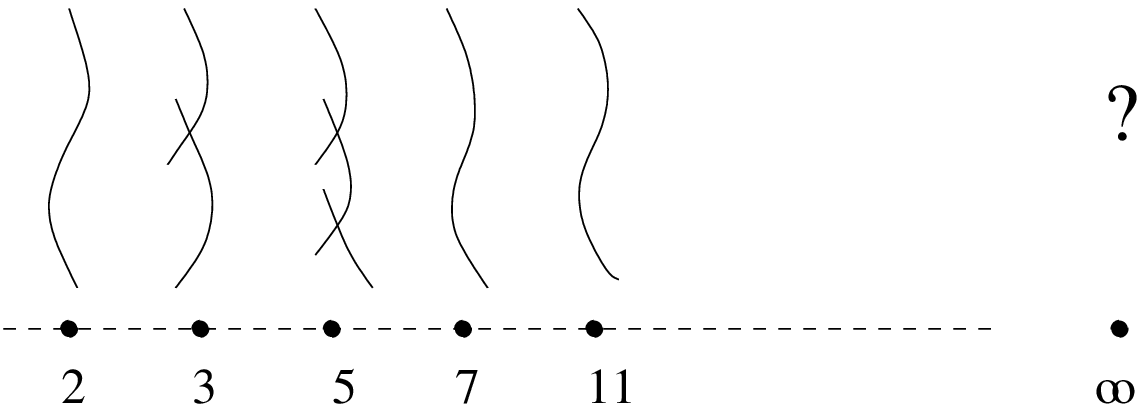}
\end{center}
\caption{Arithmetic surface over $\Sp(\Z)$
\label{Fig-arithmsurf}}
\end{figure}

Thus, the picture of an arithmetic surface over $\overline{\Sp(O_\K)}$
is as follows:
$$
\diagram
X_\wp \,\,\ar[d] \ar@{^{(}->}[r] & X_{\Sp(O_\K)} \ar[d] \ar@{^{(}->}[r]  &
X_{\overline{\Sp(O_\K)}}\ar[d] & \,\, ???\, \ar@{_{(}->}[l] \ar[d] \\
\wp \,\,\ar@{^{(}->}[r] & \Sp(O_\K) \ar@{^{(}->}[r] & \overline{\Sp(O_\K)}
& \,\,\alpha\, \ar@{_{(}->}[l]
\enddiagram
$$
where we do not have an explicit geometric description of the closed
fibers over the archimedean primes (Figure \ref{Fig-arithmsurf}).

Formally, one can enlarge the group of divisors on the arithmetic
surface by adding formal real linear
combinations of irreducible ``closed vertical fibers at infinity''
$\sum_{\alpha} \lambda_\alpha F_\alpha$. Here  the fibers
$F_\alpha$ are only treated as formal symbols, and no geometric
model of such fibers is provided. The remarkable fact is that
Hermitian geometry on the Riemann surfaces $X_\alpha(\C)$ is
sufficient to specify the contribution of such divisors to
intersection theory on the arithmetic surface, even without an
explicit knowledge of the closed fiber.

The main idea of Arakelov geometry is that it is sufficient to work
with the ``infinitesimal neighborhood'' $X_\alpha(\C)$ of the fibers
$F_\alpha$, to have well defined intersection indices.

If one thinks, by analogy, to the case of the classical geometry of a
degeneration of algebraic curves over a disk $\Delta$, with a special
fiber over $0$ the analogous statement would be saying that the
geometry of the special fiber is completely determined by the generic
fiber. This is a very strong statement on the form of the
degeneration: for instance blowing up points on the special fiber is
not seen by just looking at the generic fiber. Investigating this
analogy leads one to expect that the fiber at infinity should behave
like the {\em totally degenerate} case.
This is the case where one has maximal degeneration, where all the
components of the closed
fiber are $\P^1$'s and the geometry of the degeneration is
completely encoded by the {\em dual graph}, which describes in a
purely combinatorial way how these $\P^1$'s are joined. The dual
graph has a vertex for each component of the closed fiber and an
edge for each double point.

The local intersection multiplicities of two
finite, horizontal, irreducible divisors $D_1$, $D_2$ on
$X_{O_\K}$ is given by
\[
[D_1,D_2] = [D_1,D_2]_{fin} + [D_1,D_2]_{inf}
\]
 where the first term counts the contribution from the finite places
 (\ie what happens over $\Sp(O_\K)$) and the second term is the contribution
 of the archimedean primes, \ie the part of the intersection that happens
 over arithmetic infinity.

While the first term is computed in algebro geometric terms, from the
local equations for the divisors $D_i$ at $P$, the second term is
defined as a sum of values of Green functions $g_\alpha$ on the
Riemann surfaces $X_\alpha(\C)$,
\[
[D_1,D_2]_{inf} =
-\sum_{\alpha}\epsilon_\alpha(\sum_{\beta,\gamma}
g_\alpha(P_{1,\beta}^{\alpha},P_{2,\gamma}^{\alpha})),
\]
at points
 $$\{P_{i,\beta}^{\alpha}~|~\beta =
1,\ldots,[\K(D_i):\K]\}\subset X_{\alpha}(\C),$$
for a finite extension $\K(D_i)$ of $\K$ determined by $D_i$.
Here $\epsilon_\alpha =1$ for real embeddings and $=2$ for complex
embeddings.

For a detailed account of these notions of Arakelov geometry, one can
refer to \cite{CoSilv} \cite{Lang-Arak}.

Further evidence for the similarity between the archimedean and the
totally degenerate fibers came from an explicit
computation of the Green function at the archimedean places
derived by Manin \cite{Man-hyp} in terms of a Schottky
uniformization of the Riemann surface $X_\alpha(\C)$. Such
uniformization has an analog at a finite prime, in terms of p-adic
Schottky groups, only in the totally degenerate case.
Another source of evidence comes from a
cohomological theory of the local factors at the archimedean primes,
developed by Consani \cite{KC}, which shows that the resulting
description of the
local factor as regularized determinant at the archimedean primes
resembles mostly the case of the totally degenerate reduction at a
finite prime.

We will present both results in the light of the noncommutative space
$(\cO_A,\cH,\sD)$ introduced in the previous section. As showed by
Consani and the author in \cite{CM} \cite{CM1}
\cite{CM2} \cite{CM3} the noncommutative geometry of this space is
naturally related to both Manin's result on the
Arakelov Green function and the cohomological construction
of Consani.

\section{Arakelov geometry and hyperbolic geometry}

In this section we give a detailed account of Manin's result
\cite{Man-hyp} on the 
relation between the Arakelov Green function on a Riemann surface
$X(\C)$ with Schottky uniformization and geodesics in the
3-dimensional hyperbolic handlebody $\mX_\Gamma$.
Our exposition follows closely the seminal paper \cite{Man-hyp}.

\subsection{Arakelov Green function}

Given a divisor $A=\sum_x m_x(x)$ with support $|A|$ on a 
smooth compact Riemann surface $X(\C)$, and a choice of a positive
real--analytic 2-form $d\mu$ on $X(\C)$, the Green function
$g_{\mu ,A}=g_A$ is a real analytic function on on $X(\C) \setminus
|A|$, uniquely determined by the following conditions:

\begin{itemize}

\item {\it Laplace equation}: $g_{A}$ satisfies
$\partial\bar{\partial}\,g_A=\pi i\,(\deg (A)\,d\mu -\delta_A)$,
with $\delta_A$ the $\delta$--current $\varphi\mapsto \sum_x
m_x\varphi (x)$.

\medskip

\item {\it Singularities}: if $z$ is a local coordinate in a
neighborhood of $x$, then $g_A - m_x\log |z|$ is locally real
analytic.

\medskip

\item {\it Normalization}: $g_{A}$ satisfies $\int_X g_A d\mu =0$.

\end{itemize}

\medskip

If $B=\sum_y n_y(y)$ is another divisor, such that $|A|\cap
|B|=\emptyset$, then the expression $g_{\mu}(A,B):=\sum_y
n_yg_{\mu , A}(y)$ is symmetric and biadditive in $A,B$.
Generally, such expression $g_{\mu}$ depends on $\mu$, where the
choice of $\mu$ is equivalent to the choice of a real analytic
Riemannian metric on $X(\C)$, compatible with the complex
structure.

However, in the special case of degree zero divisors, $\deg A=
\deg B=0$, the $g_{\mu}(A,B)=g(A,B)$ are conformal invariants.

In the case on the Riemann sphere $\P^1(\C)$, if $w_A$ is a
meromorphic function with $Div(w_A)=A$, we have
\begin{equation}\label{divP1}
g(A,B)=\log \prod_{y\in |B|}|w_A(y)|^{n_y} =Re \,
\int_{\gamma_B}\frac{dw_A}{w_A},
\end{equation}
where $\gamma_{B}$ is a 1--chain with
boundary $B$.

In the case of degree zero divisors $A,B$ on a Riemann surface of higher
genus, the formula \eqref{divP1} can be generalized, replacing the
logarithmic differential $dw_A/w_A$ with a
differential of the third kind (meromorphic differential with
nonvanishing residues) $\omega_A$ with purely
imaginary periods and residues $m_x$ at $x$. This gives
\begin{equation} \label{Green} g(A,B)= Re \int_{\gamma_{B}}
\omega_A. \end{equation}

Thus, one can explicitly compute $g(A,B)$ from a basis of
differentials of the third kind with purely
imaginary periods.

\subsection{Cross ratio and geodesics}

\begin{figure}
\begin{center}
\epsfig{file=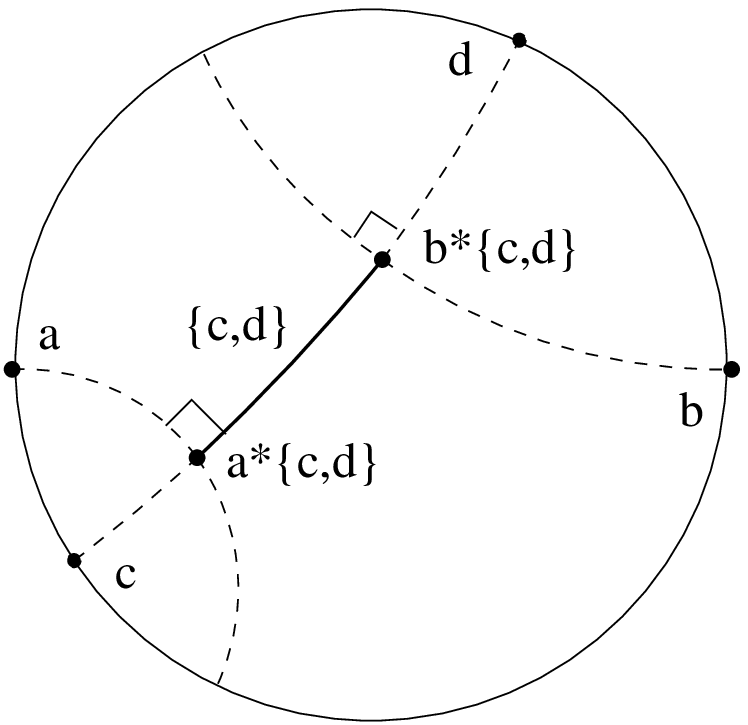}
\end{center}
\caption{Cross ratio and geodesic length
\label{Fig-cross}}
\end{figure}

The basic step leading to the result of Manin, expressing the Arakelov
Green function in terms of geodesics in the hyperbolic handlebody
$\mX_\Gamma$, is a very simple classical fact of hyperbolic
geometry, namely the fact that the cross ratio of four points on
$\P^1(\C)$ can be expressed in terms of geodesics in the
`interior' $\H^3$:
\begin{equation}\label{cr-geod}
 \log | \langle a,b,c,d \rangle | =
-{\rm ordist}\, \left( a*\{c,d\}, b*\{c,d\}
\right).
\end{equation}
Here, following \cite{Man-hyp}, $ordist$ denotes the oriented
distance, and we use the notation $a*\{c,d\}$ to indicate the
point on the geodesic $\{c, d\}$ in $\H^3$ with endpoints $c,d \in
\P^1(\C)$, obtained as the intersection of $\{c, d\}$ with the
unique geodesic from $a$ that cuts $\{c, d\}$ at a right angle
(Figure \ref{Fig-cross}).

\subsection{Differentials and Schottky uniformization}

The next important step in Manin's result \cite{Man-hyp} is to show
that, if $X(\C)=\Gamma\backslash \Omega_\Gamma$ is a Riemann surface
with a Schottky uniformization, then one obtains a basis of
differentials of the third kind with purely
imaginary periods, by taking suitable averages over the group
$\Gamma$ of expressions involving the cross ratio of points on
$\P^1(\C)$.

We denote by $C(|\gamma)$ a set of
representatives of $\Gamma /(\gamma^{\Z})$, by $C(\rho|\gamma)$ a
set of representatives for $(\rho^{\Z})\setminus\Gamma
/(\gamma^{\Z})$, and by $S(\gamma)$ the conjugacy class of
$\gamma$ in $\Gamma$.

Let $w_A$ be a meromorphic function on $\P^1(\C)$ with divisor
$A=(a)-(b)$, such that the support $|A|$ is contained in the
complement of an open neighborhood of $\Lambda_\Gamma$. We use the
notation
\begin{equation}\label{cross-ratio}
 \langle a,b,c,d \rangle := \frac{(a-b)(c-d)}{(a-d)(c-b)}
\end{equation}
for the cross-ratio of points $a,b,c,d\in \P^1(\C)$.

For a fixed choice of a base point $z_0 \in \Omega_{\Gamma}$, the
series
\begin{equation} \label{3k}
 \nu_{(a)-(b)} := \sum_{\gamma \in \Gamma} d \log \langle a,b,\gamma
z ,\gamma z_0 \rangle
\end{equation}
gives the lift to $\Omega_\Gamma$ of a differential of the third
kind on the Riemann surface $X(\C)$, endowed with the choice of
Schottky uniformization. These differentials have residues $\pm 1$
at the images of $a$ and $b$ in $X(\C)$, and they have vanishing
$a_k$ periods, where $\{ a_k, b_k \}_{k=1\ldots g}$ are the generators
of the homology of $X(\C)$.

Similarly, we obtain lifts of differentials of the first kind on
$X(\C)$, by considering the series
\begin{equation}
\omega_\gamma=\sum_{h\in C(|\gamma)}\, d\,\log \langle
hz^+(\gamma), hz^-(\gamma),z,z_0\rangle, \label{1k}
\end{equation}
where we denote by $\{ z^+(\gamma), z^-(\gamma) \}\subset
\Lambda_\Gamma$ the pair of the attractive and repelling fixed
points of $\gamma \in \Gamma$.

The series \eqref{3k} and \eqref{1k} converge absolutely on
compact sets $K\subset \Omega_\Gamma$, whenever $\dim_H
\Lambda_\Gamma <1$. Moreover, they do not depend on the choice of
the base point $z_0\in \Omega_{\Gamma}$.

In particular, given a choice $\{ \gamma_k \}_{k=1}^g$ of generators of
the Schottky group $\Gamma$, we obtain by \eqref{1k} a basis of
holomorphic differentials $\omega_{\gamma_k}$, that satisfy
\begin{equation}\label{normaliz}
 \int_{a_k}\omega_{\gamma_\ell} =2\pi \sqrt{-1} \,\delta_{k\ell}.
\end{equation}

One can then use a linear combination of the holomorphic differentials
$\omega_{\gamma_k}$ to correct the meromorphic differentials
$\nu_{(a)-(b)}$ in such a way that the resulting meromorphic
differentials have purely imaginary $b_k$--periods.
Let $X_\ell(a,b)$ be coefficients such that the differentials of
the third kind
\begin{equation}\label{IIIkind}
 \omega_{(a)-(b)}:=\nu_{(a)-(b)}-\sum_\ell
X_\ell(a,b)\omega_{g_\ell}
\end{equation}
have purely imaginary $b_k$--periods. The
coefficients $X_\ell(a,b)$ satisfy the system of equations
\begin{equation}\label{X-equation}
\sum_{\ell=1}^g X_\ell(a,b)\,Re\,\tau_{kl}= {\rm
Re}\,\int_{b_k}\nu_{(a)-(b)}= \sum_{h\in S(g_k)} \log |\langle
a,b,z^+(h),z^-(h)\rangle |.
\end{equation}

Thus, one obtains (\cite{Man-hyp}, \cf also \cite{Wer}) from
\eqref{Green} and \eqref{IIIkind} that the
Arakelov Green function for $X(\C)$ with Schottky uniformization
can be computed as
\begin{equation}\label{Greenformula}
\begin{array}{c}
g((a)-(b),(c)-(d))  =  \sum_{h\in\Gamma} \log |\langle
a,b,hc,hd\rangle | \\[2mm]
 -  \sum_{\ell=1}^g X_\ell(a,b)\, \sum_{h\in S(g_\ell)} \log
|\langle z^+(h),z^-(h), c,d\rangle |. \end{array}
\end{equation}

Notice that this result seems to indicate that there is a choice of
Schottky uniformization involved as additional data for Arakelov
geometry at arithmetic infinity. However, as we remarked previously,
at least in the case of archimedean primes that are real embeddings
(as is the case of arithmetic infinity for $\Q$), the Schottky
uniformization is determined by the real structure, by splitting
$X(\C)$ along the real locus $X(\R)$, when the latter is nonempty.

\medskip

\subsection{Green function and geodesics}

By combining the basic formula \eqref{cr-geod} with the formula
\eqref{Greenformula} for the Green function on a Riemann surface with
Schottky uniformization one can replace each term appearing in
\eqref{Greenformula} with a corresponding term which computes the
oriented geodesic length of a certain arc of geodesic in $\mX_\Gamma$:
\begin{equation}\label{Green-geod}\begin{array}{c}
  g((a)-(b),(c)-(d))=  - \sum_{h\in\Gamma} {\rm ordist}( a*\{
hc,hd \}, b* \{ hc,hd \}) \\[2mm]
+  \sum_{\ell=1}^g X_\ell(a,b)\, \sum_{h\in S(g_\ell)} {\rm
ordist}( z^+(h) *\{ c,d \}, z^-(h) * \{ c,d \}). \end{array}
\end{equation}
The coefficients
$X_\ell(a,b)$ can also be expressed in terms of geodesics, using
the equation \eqref{X-equation}.

\section{Intermezzo: Quantum gravity and black holes}

Anti de Sitter space $\AdS_{d+1}$ is a highly symmetric
space--time, which satisfies Einstein's equations with constant
curvature $R <0$. Physically it describes empty space with a
negative cosmological constant. In order to avoid time--like
closed geodesics, it is customary to pass to the universal cover
$\tilde\AdS_{d+1}$, whose boundary at infinity of the is a
compactification of $d$--dimensional Minkowski space--time. When
passing to Euclidean signature $\tilde\AdS_{d+1}$ becomes the real
hyperbolic space $\H^{d+1}$.

The $3+1$ dimensional Anti de Sitter space is a well known example
of space--time in general relativity. Topologically $\AdS_{3+1}$
is of the form $S^1\times \R^3$, while metrically it is realized
by the hyperboloid $-u^2-v^2 +x^2 +y^2 + z^2 =1$ in $\R^5$, with
the metric element $ds^2= -du^2 - dv^2 + dx^2 + dy^2 + dz^2$. The
universal cover is topologically $\R^4$. In the context of quantum
gravity, it is especially interesting to consider the case of the
$2+1$ dimensional Anti de Sitter space $\AdS_{2+1}$ and its
Euclidean counterpart, the real 3-dimensional hyperbolic space
$\H^3$.

The holography principle postulates the existence of an explicit
correspondence between gravity on a bulk space which is
asymptotically $\tilde\AdS_{d+1}$ (\eg a space obtained as a
global quotient of $\tilde\AdS_{d+1}$ by a discrete group of
isometries) and field theory on its conformal boundary at
infinity.

The relation \eqref{cr-geod}, that identifies the Green function
$$g((a)-(b),(c)-(d))$$ on $\P^1(\C)$ with the oriented length of a
geodesic arc in $\H^3$, can be thought of as an instance of this
holography principle, when we interpret one side as geodesic
propagator on the bulk space in a semiclassical approximation, and
the other side (the Green function) as the two point correlation
function of the boundary field theory. Notice that, because of the
prescribed behavior of the Green function at the singularities
given by the points of the divisor, our ``four-points''
$g((a)-(b),(c)-(d))$, when $a \to c$ and $b\to d$, gives the two
point correlator with a logarithmic divergence which is intrinsic
and does not depend on a choice of cut-off functions (unlike the
regularization often used in the physics literature).

In \cite{ManMar2} we showed that the relation \eqref{Green-geod}
between Arakelov Green functions and configurations of geodesics
in the hyperbolic handlebody, proved by Manin in \cite{Man-hyp},
provides in fact precisely the correspondence prescribed by the
holography principle, for a class of $2+1$ dimensional
space--times known as Euclidean Krasnov black holes. These include
the Ba\~nados--Teitelboim--Zanelli black holes: an important class
of space--times in the context of $(2+1)$-dimensional quantum
gravity. 

\subsection{Ba\~{n}ados--Teitelboim--Zanelli black hole}

We consider the case of a hyperbolic handlebody of genus one (a
solid torus), with conformal boundary at infinity given by an
elliptic curve.

Recall that we can describe elliptic curves via the Jacobi
uniformization. (We have already encountered it in the previous lectures,
in the context of non-commutative elliptic curves.) Let $X_q(\C)=
\C^* / q^{\Z}$ be such description of an elliptic curve, where $q$
is a hyperbolic element of $\PSL(2,\C)$ with fixed points $\{ 0,
\infty \}$ on the sphere at infinity $\P^1(\C)$ of $\H^3$, that
is, $q\in \C^*$ with $|q|<1$. The action of $q$ on $\P^1(\C)$ extends
to an action on $\bar\H^3 = \H^3\cup \P^1(\C)$ by
$$ q: (z,y) \mapsto (qz,|q|y). $$
It is easy to see that the quotient by this action
\begin{equation}\label{solid-torus}
\mX_q = \H^3 / q^\Z
\end{equation}
is topologically a solid torus, compactified at infinity by the
conformal boundary $X_q(\C)$.

\begin{center}
\begin{figure}
\epsfig{file=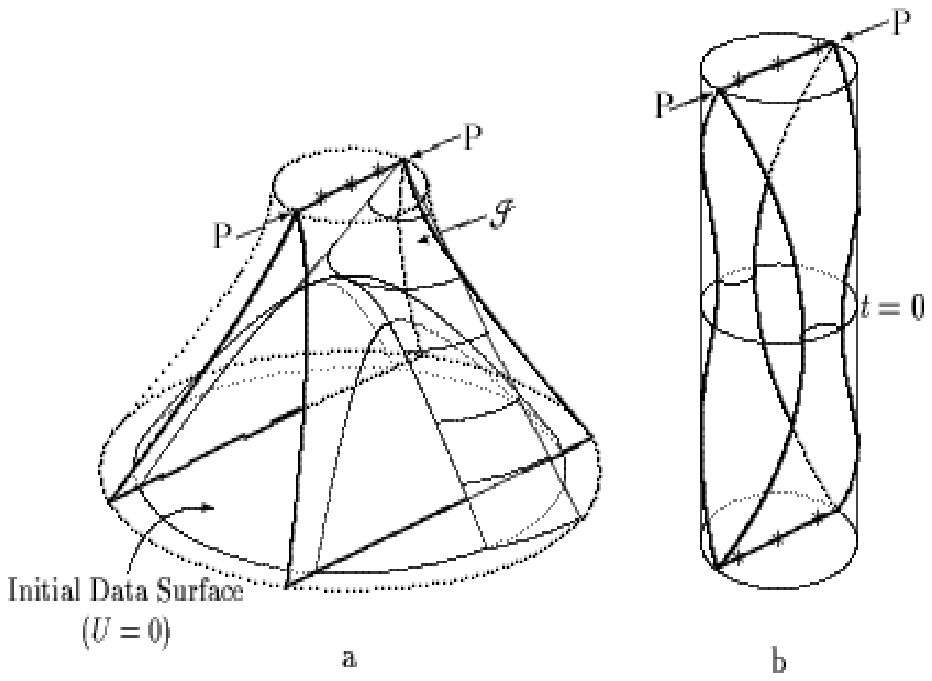}\caption{Minkowskian BTZ black
hole\label{Fig-minkBTZ}}
\end{figure}
\end{center}

The space $\mX_q$ is well known in the physics literature as the
Euclidean Ba\~{n}ados--Teitelboim--Zanelli black hole, where the
parameter $q\in \C^*$ is written in the form
$$ q= \exp\left( \frac{2\pi ( i|r_-| - r_+ ) }{\ell} \right),$$
with
$$
r_{\pm}^2 = \frac{1}{2} \left( M\ell \pm \sqrt{ M^2\ell^2+ J }
\right).
$$
Here $M$ and $J$ are the mass and angular momentum of the rotating
black hole, and $-1/\ell^2$ is the cosmological constant. The
corresponding black hole in Minkowskian signature is illustrated in
Figure \ref{Fig-minkBTZ}\footnote{Figure \ref{Fig-minkBTZ} is taken
from \cite{ABBHP}.}.

\begin{center}
\begin{figure}
\epsfig{file=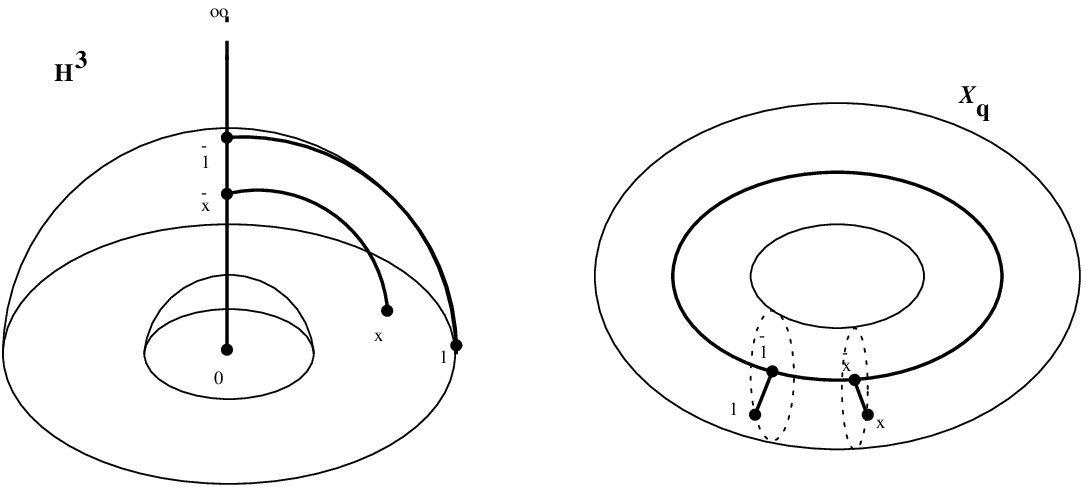}\caption{Geodesics in Euclidean BTZ black
holes\label{fig3a}}
\end{figure}
\end{center}

\begin{center}
\begin{figure}
\epsfig{file=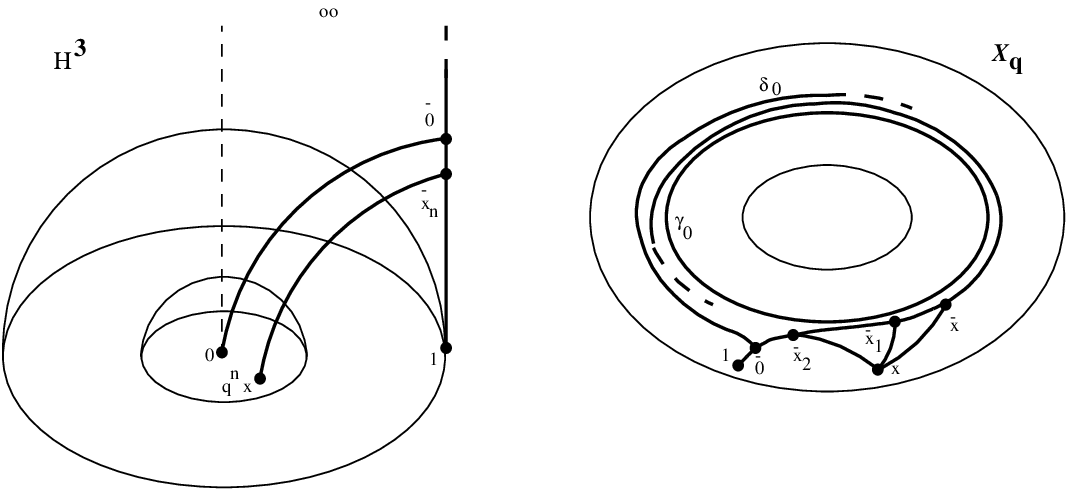}\caption{Geodesics in Euclidean BTZ black
holes\label{fig3b}}
\end{figure}
\end{center}

In the case of the elliptic curve $X_q(\C)=\C^* /q^\Z$, a formula
of Alvarez-Gaume, Moore, and Vafa gives the operator product
expansion of the path integral for bosonic field theory as
$$ g(z,1)= \log \left( |q|^{B_2(\log |z|/\log |q|)/2} |1-z| \,
\prod_{n=1}^\infty | 1- q^n z | \, \, |1-q^n z^{-1} | \right). $$
This is in fact the Arakelov Green function on $X_q(\C)$. In terms
of geodesics in the Euclidean BTZ black hole this becomes (\cf
\cite{Man-hyp})
$$
g(z,1) =-\frac{1}{2} \ell(\gamma_0)\, B_2\left(
\frac{\ell_{\gamma_0}(\bar{z},\bar{1})}{\ell(\gamma_0)}\right)
 +\sum_{n\ge 0} \ell_{\gamma_1}(\bar{0},\bar{z}_n)
+\sum_{n\ge 1} \ell_{\gamma_1}(\bar{0},\tilde{z}_n).
$$
Here we are using the notation $\bar x =x * \{ 0, \infty \}$;
$\bar z_n =q^n z * \{ 1,\infty \}$, $\tilde z_n = q^n z^{-1} * \{
1, \infty \}$ as in \cite{Man-hyp}, as illustrated in Figure
\ref{fig3a} and \ref{fig3b}. These terms describe gravitational
properties of the Euclidean BTZ black hole. For instance,
$\ell(\gamma_0)$ measures the black hole entropy. The whole
expression is a combination of geodesic propagators.

\subsection{Krasnov black holes}

The problem of computing the bosonic field propagator on an
algebraic curve $X_{\C}$ can be solved by providing differentials
of the third kind with purely imaginary periods
$$ \omega_{(a)-(b)}:=\nu_{(a)-(b)}-\sum_l X_l(a,b)\omega_{g_l}, $$
hence it can be related directly to the problem of computing the
Arakelov Green function.

Differentials as above then determine all the higher correlation
functions
$$ G(z_1,\ldots,z_m;w_1,\ldots,w_\ell)=\sum_{j=1}^m \sum_{i=1}^\ell
q_i \langle \phi(z_i,\bar z_i) \phi(w_j,\bar w_j) \rangle
q_j^\prime,
$$
for $q_i$ a system of charges at positions $z_i$ interacting with
$q_j^\prime$ charges at positions $w_j$, from the basic two point
correlator $G_\mu(a-b,z)$ given by the Green function expressed in
terms of the differentials $\omega_{(a)-(b)}(z)$. When we use a
Schottky uniformization, we obtain the differentials
$\omega_{(a)-(b)}$ as in \eqref{IIIkind}.

The bulk space corresponding to the conformal boundary $X(\C)$ is
given by the hyperbolic handlebody $\mX_\Gamma$. As in the case of
the BTZ black hole, it is possible to interpret these real
hyperbolic 3-manifolds as {\em analytic continuations} to
Euclidean signature of Minkowskian black holes that are global
quotients of $\AdS_{2+1}$. This is not just the effect of the
usual rotation from Minkowskian to Euclidean signature, but a more
refined form of ``analytic continuation'' which is adapted to the
action of the Schottky group, and which was introduced by Kirill
Krasnov \cite{Kras} \cite{Kras2} in order to deal with this class of
space-times.

The formula \eqref{Green-geod} then gives the explicit
bulk/boundary correspondence of the holography principle for this
class of space-times: each term in the Bosonic field propagator
for $X_{\C}$ is expressed in terms of geodesics in the Euclidean
Krasnov black hole $\mX_\Gamma =\H^3/\Gamma$.

\section{Dual graph and noncommutative geometry}

Manin's result on the Arakelov Green function and hyperbolic geometry
suggests a geometric model for the dual graph of the mysterious fiber
at arithmetic infinity.

In fact, the result discussed above on the Green function has an
analog, due to Drinfel'd and Manin \cite{DriMan} in the case of a
finite prime with a totally split fiber. This
is the case where the $p$-adic completion admits a Schottky
uniformization by a $p$-adic Schottky group.

\subsection{Schottky--Mumford curves}

If $K$ is a given finite extension of $\Q_p$, we denote by
$\O\subset K$ its ring of integers, by $\m\subset \O$ the maximal
ideal, and by $k$ the residue field $k = \O/\m$. This is a finite
field of cardinality $q = \rm{card}(\O/\m)$.

It is well known that a curve $X$ over a finite extension $K$ of
$\Q_p$, which is $k$-split degenerate, for $k$ the residue field,
admits a $p$-adic uniformization by a $p$-adic Schottky group
$\Gamma$ acting on the Bruhat-Tits tree $\Delta_K$.

The Bruhat-Tits tree is obtained by considering free $\sO$-modules
of rank $2$, with the equivalence relation that $M_1 \sim M_2$ iff
$\exists\lambda\in K^\ast$, such that $M_1 = \lambda M_2$. The
vertices of the tree consist of these equivalence classes. A
distance is introduced on this set by prescribing that if $
\{M_1\}, \{M_2\} \in\Delta^0_K$, with $M_1\supset M_2 $, then
$$ M_1/M_2 \simeq \sO/\m^l \oplus \sO/\m^k,\quad l,k\in\N $$
and we set the distance to be $d(\{M_1\},\{M_2\}) = \vert
l-k\vert$. We form the tree $\Delta_K$, by adding an edge between
any two pairs of vertices with $d(\{M_1\},\{M_2\}) =1$. This gives
a connected, locally finite tree with $q+1$ edges departing from
each vertex. The group $\PGL(2,K)$ acts transitively (on the left)
by isometries.

The Bruhat-Tits tree $\Delta_K$ is the analog of the 3-dimensional
real hyperbolic space $\H^3$ at the infinite primes. The set of
ends of $\Delta_K$ is identified with $\P^1(K)$, just as we have
$\P^1(\C)=\partial \H^3$ in the case at infinity.

\begin{center}
\begin{figure}
\epsfig{file=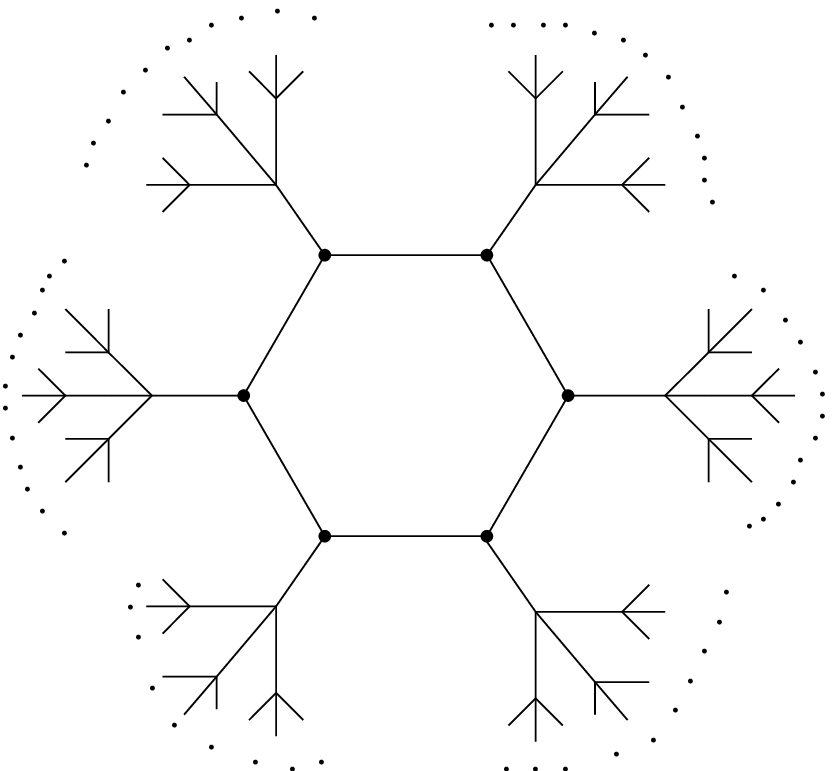}\caption{Mumford curve of genus one:
Jacobi--Tate uniformized elliptic curve
\label{Fig-JacTate}}
\end{figure}
\end{center}

A p-adic Schottky group $\Gamma$ is a discrete subgroup of
$\PGL(2,K)$ which consists of hyperbolic elements (the eigenvalues
of in $K$ have different valuation), and which is isomorphic to a
free group in $g$ generators. We still denote by
$\Lambda_\Gamma\subset \P^1(K)$ the limit set, that is, the
closure of the set of points in $\P^1(K)$ that are fixed points of
some $\gamma\in\Gamma\setminus\{1\}$. As in the case at infinity,
we have $\rm{card}(\Lambda_\Gamma) < \infty$ if and only if
$\Gamma = (\gamma)^\Z$, for some $\gamma\in\Gamma$ (the genus one
case). We denote by $\Omega_\Gamma=\Omega_\Gamma(K)$ the
complement $\Omega_\Gamma = \P^1(K)\setminus\Lambda_\Gamma$, that
is, the domain of discontinuity of $\Gamma$.

In the case of genus $g\geq 2$ the quotient
$$ X_\Gamma := \Omega_\Gamma /\Gamma $$
is a Schottky--Mumford curve, with $p$-adic Schottly
uniformization. The case of genus one gives the Jacobi-Tate
uniformization of the elliptic curve (Figure \ref{Fig-JacTate}).

A path in $\Delta_K$, infinite in both directions and with no
back-tracking, is called an {\it axis} of $\Delta_K$. Any two
points $z_1,z_2\in \P^1(K)$ uniquely define their connecting axis
with endpoints at $z_1$ and $z_2$ in $\partial\Delta_K$. The
unique axis of $\Delta_K$ whose ends are the fixed points of an
hyperbolic element $\gamma$ is called the axis of $\gamma$. The
element $\gamma$ acts on its axis as a translation. We denote by
$\Delta'_\Gamma \subset \Delta_K$ the smallest subtree containing
the axes of all elements of $\Gamma$.

\begin{center}
\begin{figure}
\epsfig{file=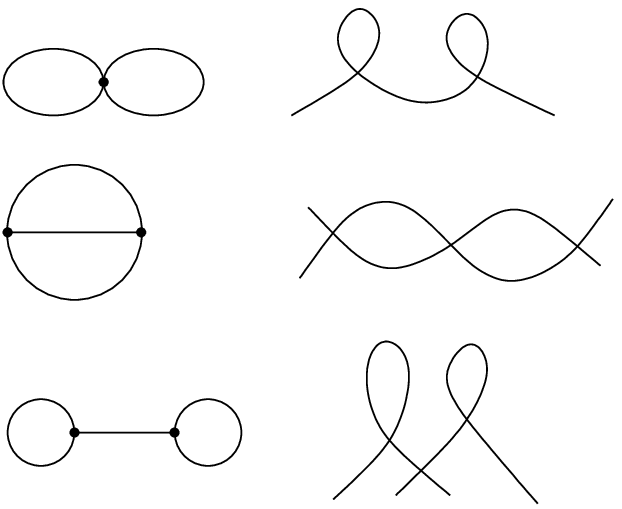}\caption{Genus two:
special fibers and dual graphs
\label{Fig-Mumford2graphs}}
\end{figure}
\end{center}

\begin{center}
\begin{figure}
\epsfig{file=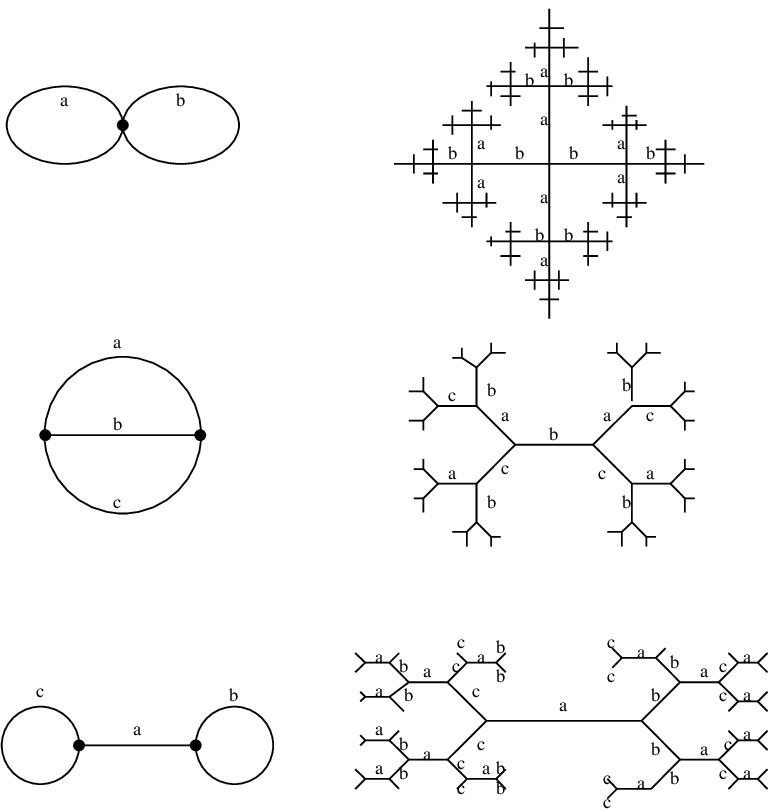}\caption{Genus two: trees and dual graphs
\label{Fig-Mumford2trees}}
\end{figure}
\end{center}

This subtree is $\Gamma$-invariant, with set of ends
$\Lambda_\Gamma$. The quotient $\Delta'_\Gamma /\Gamma$ is a {\it
finite graph}, which is the dual graph of the closed fiber of the
{\it minimal smooth model} over $\O$ ($k$-split degenerate
semi-stable curve) of $X_\Gamma$. Figure \ref{Fig-Mumford2graphs}
shows the special fiber and the corresponding dual graph for all the
possible cases of maximal degenerations special for genus two. Figure
\ref{Fig-Mumford2trees} shows the corresponding trees
$\Delta'_\Gamma$.

For each $n\geq 0$, we can also consider a subgraph
$\Delta_{K,n}$ of the Bruhat-Tits tree $\Delta_K$ defined by
setting
$$ \Delta_{K,n}^0 := \{ v\in \Delta_K^0 : \,  d(v,
\Delta_\Gamma ')\leq n \}, $$ where $d$ is the distance on
$\Delta_K$ and $d(v,\Delta_\Gamma ') := \inf \{ d(v,\tilde v): \,
\tilde v \in (\Delta_\Gamma ')^0 \}$, and
$$ \Delta_{K,n}^1 :=\{ w\in \Delta_K^1 : \, s(w), r(w)\in
\Delta_{K,n}^0 \}. $$ In particular, we have $\Delta_{K,0} =
\Delta_\Gamma '$.

For all $n\in \N$, the graph $\Delta_{K,n}$ is invariant under the
action of the Schottky group $\Gamma$ on $\Delta$, and the finite
graph $\Delta_{K,n}/\Gamma$ gives the dual graph of the reduction
$X_K \otimes \O/ \m^{n+1}$.

For a more detailed account of Schottky--Mumford curves see \cite{Mu}
and \cite{Ma}.

\subsection{Model of the dual graph}

The dictionary between the case of Mumford curves and the case at
arithmetic infinity is then summarized as follows:

\begin{tabular}{|c|c|}\hline
Bruhat--Tits $\Delta_K$ & Hyperbolic 3-space $\H^3$ \\
\hline
$\P^1(K)=\partial \Delta_K$ & $\P^1(\C)=\partial \H^3$  \\ 
\hline
Schottky $\Gamma\subset \PGL(2,K)$ & Schottky $\Gamma\subset \PSL(2,\C)$
\\ 
\hline
paths in $\Delta_K$ & geodesics in $\H^3$
\\ 
\hline
Mumford curve $\Omega_\Gamma/\Gamma$  &  Riemann surface
$\Omega_\Gamma/\Gamma$ \\ 
\hline
tree $\Delta_\Gamma'$ & convex core in $\H^3$ \\ 
\hline
graph $\Delta_K/\Gamma$ & handlebody $\mX_\Gamma$
 \\ 
\hline
dual graph $\Delta_\Gamma'/\Gamma$ & bounded geodesics in $\mX_\Gamma$
 \\ 
\hline
\end{tabular}

\bigskip

Since we can identify bounded geodesics in $\mX_\Gamma$ with infinite
geodesics in $\H^3$ with endpoints on $\Lambda_\Gamma\subset
\P^1(\C)$, modulo the action of $\Gamma$, these are parameterized by
the complement of the diagonal in $\Lambda_\Gamma \times_\Gamma
\Lambda_\Gamma$. This quotient is identified with the quotient of the
totally disconnected space $\cS$ of \eqref{Sshift} by the action of
the invertible shift $T$ of \eqref{Tshift}. Thus we obtain the
following model for the dual graph of the fiber at infinity:

\begin{itemize}
\item The solenoid ${\mathcal S}_T$ of \eqref{mapTorus} is a geometric
model of the dual graph of the fiber at infinity of an arithmetic
surface.
\item The noncommutative space $\cO_A$, representing the algebra of
coordinates on the quotient $\Lambda_\Gamma /\Gamma$, corresponds to
the set of ``vertices of the dual graph''(set of components of the
fiber at infinity), while the noncommutative space \eqref{CSTprod},
corresponding to the quotient $\Lambda_\Gamma \times_\Gamma
\Lambda_\Gamma$ gives the set of ``edges of the dual graph''.
\end{itemize}

Moreover, using noncommutative geometry, one can also give a notion of
``reduction mod $\infty$'' analogous to the reduction maps mod $p^m$
defined by the graphs $\Delta_{K,n}$ in the case of Mumford curves.

In fact, the reduction map corresponds to the paths connecting
ends of the graph $\Delta_{K,n}/\Gamma$ to corresponding vertices of the
graph $\Delta_\Gamma'/\Gamma = \Delta_{K,0} /\Gamma$. The analog at
arithmetic infinity consists then of
geodesics in $\mX_\Gamma$ which are the image of geodesics in $\H^3$
starting at some point $x_0\in \H^3\cup \Omega_\Gamma$ and having the
other end at a point of $\Lambda_\Gamma$. These are parameterized by
the set
$$ \Lambda_\Gamma \times_\Gamma (\H^3\cup \Omega_\Gamma). $$

Thus, in terms of noncommutative geometry, the reduction mod $\infty$
corresponds to a compactification of the homotopy quotient
$$ \Lambda_\Gamma \times_\Gamma \H^3 = \Lambda_\Gamma \times_\Gamma
\underline{E}\Gamma $$
where
$$ \underline{B}\Gamma = \H^3/\Gamma =\mX_\Gamma. $$

Thus, we can view $\Lambda_\Gamma /\Gamma$ as the quotient of a
foliation on the homotopy quotient with contractible leaves
$\H^3$. The reduction mod $\infty$ is then given by the $\mu$-map
 $$ \mu: K^{*+1} (\Lambda_\Gamma \times_\Gamma \H^3) \to
K_*({\rm C}(\Lambda_\Gamma)\rtimes \Gamma). $$

This shows that the noncommutative space $(\cO_A,\cH,\sD)$ is closely
related ot the geometry of the fiber at arithmetic infinity of an
algebraic variety. One can ask then what arithmetic information is
captured by the Dirac operator $\sD$ of this spectral triples. We'll
see in the next section that (as proved in \cite{CM}), the Dirac
operator gives another important arithmetic invariant, namely the
local $L$-factor at the archimedean prime.

\section{Arithmetic varieties and $L$--factors}

An important invariant of arithmetic varieties is the
$L$-function. This is written as a product of contributions from
the finite primes and the archimedean primes,
\begin{equation}\label{Lfunction}
\prod_{\wp\in \overline{\Sp\, O_\K}} L_\wp (H^m (X),s).
\end{equation}
We do not plan to give a detailed account on the subject, but we refer
the interested reader to \cite{Serre}. Here we only try to convey some
basic ideas.

The reason why one needs to consider also the contribution of the
archimedean primes can be seen in the case of the ``affine line''
$\Sp(\Z)$, where one has the Riemann zeta function, which is
written as the Euler product
\begin{equation}\label{zetaEuler}
 \zeta(s)=\prod_p (1-p^{-s})^{-1}.
\end{equation}
However, to have a nice functional equation, one needs to consider
the product
\begin{equation}\label{zetaGamma}
\zeta(s)\, \Gamma(s/2) \pi^{-s/2},
\end{equation}
which includes a contribution of the archimedean prime, expressed
in terms of the Gamma function
\begin{equation}\label{Gammafunct}
\Gamma(s)= \int_0^{\infty}e^{-t}t^{s-1}\, dt.
\end{equation}

An analogy with ordinary geometry suggests to think of the
functional equation as a sort of ``Poincar\'e duality'', which
holds for a compact manifold, hence the need to ``compactify'' 
arithmetic varieties by adding the archimedean primes and the
corresponding archimedean fibers. 

When one looks at an arithmetic variety over a finite prime $\wp \in
\Sp(O_\K)$, the fact that the reduction lives over a residue field
of positive characteristic implies that there is a special
operator, the {\em geometric Frobenius} $Fr_\wp^*$, acting on a
suitable cohomology theory (\'etale cohomology), induced by the 
Frobenius automorphism $\phi_\wp$ of $\Gal(\bar{\mathbb
F}_p/{\mathbb F}_p)$. 

The local $L$-factors of \eqref{Lfunction} at finite primes 
encodes the action of the geometric Frobenius in the form \cite{Serre}
\begin{equation}\label{L-factor}
  L_\wp (H^m (X),s) = 
\det\left( 1-Fr_\wp^* N(\wp)^{-s} | H^m(\bar X,
\Q_\ell)^{I_\wp} \right)^{-1}.
\end{equation}

Here we are considering the action of the geometric Frobenius
$Fr_\wp^*$ on the {\em inertia invariants} $H^m(\bar X,
\Q_\ell)^{I_\wp}$ of the \'etale cohomology. An introduction to
\'etale cohomology and a precise definition of these arithmetic
structures is  beyond the purpose of these notes. In fact, our primary
concern will only be the contribution of the archimedean primes to
\eqref{Lfunction}, where the construction will be based on ordinary de
Rham cohomology. Thus, we only give a quick and somewhat heuristic
explanation of \eqref{L-factor}. We refer to \cite{Serre} \cite{Weil}
for a detailed and rigorous account and for the precise hypothesis
under which the following holds.

For $X$ a smooth projective algebraic variety (in any dimension)
defined over $\Q$, the notation $\bar X$ denotes 
$$ \bar X := X \otimes \Sp(\bar{\Q}), $$
where $\bar\Q$ is an algebraic closure.
For $\ell$ a prime, the cohomology $H^*(\bar X,\Q_\ell)$ is a
finite dim $\Q_\ell$-vector space satisfying
\begin{equation}\label{comparison}
 H^i(X(\C),\C) \simeq H^i(\bar X,\Q_\ell) \otimes \C. 
\end{equation}
The absolute Galois group $\Gal(\bar \Q/\Q)$ acts on $H^*(\bar
X,\Q_\ell)$.

Similarly, we can consider $H^*(\bar X,\Q_\ell)$ for $X$ defined over
a number field $\K$. For $\wp\in \Sp(O_\K)$ and $\ell$ a prime such
that $(\ell,q)=1$, where $q$ is the cardinality of the residue field
at $\wp$, the inertia invariants
\begin{equation}\label{inertiaIv}
H^*(\bar X,\Q_\ell)^{I_\wp} \subset H^*(\bar X,\Q_\ell)
\end{equation}
are the part of the $\ell$-adic cohomology where the inertia group
at $\wp$ acts trivially. The latter is defined as
\begin{equation}\label{inertiaG}
I_\wp = \Ker(D_\wp \to \Gal(\bar{\mathbb F}_p/\mathbb
F_p)),
\end{equation}
with $D_\wp = \{\sigma \in \Gal(\bar\Q/\Q)\,\,|\, \sigma(\wp) =
\wp \}$. 
The Frobenius automorphism of $\Gal(\bar{\mathbb
F}_p/\mathbb F_p)$ lifts to $\phi_\wp \in D_\wp/I_\wp$ and induces the
geometric Frobenius $F_\wp^* := (\phi_\wp^{-1})^*$ acting on $H^*(\bar
X,\Q_\ell)^{I_\wp}$. We use the notation $N$ in \eqref{L-factor} for
the norm map.  
Thus, we can write the local $L$-factor \eqref{L-factor} equivalently as
\begin{equation}\label{loc-factor}
L_\wp (H^m (X),s) = \prod_{\lambda\in \Sp(Fr_\wp^*)} (1-\lambda
q^{-s})^{-\dim H^m(X_\Gamma)_\lambda^{I_\wp}},
\end{equation}
where $H^m(X_\Gamma)_\lambda^{I_\wp}$ is the eigenspace of
the Frobenius with eigenvalue $\lambda$.

For our purposes, what is most important to retain from the discussion
above is that the local $L$-factor \eqref{L-factor} depends upon the
data 
\begin{equation}\label{pairHFr}
\left( H^*(\bar X,\Q_\ell)^{I_\wp}, Fr_\wp^* \right)
\end{equation}
of a vector space, which has a cohomological interpretation, together
with a linear operator.

\subsection{Archimedean $L$-factors}

Since \'etale cohomology satisfies the compatibility
\eqref{comparison}, if we again resort to the general philosophy,
according to which we can work with the smooth complex manifold
$X(\C)$ and gain information on the ``closed fiber'' at airhtmetic
infinity, we are lead to expect that the contribution of the
archimedean primes to the $L$-function may be expressed in terms of
the cohomology $H^*(X(\C),\C)$, or equivalently in terms of de Rham
cohomology. 

In fact, Serre showed (\cite{Serre}) that the expected contribution of
the archimedean primes depends upon the Hodge structure
\begin{equation}\label{Hodge-str}
 H^m(X(\C))= \oplus_{p+q=m} H^{p,q}(X(\C))
\end{equation}
and is again expressed in terms of Gamma functions, as in the
case of \eqref{zetaGamma}. Namely, one has a product of Gamma
functions according to the Hodge numbers $h^{p,q}$,
\begin{equation}\label{archL-factor}
L(H^*,s)= \left\{ \begin{array}{l} \prod_{p,q}
\Gamma_\C(s-\text{min}(p,q))^{h^{p,q}} \\[2mm]
\prod_{p<q}\Gamma_\C(s-p)^{h^{p,q}}\prod_p
\Gamma_\R(s-p)^{h^{p+}}\Gamma_\R(s-p+1)^{h^{p-}}
\end{array}\right.
\end{equation}
where the two cases correspond, respectively, to the complex and the
real embeddings. Here $h^{p,\pm}$ is the dimension of the
$\pm(-1)^p$-eigenspace of the involution on $H^{p,p}$ induced by the
real structure and
\begin{equation}\label{GammaCR}
\Gamma_\C(s) := (2\pi)^{-s}\Gamma(s), \ \ \ 
\Gamma_\R(s) :=2^{-1/2}\pi^{-s/2}\Gamma(s/2).
\end{equation}

One of the general ideas in arithmetic geometry is that one should
always seek a unified picture of what happens at the finite and at the
infinite primes. In particular, there should be a suitable
reformulation of the local factors \eqref{L-factor} and
\eqref{archL-factor} where both formulae can be expressed in the same
way. 

Seeking a unified description of local $L$-factors at finite and
infinite primes, Deninger in \cite{Den1}, \cite{Den}, \cite{Den2}
expressed both \eqref{L-factor} and \eqref{archL-factor} as infinite
determinants. 

Recall that the Ray-Singer determinant of an operator $T$ with pure
point spectrum with finite multiplicities $\{ m_\lambda \}_{\lambda
\in \Sp(T)}$ is defined as
\begin{equation}\label{RSdet}
\det_\infty(s-T) := \exp\left(-\frac{d}{dz}\zeta_T(s,z)_{|z=0}\right)
\end{equation}
where the zeta function of $T$ is defined as
\begin{equation}\label{zetaT}
\zeta_T(s,z) =
\sum_{\lambda\in \Sp(T)}m_\lambda(s-\lambda)^{-z}.
\end{equation}
Suitable conditions for the convergence of these expressions in the
case of the local factors are described in \cite{Man4}.

Deninger showed that \eqref{loc-factor} can be written equivalently in
the form  
\begin{equation}\label{detL-factor}
 L_\wp (H^m(X),s)^{-1}= det_\infty(s-\Theta_q), 
\end{equation}
for an operator with spectrum
\begin{equation}\label{Thetaq-Sp}
\Sp(s-\Theta_q)=\left\{ \frac{2\pi i}{\log q} \left( \frac{\log
q}{2\pi i}(s- \alpha_\lambda) + n \right) : n\in \Z,\,\, \lambda \in
\Sp(Fr_v^*) \right\},
\end{equation}
with multiplicities $d_\lambda$ and with
$q^{\alpha_\lambda}=\lambda$. 

Moreover, the local factor \eqref{archL-factor}
at infinity can be written similarly in the form
\begin{equation}\label{detL-arch} 
 L(H^q(X),s) =
\det_{\infty}\left(\frac{1}{2\pi}(s-\Phi)|_{{\mathcal H}^m}
\right)^{-1} 
\end{equation}
where ${\mathcal H}^m$ is an infinite dimensional vector space and 
$\Phi$ is a linear operator with spectrum $\Sp(\Phi)=\Z$ and finite
multiplicites. This operator is regarded as a ``logarithm of
Frobenius'' at arithmetic infinity.

Given Deninger's formulae \eqref{detL-factor} and \eqref{detL-arch}, 
it is natural to ask for a cohomological interpretation of the
data 
\begin{equation}\label{pairHPhi}
\left( {\mathcal H}^m, \Phi \right)
\end{equation}
somewhat analogous to \eqref{pairHFr}.

\subsection{Arithmetic surfaces: $L$-factor and Dirac operator} 

Let us now return to the special case of arithmetic surfaces, in the
case of genus $g\geq 2$.

At an archimedean prime, we consider the Riemann surface
$X_\alpha(\C)$ with a Schottky uniformization
$$ X(\C)=\Omega_\Gamma /\Gamma. $$
In the case of a real embedding we can assume that the choice of
Schottky uniformization is the one that corresponds to the real
structure, obtained by cutting $X(\C)$ along the real locus $X(\R)$. 

Consider the spectral triple $(\cO_A, \cH, \sD)$ associated to the
Schottky group $\Gamma$ acting on its limit set
$\Lambda_\Gamma$. Since the spectral triple is not finitely summable,
we do not have zeta functions of the spectral triple in the form
$\Tr(a |\sD|^z)$. However, we can consider the restriction of $\sD$ to
a suitable subspace of $\cH$, where it becomes of trace class.

In particular, consider the zeta function
\begin{equation}\label{zetaDOA}
\zeta_{\pi_{{\mathcal V}},
{\mathcal D}} (s,z):= \sum_{\lambda \in \Sp({\mathcal D})} {\rm Tr}
\left( \pi_{{\mathcal V}} \Pi(\lambda, {\mathcal D})
\right) (s-\lambda)^{-z}
\end{equation}
where $\pi_{{\mathcal V}}$ is the orthogonal projection of $\cH$ onto
the subspace ${\mathcal V}$ of $0\oplus {\mathcal
L}\subset {\mathcal H}$ defined by
\begin{equation}\label{PinV}
\Pi_n\, \pi_{{\mathcal V}}\, \Pi_n= \sum_{i=1}^g S_i^n
{S_i^*}^n,
\end{equation}
on ${\mathcal
L}$, where $\Pi_n$ are the spectral projections of the Dirac operator,
$\Pi_n = \Pi(-n, {\mathcal D})$, for $n\geq 0$.  

In the case of an arithmetic surface the interesting local factor is
the one for the first cohomology, $L(H^1(X),s)$. This can be computed
(\cite{CM}, \cite{CM2}) from the zeta function \eqref{zetaDOA} of the
spectral triple $(\cO_A, \cH, \sD)$.

\begin{thm}\label{regdetSP3}
The local $L$-factor \eqref{archL-factor} is given by 
\begin{equation}\label{LfactSP3}
L(H^1(X),s) =\exp\left( - \frac{d}{dz} \zeta_{\pi_{{\mathcal V}},
\frac{{\mathcal 
D}}{2\pi}} \left(\frac{s}{2\pi},z\right)|_{z=0}
\right)^{-1}. 
\end{equation}
In case of a real embedding, the same holds, with the projection 
$\pi_{{\mathcal V}, \bar F_\infty=id}$ onto $+1$ eigenspace of the
involution $\bar F_\infty$ induced by the real structure on ${\mathcal
V}$.
\end{thm}

This result shows in particular that, for the special case of
arithmetic surfaces with $X(\C)$ of genus $g\geq 2$, the pair
\begin{equation}\label{pairVD}
\left({\mathcal V}\subset {\mathcal L}, {\mathcal
D}|_{\mathcal V}\right)
\end{equation}
is a possible geometric construction of the pair $({\mathcal
H}^1,\Phi)$ of \eqref{pairHPhi}. In particular, the Dirac operator of
the spectral triple has an arithmetic meaning, in as it recovers the
``logarithm of Frobenius''
\begin{equation}\label{DPhi}
\sD|_\cV = \Phi.
\end{equation}

If we look more closely at the subspace $\cV$ of the Hilbert space
$\cH$, we see that it has a simple geometric interpretation in terms
of the geodesics in the handlebody $\mX_\Gamma$. Remeber that the
filtered subspace $\cP$ of $\cL=L^2(\Lambda_\Gamma,d\mu)$ describes
1-cochains on the mapping torus $\cS_T$, with
$\cP/\delta\cP=H^1(\cS_T)$. 

The mapping torus $\cS_T$ is a copy of the tangle of bounded geodesics
inside $\mX_\Gamma$. Among these geodesics there are $g$
``fundamental'' closed geodesics that correspond to the generators of
$\Gamma$ (they correspond to geodesics in $\H^3$ connecting the fixed
points $\{ z^\pm(\gamma_i) \}$, for $i=1,\ldots,g$). Topologically,
these are the $g$ core handles of the handlebody $\mX_\Gamma$ and they
generate the homology $H_1(\mX_\Gamma)$. 

Consider the cohomology $H^1(\cS_T)$. Suppose that we wish to find
elements of $H^1(\cS_T)$ supported on these fundamental closed
geodesics. This is impossible, because forms on $\cS_T$ are defined by
functions in $C(\cS,\Z)$, that are supported on some clopen set
covering the totally disconnected set $\cS$, which contains no
isolated points. However, it is possible to choose a sequence of
1-cochains on $\cS_T$ whose supports are smaller and smaller clopen
sets containing the infinite word in $\cS$ corresponding to one of the
$g$ fundamental geodesics.
The finite dimensional subspace $\cV_n=\cV\cap \cP_n$ 
$$ \cV_n \subset \cP_n \ \ \  \dim \cV_n = 2g $$
gives representatives of exactly such cohomology classes in
$H^1(\cS_T)$. (We get $2g$ instead of $g$ because we take into account
the two possible choices of orientation.)

Thus, this gives us, in the case of arithmetic surfaces, a cohomological
interpretation of the space $\cH^1=\cV$ in the pair \eqref{pairHPhi}.
Moreover, the Schottky uniformization also provides us with a way of
expressing this cohomology $\cH^1=\cV$ in terms of the de Rham cohomology
$H^1(X(\C))$. In fact, we have already seen in the calculation of the
Green function that, under the hypothesis that
$\dim_H(\Lambda_\Gamma)<1$, to each generator $\gamma_i$ of the
Schottky group we can associate a holomorphic differential on the
Riemann surface $X(\C)$ by \eqref{1k}. The map
$$ \gamma_i \mapsto \omega_{\gamma_i}=\sum_{h\in C(|\gamma_i)}\,
d\,\log \langle 
hz^+(\gamma_i), hz^-(\gamma_i),z,z_0\rangle
$$
thus gives an identification
\begin{equation}\label{Varchcoh}
\cV \simeq \oplus_{n\in \Z_{\geq 0}} H^1(X(\C)).
\end{equation}
We'll see in the next section that, in fact, the right hand side of
\eqref{Varchcoh} is a particular case of a more general construction
that works for arithmetic varieties in any dimensions and that gives a
cohomological interpretation of \eqref{pairHPhi}.

\section{Archimedean cohomology}

Consani gave in \cite{KC}, for general arithmetic varieties (in any
dimension), a cohomological 
interpretation of the pair $(\cH^m,\Phi)$ in Deninger's
calculation of the archimedean $L$-factors as regularized
determinants. 

\begin{figure}
\begin{center}
\epsfig{file=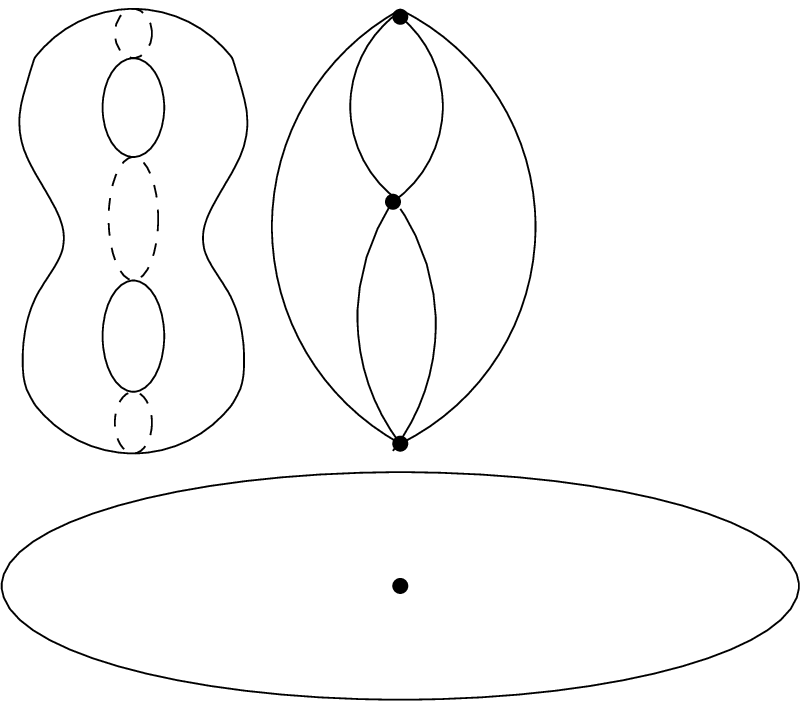}
\end{center}
\caption{Vanishing cycles
\label{Fig-vanishing}}
\end{figure}

Her construction was motivated by the analogy between geometry at
arithmetic infinity and the classical geometry of a degeneration over
a disk. She introduced a double complex of differential forms
with an endomorphism $N$ representing the ``logarithm of the
monodromy'' around the special fiber at arithmetic infinity, which is
modelled on (a resolution of) the complex of {\em nearby cycles} in the
geometric case. The definition of the complex of nearby cycles and of
its resolution, on which the following construction is modelled is
rather technical. What is easier to visualize geometrically is the
related complex of the {\em vanishing cycles} of a geometric
degeneration (Figure \ref{Fig-vanishing}).

We describe here the construction of \cite{KC} using the notation of
\cite{CMrev}. We construct the cohomology theory underlying
\eqref{pairHPhi} in several steps.

In the following we let $X=X(\C)$ be a complex compact K\"ahler
manifold. 

\medskip

{\bf Step 1:} We begin by considering a doubly infinite graded complex
\begin{equation}\label{Ccomplex}
C^\cdot = \Omega^\cdot(X) \otimes \C [U, U^{-1}]\otimes \C [\hbar,
\hbar^{-1}],
\end{equation}
where $\Omega^\cdot(X)$ is the de Rham complex of differential forms
on $X$, while $U$ and $\hbar$ are formal variables, with $U$ of degree
two and $\hbar$ of degree zero. 

On this complex we consider differentials 
\begin{equation}\label{dC}
 d'_C:=\hbar \, d, \ \ \ d_C ''= \sqrt{-1}(\bar\partial -\partial),
\end{equation}
with total differential $\delta_C=d'_C + d_C ''$.

We also have an inner product
\begin{equation}\label{inner}
 \langle \alpha \otimes U^r\otimes
\hbar^k, \eta \otimes U^s \otimes \hbar^t \rangle := \langle
\alpha, \eta \rangle \, \, \delta_{r,s} \delta_{k,t}
\end{equation}
where $\langle \alpha, \eta \rangle$ is the usual Hodge inner product
of forms,
\begin{equation}\label{innerHodge}
\langle \alpha, \eta \rangle= \int_X \alpha \wedge *
\,C(\bar\eta),
\end{equation}
with $C(\eta)=(\sqrt{-1}\, )^{p-q}$, for $\eta\in \Omega^{p,q}(X)$. 

\medskip

{\bf Step 2:} We use the Hodge filtration 
\begin{equation}\label{filtrHodge}
F^p \Omega^m(X) :=
\oplus_{p'+q=m, \, p'\geq p} \, \Omega^{p',q}(X)
\end{equation}
to define linear subspaces of \eqref{Ccomplex} of the form 
\begin{equation}\label{fCmr}
\fC^{m,2r} = \displaystyle{\bigoplus_{\substack{p+q=m \\k\geq
\max\{ 0, 2r+m\}}}} \,\,\, F^{m+r-k} \, \Omega^m(X) \, \otimes U^r
\otimes \hbar^k
\end{equation}
and the $\Z$-graded vector space
\begin{equation}\label{fCcomplex}
\fC^\cdot=\oplus_{\cdot=m+2r} \fC^{m,2r}.
\end{equation}

\medskip

{\bf Step 3:} We pass to a real vector space by considering 
\begin{equation}\label{Tcomplex}
\T^\cdot=(\fC^\cdot)^{c=id},
\end{equation}
where $c$ denotes complex conjugation.

In terms of the intersection of the Hodge filtrations 
\begin{equation}\label{gammaFiltr}
\gamma^\cdot =F^\cdot \cap \bar F^\cdot
\end{equation}
we can write \eqref{Tcomplex} as 
$$ \T^\cdot= \oplus_{\cdot=m+2r} \T^{m,2r} $$
where
\begin{equation}\label{Tmr}
\T^{m,2r} = \displaystyle{\bigoplus_{\substack{p+q=m
\\k\geq \max\{ 0, 2r+m\}}}} \,\,\, \gamma^{m+r-k} \,
\Omega^m(X)\,  \otimes U^r \otimes \hbar^k.
\end{equation}

The $\Z$-graded complex vector space $\fC^\cdot$ is a subcomplex of
$C^\cdot$ with respect to the differential $d'_C$, and for $P^\perp$
the orthogonal projection onto $\fC^\cdot$ in the inner product
\eqref{inner}, we obtain a second differential $d''=P^\perp
d_C''$. Similarly, $d'=d_C'$ and $d''=P^\perp d_C''$ define
differentials on the $\Z$-graded real vector space $\T^\cdot$, since
the inner product \eqref{inner} is real on real forms and induces an
inner product on $\T^\cdot$. We write $\delta=d'+d''$ for the total
differential.

\begin{figure}
\begin{center}
\epsfig{file=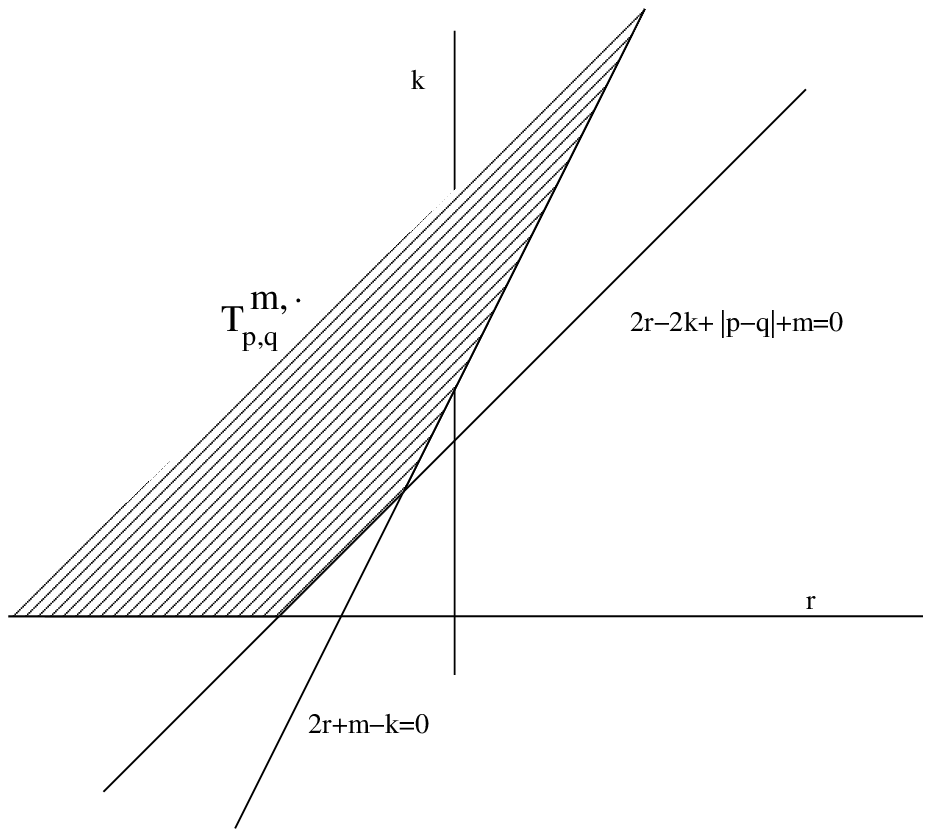}
\end{center}
\caption{Cutoffs defining the complex at arithmetic infinity
\label{Fig-Figkappa}}
\end{figure}

We can describe the real vector spaces $\T^\cdot$ in terms of certain
cutoffs on the indices of the complex $C^\cdot$. Namely, for
\begin{equation}\label{Lambdapq}
\Lambda_{p,q}=\{ (r,k)\in \Z^2 : \,\, k\geq
\kappa(p,q,r) \}
\end{equation}
with
\begin{equation}\label{kappapq}
\kappa(p,q,r):= \max\left\{0,2r+m,\frac{|p-q|+2r+m}{2}\right\}
\end{equation}
(\cf Figure \ref{Fig-Figkappa}),
we identify $\T^\cdot$ as a real vector space with the span
\begin{equation}\label{Tcutoff}
\T^\cdot = \R \langle \alpha\otimes U^r \otimes \hbar^k \rangle,
\end{equation}
where $(r,k)\in \Lambda_{p,q}$ for $\alpha=\xi+\bar\xi$, with $\xi \in
\Omega^{p,q}(X)$. 

\subsection{Operators}

The complex $(\T^\cdot,\delta)$ has some interesting structures given
by the action of certain linear operators. 

We have the operators $N$ and $\Phi$ that correspond to the
``logarithm of the monodromy'' and the ``logarithm of
Frobenius''. These are of the form
\begin{equation}\label{NPhi}
N=U \hbar  \ \ \ \ \  \Phi = -U \frac{\partial}{\partial U}
\end{equation}
and they satisfy $[N,d']=[N, d'']=0$ and $[\Phi,d']=[\Phi, d'']=0$,
hence they induce operators in cohomology.

Moreover, there is another important operator, which corresponds to
the Lefschetz operator on forms,
\begin{equation}\label{Lefschetz}
\hL : \eta \otimes U^r \otimes \hbar^k \mapsto  \eta \wedge \omega
\otimes U^{r-1} \otimes \hbar^k, 
\end{equation}
where $\omega$ is the K\"ahler form on the manifold $X$.
This satisfies $[\hL,d']=[\hL, d'']=0$, so it also descends on the
cohomology. 

The pairs of operators $N$ and $\Phi$ or $\hL$ and $\Phi$ satisfy the
interesting commutation relations
$$ [\Phi,N]=-N, \ \ \ \ [\Phi,\hL]=\hL $$
that can be seen as an action of the ring of differential operators
$$ \C[P,Q]/(PQ - QP =  Q). $$

\subsection{$\SL(2,\R)$ Representations}

Another important piece of the structure of $(\T^\cdot,\delta)$ are
two involutions
\begin{equation}\label{dual1}
 S: \alpha \otimes U^r \otimes \hbar^{2r+m+\ell} \mapsto \alpha \otimes
U^{-(r+m)} \otimes \hbar^\ell 
\end{equation}
\begin{equation}\label{dual2}
 \tilde S : \alpha \otimes U^r \otimes \hbar^k \mapsto C(*\alpha)
\otimes U^{r-(n-m)} \otimes \hbar^k .
\end{equation}

These maps, together with the nilpotent operators $N$ and $\hL$ define
two representations of $\SL(2,\R)$.
In terms of 
\begin{equation}\label{genSL2R}
\begin{array}{l}
  \nu(s) := \left(\begin{array}{cc} s & 0 \\ 0 & s^{-1}
\end{array} \right)\,\, s \in \R^* \\[3mm]
  u(t):= \left(\begin{array}{cc} 1 & t \\ 0 & 1 \end{array}
\right)\,\, t\in \R \\[3mm]  w:= \left(\begin{array}{cc} 0 & 1 \\
-1 & 0 \end{array} \right),
\end{array}
\end{equation}
the representations $\sigma^L$ and $\sigma^R$ are given by
\begin{equation}\label{sigmaLR}
\begin{array}{lr}
\sigma^L(\nu(s)) = s^{-n+m} & \sigma^R(\nu(s))
= s^{2r+m} \\[3mm] \sigma^L(u(t))= \exp(t\,\hL) & \sigma^R(u(t))=
\exp(t\,N)
\\[3mm] \sigma^L(w)= (\sqrt{-1})^n\,C\, \tilde S & \sigma^R(w)= C\,
S. \end{array} 
\end{equation}

Of these representations, $\sigma^L$ extends to an action by bounded
operators on the Hilbert completion of $\T^\cdot$ in the inner product
\eqref{inner}, while the action of the subgroup $\nu(s)$, $s\in \R^*$
of $\SL(2,\R)$ via the representation $\sigma^R$ on this Hilbert space
is by unbounded densely defined operators. 

\subsection{Renormalization group and monodromy}

This very general structure, which exists for varieties in any
dimension, also has interesting connections to noncommutative
geometry. For instance, we can see that in fact the map $N$ does play
the role of the ``logarithm of the monodromy'' using an analog in our
context of the theory of renormalization \`a la Connes--Kreimer
\cite{CK2}. 

In the classical case of a geometric degeneration on a disk, the
monodromy around the special fiber is defined as the map
\begin{equation}\label{Tmonodromy}
T = \exp(-2\pi\sqrt{-1}\,\,\text{Res}_0(\nabla))
\end{equation}
where 
\begin{equation}\label{Nmonodromy}
N = \text{Res}_0(\nabla) 
\end{equation}
is the residue at zero of the connection, acting as an endomorphism of
the cohomology.

In our setting, 
we consider loops $\phi_\mu$ with values in the group $G=
\Aut(\T^\cdot_\C,\delta)$, depending on a ``mass parameter'' $\mu\in
\C^*$. Here $\T^\cdot_\C$ is the complexification of the real vector
space $\T^\cdot$.
The Birkhoff decomposition of a loop  $\phi_\mu$ consists in
the multiplicative decomposition
\begin{equation}\label{Birkhoff}
\phi_\mu(z)=   \phi_\mu^-(z)^{-1}\, \phi^+_\mu(z),
\end{equation}
for $z\in \partial\Delta \subset \P^1(\C)$, where $\Delta$ is a small
disk centered at zero. Of the two terms in the right hand side of
\eqref{Birkhoff}, $\phi^+_\mu$ extends to a holomorphic function on
$\Delta$ and $\phi_\mu^-$ to a holomorphic functions on
$\P^1(\C)\smallsetminus \Delta$ with values in $G$. We normalize
\eqref{Birkhoff} by requiring that $\phi^-_\mu(\infty)=1$. 

By analogy with the Connes--Kreimer theory of renormalization we
require the following properties of \eqref{Birkhoff}:
\begin{itemize}
\item The time evolution
\begin{equation}\label{timeevol}
\theta_t:a\mapsto e^{-t\Phi} a\, e^{t\Phi}
\end{equation}
acts by scaling
\begin{equation}\label{timescaling}
\phi_{\lambda \mu}(\epsilon) = \theta_{t\epsilon}\,
\phi_\mu(\epsilon), 
\end{equation}
for $\lambda=e^t \in \R^*_+$ and $\epsilon\in
\partial\Delta$.
\item The term $\phi^-_\mu= \phi^-$ in the Birkhoff decomposition is
independent of the energy scale $\mu$. 
\end{itemize}

The residue of a loop $\phi_\mu$, as in \cite{CK2}, is given by
\begin{equation}\label{ResCK}
{\rm Res}\, \phi = \frac{d}{d z}\left( \phi^- (1/z)^{-1}
 \right)|_{z=0}
\end{equation}
and the beta function of renormalization is 
\begin{equation}\label{betafunct}
\beta = \Upsilon\, {\rm Res}\, \phi, \ \ \ \text{ with } \ \ 
\Upsilon =\frac{d}{dt}\theta_t |_{t=0}.
\end{equation}

In our setting, there is a natural choice of the time evolution, given
by the ``geodesic flow'' associated to the ``Dirac operator'' $\Phi$,
namely
\begin{equation}\label{timePhi}
\theta_t(a) = e^{-t\Phi} \, a \, e^{t\Phi}.
\end{equation}
This gives
\begin{equation}\label{Upsilon}
 \Upsilon(a) =\frac{d}{dt} \theta_t(a) |_{t=0} = [a,\Phi]. 
\end{equation}

There is a {\em scattering formula} (\cite{CK2}), by which one can
reconstruct $\phi^-$ from the residue. Namely, we can write
$$ \phi^-(z)^{-1}= 1+ \sum_{k\ge 1}
d_k\, z^{-k}, $$
with
\begin{equation}\label{scattering}
d_k = \int_{s_1\ge \cdots \ge s_k \ge 0} \theta_{-s_1}(\beta) \cdots
\theta_{-s_k}(\beta) \, ds_1 \cdots ds_k .
\end{equation}

The {\em renormalization group} is given by
\begin{equation}\label{renormgroup}
\rho(\lambda)= \lim_{\epsilon \to 0 }\;\phi^-(\epsilon) \,
\theta_{t\epsilon} (\phi^- (\epsilon)^{-1}),
\end{equation}
for $\lambda=e^t \in
\R^*_+$. 

Thus, we only need to specify the residue in order to have the
corresponding renormalization theory associated to
$(\T^\cdot_\C,\delta)$. By analogy to the case of the geometric
degeneration \eqref{Nmonodromy} it is natural to require that ${\rm
Res}\, \phi =N$. We then have (\cite{CMrev}):

\begin{prop}\label{renormTdelta}
A loop $\phi_\mu$ in $G=
\Aut(\T^\cdot_\C,\delta)$ with ${\rm Res}\, \phi_\mu =N$, subject to
\eqref{timescaling} and with $\phi^-$ independent of $\mu$, satisfies
$$ \phi_\mu(z) = \exp\left( \frac{\mu^z}{z}\, N \right) $$
with Birkhoff decomposition
$$ \phi_\mu(z) = \exp(-N/z)
\exp\left( \frac{\mu^z-1}{z}\, N \right). $$ 
\end{prop}

In fact, by \eqref{Upsilon} and ${\rm Res}\, \phi_\mu =N$ we have 
$\beta = [N,\Phi]=N$ and $\theta_t(N)= e^t \,\, N$, hence the
scattering formula \eqref{scattering} gives 
$$ \phi^-(z)=\exp(-N/z) $$
and the scaling property \eqref{timescaling} determines
$$ \phi_\mu(z) = \exp\left( \frac{\mu^z}{z}\, N \right). $$

The part of the Birkhoff decomposition that is regular at $z=0$
satisfies $\phi_\mu^+(0)=\mu^N$.

The renormalization group is of the form
\begin{equation}\label{renhoroc}
\rho(\lambda) = \lambda^N = \exp(t N),
\end{equation}
which, through the representation $\sigma^R$ of $\SL(2,\R)$,
corresponds to the horocycle flow on $\SL(2,\R)$
$$ \rho(\lambda) = u(t)=\left(\begin{array}{cc} 1 & t \\ 0 & 1
\end{array} \right). $$

The Birkhoff decomposition \eqref{Birkhoff} gives a trivialization of
a principal $G$-bundle over $\P^1(\C)$. We can consider the associated
vector bundle $\sE_\mu^\cdot$ with fiber $\T^\cdot_\C$. 

Moreover, we obtain a Fuchsian connection $\nabla_\mu$ on this bundle,
\begin{equation}\label{Fuchsian}
\nabla_\mu : \sE_\mu^\cdot \to
\sE_\mu^\cdot \otimes_{\sO_\Delta} \Omega^\cdot_{\Delta}(\log 0),
\end{equation}
where the notation $\Omega^\cdot_{\Delta}(\log 0)$ denotes forms with
logarithmic poles at $0$. This has the form
$$ \nabla_\mu = N\, \left( \frac{1}{z} + \frac{d}{dz} \frac{\mu^z-1}{z}
\right) \, dz, $$
with local gauge potentials
$$ -\phi^+(z)^{-1} \frac{\log
\pi(\gamma)\, dz}{z} \phi^+(z) + \phi^+(z)^{-1}\, d\phi^+(z) $$
with respect to the monodromy representation 
$$ \pi: \pi_1(\Delta^*)=\Z \to G $$
given by
\begin{equation}\label{repMonodromy}
 \pi(\gamma) =\exp(-2\pi\sqrt{-1} \, N) 
\end{equation}
for $\gamma$ the generator of $\pi_1(\Delta^*)$. This corresponds to
\eqref{Tmonodromy} in the classical geometric case. This Fuchsian
connection has residue
$$ {\rm Res}_{z=0}
\nabla_\mu =N $$
as in the geometric case \eqref{Nmonodromy}.

\subsection{Local factor and archimedean cohomology}

In \cite{KC} Consani showed that the data $(\cH^m,\Phi)$ of
\eqref{pairHPhi} can be identified with 
$$ \left( \H^\cdot(\T^\cdot,\delta)^{N=0}, \Phi \right), $$
where $\H^\cdot(\T^\cdot,\delta)$ is the hypercohomology (the
cohomology with respect to the total differential $\delta$) of the
complex $\T^\cdot$ and $\H^\cdot(\T^\cdot,\delta)^{N=0}$ is the kernel
of the map induced by $N$ on cohomology. The operator $\Phi$ is the
one induced on cohomology by that of \eqref{NPhi}. She called 
$\cH^m\simeq \H^\cdot(\T^\cdot,\delta)^{N=0}$ the {\em archimedean
cohomology}. 

This can also be viewed (\cf \cite{KC}) as a piece of the cohomology
of the cone of the monodromy $N$. This is the complex
\begin{equation}\label{ConeN}
 {\rm Cone}(N)^\cdot =
\T^\cdot \oplus \T^\cdot [+1] 
\end{equation} 
with differential
$$ D=\left(\begin{array}{cc} \delta & -N \\ 0 & \delta
\end{array}\right). $$

The complex \eqref{ConeN} inherits a positive definite inner product
from $\T^\cdot$, which 
descends on cohomology. The representation $\sigma^L$ of $\SL(2,\R)$
on $\T^\cdot$ induces a representation on ${\rm Cone}(N)^\cdot$. The
corresponding representation $d\sigma^L : {\mathfrak g} \to {\rm
End}(\T^\cdot)$ of 
the Lie algebra ${\mathfrak g}=sl(2,\R)$ extends to a representation of the
universal enveloping algebra $U({\mathfrak g} )$
on $\T^\cdot$ and on ${\rm Cone}(N)^\cdot$. This gives a
representation in the algebra of bounded operators on the Hilbert
completion of ${\rm Cone}(N)^\cdot$ in the inner product. 

\begin{thm}\label{spectral3H}
The triple
$$ (\cA,\cH,\sD)=(U({\mathfrak g}),\H^\cdot({\rm Cone}(N)), \Phi) $$
has the properties that $\sD=\sD^*$ and that $(1+\sD^2)^{-1/2}$ is a
compact operator. The commutators $[\sD,a]$ are bounded operators for
all $a\in U({\mathfrak g})$ and the triple is $1^+$-summable.
\end{thm}

Thus, $(\cA,\cH,\sD)$ has most of the properties of a spectral triple,
confirming the fact that the logarithm of Frobenius $\Phi$ should be
thought of as a Dirac operator. However, we are not dealing here with an
involutive subalgebra of a $C^*$-algebra. 

In any case, the structure is sufficient to consider zeta functions
for this ``spectral triple''. In particular, we can recover the
alternating products of the local $L$-factors at infinity from a zeta
function of the spectral triple.

\begin{thm}\label{Lspectral3H}
Consider the zeta function 
$$ \zeta_{a,\Phi}(z)= \Tr(a |\Phi|^{-z}) $$
with $a=\sigma^L(w)$. This gives 
$$
{\det_{\infty}}_{\sigma^L(w),\Phi}(s) = \prod_{m=0}^{2n}
L(H^m(X),s)^{(-1)^m}.
$$
\end{thm}

\subsection{Arithmetic surfaces: homology and cohomology}

In the particular case of arithmetic surfaces, there is an
identification (\cite{KC}, \cite{CM})
\begin{equation}\label{HCone}
\H^\cdot({\rm Cone}(N)) \simeq \cH^\cdot \oplus \check{\cH}^\cdot,
\end{equation}
where $\cH^\cdot$ is the archimedean cohomology and
$\check{\cH}^\cdot$ is its dual under the involution $S$ of
\eqref{dual1}. 

We can then extend the identification 
$$ U: \cH^1\stackrel{\simeq}{\to} \cV \subset \cL $$ 
of \eqref{Varchcoh}, by considering a subspace $\cW$ of the homology
$H_1(\cS_T)$ with $\cW\simeq \check{\cH}^\cdot$. The homology
$H_1(\cS_T)$ can also be computed as a direct limit
$$ H_1({\mathcal
S}_T,\Z)= \varinjlim_N {\mathcal K}_N, $$
where the ${\mathcal K}_N$ are free abelian of rank
$(2g-1)^N +1$ for $N$ even and $(2g-1)^N + (2g-1)$ for $N$ odd. The
$\Z$-module ${\mathcal K}_N$ is generated by the closed geodesics
represented by periodic sequences in $\cS$ of period $N+1$. These need
not be primitive closed geodesics. In terms of primitive closed
geodesics we can write equivalently 
$$ H_1({\mathcal S}_T,\Z) = \oplus_{N=0}^\infty {\mathcal R}_N $$
where the ${\mathcal R}_N$ are free abelian groups with
$$ {\rm rk}\,{\mathcal R}_N=  \frac{1}{N}\sum_{d|N} \mu(d) {\rm rk}\,
{\mathcal K}_{N/d},$$
with $\mu(d)$ the M\"obius function satisfying $\sum_{d|N}\mu(d)
=\delta_{N,1}$. 

The pairing of homology and cohomology is given by 
$$ \langle \cdot, \cdot \rangle : F_n \times
{\mathcal K}_N \to \Z   \ \ \ \ \  \langle [f], x \rangle = N
\cdot f(\bar x). $$
This determines a graded subspace $\cW\subset H_1({\mathcal S}_T,\Z)$
dual to $\cV\subset H^1(\cS_T)$. With the identification
$$ \diagram \cH^1 \rrto^{
U} \dto^{S} & & \cV \subset H^1({\mathcal
S}_T) \dto^{<,>} \\
\check{\cH}^1 \rrto^{ U} & & {\mathcal W} \subset
H_1({\mathcal S}_T),
\enddiagram $$ 
we can identify the Dirac operator of \eqref{DiracDH} with
the logarithm of Frobenius 
\begin{equation}\label{DiracFrobeniusCone}
\sD|_{\cV\oplus\cV} = \Phi_{\cH^1 \oplus \check{\cH^1}}.
\end{equation}


\begin{thebibliography}{9999}

\bibitem{Alling} N.~Alling, N.~Greenleaf, {\em Foundations of the
theory of Klein surfaces}, Lecture Notes in Mathematics Vol. 219,
Springer Verlag 1971.

\bibitem{ABBHP}  S.~Aminneborg, I.~Bengtsson, D.~Brill, S.~Holst,
P.~Peldan, {\em Black holes and wormholes in $2+1$ dimensions},
Class. Quantum Grav. 15 (1998) 627--644.

\bibitem{alr} J.~Arledge, M.~Laca, I.~Raeburn, {\em Semigroup crossed
products and Hecke algebras arising from number fields},
Doc. Math. 2 (1997) 115--138.

\bibitem{Babenko} K.~I.~Babenko, {\it On a problem of Gauss}.
Dokl. Akad. Nauk SSSR, Tom 238 (1978) No. 5, 1021--1024.

\bibitem{Barrow} J.~D.~Barrow. {\it Chaotic behaviour and the Einstein
equations.} In: Classical General Relativity, eds. W.~Bonnor et
al., Cambridge Univ. Press, Cambridge, 1984, 25--41.

\bibitem{BaumConnes} P.~Baum, A.~Connes, {\em Geometric K-theory
for Lie groups and  foliations}. Preprint IHES 1982;
l'Enseignement  Mathematique, t. 46, 2000, 1-35.

\bibitem{Boca} F.P.~Boca, {\em Projections in rotation algebras
and theta functions}, Commun. Math. Phys. 202 (1999) 325--357.

\bibitem{BC} J.B.~Bost, A.~Connes, {\em Hecke algebras, Type III
factors and phase transitions with spontaneous symmetry breaking
in number theory}, Selecta Math. (New Series) Vol.1 (1995) N.3,
411--457.

\bibitem{Bo} R.~Bowen, {\em Hausdorff dimension of quasi--circles}, 
Publ.Math. IHES 50 (1979) 11--25. 
 
\bibitem{BoHa} M.~Boyle, D.~Handelman, {\em Orbit equivalence, flow 
equivalence, and ordered cohomology}, Israel J. Math. 95 (1996) 
169--210. 

\bibitem{BR} O.~Bratteli, D.W.~Robinson, {\em Operator algebras and
quantum statistical mechanics I,II}, Springer Verlag, 1981.

 
\bibitem{ChaPal} P.~Chakraborty, A.~Pal, {\em Equivariant spectral
triples on the quantum ${\rm SU}(2)$ group},  $K$-Theory  28  (2003),
no. 2, 107--126.

\bibitem{ChMayer} C.H.Chang and D.Mayer, {\em Thermodynamic formalism and
Selberg's zeta function for modular groups}, Regular and chaotic
dynamics 15 (2000) N.3 281--312.

\bibitem{Cohen} P.B.~Cohen, {\em A $C^*$-dynamical system with
Dedekind zeta partition function and spontaneous symmetry breaking},
Journal de Th\'eorie des Nombres de Bordeaux 11 (1999) 15--30.


\bibitem{ConnesCR} A.~Connes, {\it  $C^*$ alg\`ebres et
g\'eom\'etrie differentielle}.  C.R. Acad. Sci. Paris, Ser.~A-B ,
290 (1980) 599--604.

\bibitem{Co-Thom} A.~Connes, {\em  An analogue of the Thom isomorphism
for crossed products of a $C^*$-algebra by an action of $\R$}, 
Adv. in Math. 39 (1981), no. 1, 31--55.

\bibitem{Co} A.~Connes, {\em Non--commutative differential
geometry}, Publ.Math. IHES N.62 (1985) 257--360.

\bibitem{Co-transv} A.~Connes, {\it Cyclic cohomology and the
transverse fundamental class of a foliation}. In: Geometric
methods in operator algebras (Kyoto, 1983). Pitman Res. Notes in
Math., 123, Longman, Harlow 1986, 52--144.

\bibitem{Co-fredh}  A.~Connes, {\it Compact metric spaces, Fredholm
modules, and hyperfiniteness}, Ergod. Th. Dynam. Sys. (1989) 9,
207--220.

\bibitem{Co94}  A.~Connes, {\em Noncommutative geometry}, Academic
Press, 1994.

\bibitem{Connes} A.~Connes, {\em Geometry from the spectral point of
view}. Lett. Math. Phys. 34 (1995), no. 3, 203--238.


\bibitem{Connes-Zeta} A.~Connes,  {\it Trace formula in Noncommutative
Geometry and the zeros of the Riemann zeta function}. Selecta
Mathematica. New Ser. 5 (1999) 29--106.


\bibitem{Co2} A.~Connes, {\em A short survey of noncommutative
geometry}, J. Math. Phys. 41 (2000), no. 6, 3832--3866.


\bibitem{CoQgr} A.~Connes, {\em Cyclic cohomology,
 Quantum group Symmetries and the Local
Index Formula for $SU_q(2)$},  J. Inst. Math. Jussieu  3  (2004),
no. 1, 17--68.

\bibitem{CDS} A.~Connes, M.~Douglas, A.~Schwarz,
{\it Noncommutative geometry and  Matrix theory: compactification
on tori}. J. High Energy Phys. (1998) no. 2, Paper 3, 35 pp.
(electronic)

\bibitem{CDV1}
 A.~Connes, M.~Dubois-Violette, {\em Noncommutative
finite-dimensional manifolds. I. spherical manifolds and related
examples}, Comm. Math. Phys.  Vol.230  (2002) N.3, 539--579.

\bibitem{CDV2} A.~Connes, M.~Dubois-Violette,
{\em Moduli space and
structure of noncommutative 3-spheres}, Lett. Math. Phys., Vol.66 (2003)
N.1-2, 91--121.

\bibitem{CK2} A.~Connes, D.~Kreimer, {\em Renormalization in quantum
field theory and the Riemann-Hilbert problem. II. The
$\beta$-function, diffeomorphisms and the renormalization group}.
Comm. Math. Phys.  216  (2001),  no. 1, 215--241.

\bibitem{CoMar} A.~Connes, M.~Marcolli, {\em Quantum Statistical
Mechanics of $\Q$-lattices, (From Physics to Number Theory via
Noncommutative Geometry, Part I)}, preprint arXiv:math.NT/0404128.

\bibitem{cmln} A.~Connes, M.~Marcolli, {\em $\Q$-lattices: quantum
statistical mechanics and Galois theory}, to appear in Journal of
Geometry and Physics.

\bibitem{CMR} A.~Connes, M.~Marcolli, N.~Ramachandran {\em KMS states
and complex multiplication}, preprint 2004.

\bibitem{ConnesMosc} A.~Connes, H.~Moscovici, {\it The local index
formula in noncommutative geometry}. Geom. Funct. Anal. 5 (1995),
no. 2, 174--243.

\bibitem{CoMo1} A.~Connes, H.~Moscovici, {\em Modular
Hecke algebras and their Hopf symmetry}, Moscow Math. Journal,
Vol.4 (2004) N.1, 67--109.

\bibitem{CoMo2} A.~Connes, H.~Moscovici,
{\em Rankin-Cohen Brackets and the Hopf Algebra of Transverse
Geometry}, Moscow Math. Journal, Vol.4 (2004) N.1, 111--130.

\bibitem{KC} C.~Consani, {\em Double complexes and Euler $L$--factors},
Compositio Math. 111 (1998) 323--358.

\bibitem{CM} C.~Consani, M.~Marcolli, {\em Noncommutative
geometry, dynamics and $\infty$-adic Arakelov geometry}, to appear in
Selecta Mathematica

\bibitem{CM1} C.~Consani, M.~Marcolli, {\em Triplets spectreaux in
geometrie d'Arakelov}, C.R.Acad.Sci. Paris, Ser. I 335 (2002) 779--784.

\bibitem{CM2} C.~Consani, M.~Marcolli, {\em New perspectives in
Arakelov geometry}, to appear in proceedings of CNTA7 meeting Montreal 2002.

\bibitem{CM3} C.~Consani, M.~Marcolli, {\em Spectral triples from
Mumford curves}, International Mathematics Research Notices 36 (2003)
1945--1972.

\bibitem{CMrev} C.~Consani, M.~Marcolli, {\em Archimedean cohomology
revisited}, preprint 2004.

\bibitem{CoSilv} G.~Cornell and J.H.~Silverman (eds.) {\em Arithmetic
geometry}, Springer-Verlag, New York, 1986

\bibitem{CuKri} J.~Cuntz, W.~Krieger, {\it A class of $C^*$--algebras
and topological Markov chains}, Invent. Math. 56 (1980) 251--268.

\bibitem{Den1} C.~Deninger, {\it On the $\Gamma$--factors attached to
motives}, Invent. Math. 104 (1991) 245--261.

\bibitem{Den} C.~Deninger, {\it Local $L$-factors of motives and
regularized determinants}. Invent. Math. 107 (1992), no. 1,
135--150.

\bibitem{Den2} C.~Deninger, {\it Motivic $L$-functions and regularized
determinants}, in ``Motives'', Proceedings of Symposia in Pure
Mathematics, Vol. 55 (1994) Part I, 707--743.

\bibitem{DriMan} V.~Drinfeld, Yu.I.~Manin, {\em Periods of $p$-adic
Schottky groups}.  J. Reine Angew. Math.  262/263  (1973), 239--247.

\bibitem{Haag} R.~Haag, {\em Local Quantum Physics},
Springer, Berlin 1992.

\bibitem{HHW} R. Haag, N. M. Hugenholtz, M. Winnink
 {\em On the equilibrium states in quantum statistical mechanics},
Comm. Math. Phys. 5 (1967), 215, 236.

\bibitem{HaLe} D.~Harari, E.~Leichtnam, {\em Extension du
ph\'enom\`ene de brisure spontan\'ee de sym\'etrie de Bost--Connes au
cas des corps globaux quelconques}, Selecta Math. (New Series) Vol.3
(1997) 205--243.

\bibitem{Ha} G.~Harder, {\em General aspects in the theory of modular
symbols}, Seminar on number theory, Paris 1981--82 (Paris,
1981/1982), 73--88, Progr. Math., 38, Birkh\"auser Boston, Boston,
MA, 1983.

\bibitem{Hen2} D.~Hensley, {\em Continued fraction Cantor sets, Hausdorff
dimension, and functional analysis}, J. Number Theory 40 (1992)
336--358.

\bibitem{KLKSS} I.~M.~Khalatnikov, E.~M.~Lifshitz, K.~M.~Khanin, L.~N.~Schur,
Ya.~G.~Sinai. {\it On the stochasticity in relativistic
cosmology.} J.~Stat.~Phys., 38:1/2 (1985), 97--114.

\bibitem{Kras} K.~Krasnov, {\em Holography and Riemann Surfaces},
Adv. Theor. Math. Phys. 4 (2000) 929--979

\bibitem{Kras2} K.~Krasnov, {\em Analytic Continuation for
Asymptotically AdS 3D Gravity}, gr-qc/0111049

\bibitem{Laca} M.~Laca, {\em Semigroups of $\sp *$-endomorphisms,
Dirichlet series, and phase transitions}. J. Funct. Anal. 152 (1998),
no. 2, 330--378.

\bibitem{laca-end} M.~Laca, {\em From endomorphisms to automorphisms
and back: dilations and full corners}, J. London Math. Soc. (2) 61
(2000) 893--904.

\bibitem{Lang} S.~Lang, {\em Elliptic Functions}, (Second Edition),
  Graduate Texts in Mathematics, Vol.112, Springer-Verlag 1987.

\bibitem{Landi} G.~Landi, {\em An introduction to noncommutative
spaces and their geometries}, Lecture Notes in Physics, Vol. m-51,
Springer Verlag 1997.

\bibitem{Lang-Arak} S.~Lang, {\em
Introduction to Arakelov Theory},
Springer-Verlag, New York, 1988.

\bibitem{Lou} R.~Loudon, {\em The Quantum Theory of Light}, third
edition, Oxford University Press, 2000.

\bibitem{mandel} L.~Mandel, E.~Wolf, {\em Optical Coherence and
Quantum Optics}, Cambridge University Press, 1995.

\bibitem{Man3} Yu.I.~Manin, {\em Von Zahlen und Figuren},
arXiv:math.AG/0201005.

\bibitem{Man5} Yu.I.~Manin, {\em Real Multiplication and
noncommutative geometry}, arXiv:math.AG/0202109.

\bibitem{Man-theta} Yu.I.~Manin, {\em Theta functions, quantum tori,
and Heisenberg groups}, Lett. in Math. Phys. 56 (2001) 295--320.

\bibitem{Man4} Yu.I.~Manin, {\em Lectures of zeta functions and motives
(according to Deninger and Kurokawa)}. Columbia University Number
Theory Seminar (New York, 1992). Ast\'erisque No. 228 (1995), 4,
121--163.

\bibitem{Man-hyp} Yu.I.~Manin, {\em Three--dimensional hyperbolic geometry
as $\infty$--adic Arakelov geometry}, Invent. Math. 104 (1991)
223--244.

\bibitem{Ma} Yu.I.~Manin, {\em $p$-adic automorphic 
functions}. Journ. of Soviet Math., 5 (1976) 279-333. 

\bibitem{Man-sym} Yu.I.~Manin, {\em Parabolic points and zeta
functions of modular curves}, Math. USSR Izvestija, vol. 6 N. 1
(1972) 19--64. Selected Papers, World Scientific, 1996, 202--247.

\bibitem{ManMar2} Yu.I.~Manin, M.~Marcolli, {\em Holography principle
and arithmetic of algebraic curves}, Adv. Theor. Math. Phys. Vol.3
(2001) N.5, 617--650.

\bibitem{ManMar} Yu.I.~Manin, M.~Marcolli, {\em Continued fractions,
modular symbols, and noncommutative geometry}, Selecta Mathematica
(New Series) Vol.8 N.3 (2002) 475-520.

\bibitem{Mar-lyap} M.~Marcolli, {\em Limiting modular symbols and the
Lyapunov spectrum}, Journal of Number
Theory, Vol.98 N.2 (2003) 348-376.

\bibitem{Mar-cosm} M.~Marcolli, {\em Modular curves,
$C^*$-algebras and chaotic cosmology}, preprint.

\bibitem{MaTa} K.~Matsuzaki, M.~Taniguchi, {\em Hyperbolic manifolds
and Kleinian groups}, Oxford Univ. Press, 1998.

\bibitem{Mayer} D.H.~Mayer, {\em Relaxation properties of the
mixmaster universe},  Phys. Lett. A 121 (1987), no. 8-9,
  390--394.

\bibitem{Mayer2} D.H.~Mayer, {\em Continued fractions and related
transformations}, 1991.

\bibitem{Maz} B.~Mazur, {\em Courbes elliptiques et symboles
modulaires}, S\'eminaire Bourbaki vol.1971/72, Expos\'es 400--417,
LNM 317, Springer Verlag 1973, pp. 277--294.

\bibitem{Merel} L.~Merel, {\em Intersections sur les courbes
modulaires}, Manuscripta Math., 80 (1993) 283--289.

\bibitem{Mi} J.S.~Milne, {\em Canonical models of Shimura curves},
  manuscript, 2003 (www.jmilne.org)

\bibitem{Mum-Tata} D.~Mumford (with M.Nori and P.Norman), {\em Tata
lectures on theta, III}, Progress in Mathematics Vol. 97,
Birkh\"auser 1991.

\bibitem{Mu} D.~Mumford, {\em An analytic construction of 
degenerating curves over complete local rings}, Compositio Math. 24 
(1972) 129--174. 

\bibitem{MSW} D.~Mumford, C.~Series, D.~Wright, {\em Indra's
Pearls. The Vision of Felix Klein}. Cambridge University Press, New
York, 2002.

\bibitem{PaTu} W.~Parry, S.~Tuncel, {\it Classification problems in 
ergodic theory}, London Math. Soc. Lecture Notes Series 67, 1982. 

\bibitem{PoWei} M.~Pollicott, H.~Weiss, {\em Multifractal analysis of
Lyapunov exponent for continued fraction and Manneville-Pomeau
transformations and applications to Diophantine approximation},
Comm. Math. Phys. 207 (1999), no. 1, 145--171.

\bibitem{RS} D.B.~Ray, I.M.~Singer, {\em Analytic torsion for complex
manifolds}. Ann. of Math. (2) 98 (1973), 154--177.

\bibitem{Rief1} M.A.~Rieffel, {\em ${\rm C}^*$-algebras associated
to irrational rotations}, Pacific J. Math. 93 (1981) 415-429.

\bibitem{Rosenberg98}
A.L.~Rosenberg, {\em Noncommutative schemes}.
Compositio Math. 112 (1998) N.1, 93--125.

\bibitem{Serre} J.~P.~Serre, {\it Facteurs locaux des fonctions z\^eta
des vari\'et\'es alg\'ebriques (d\'efinitions et conjectures)}.
S\'em. Delange-Pisot-Poitou, exp.~19, 1969/70.

\bibitem{Sh} G.~Shimura, {\em Arithmetic Theory of Automorphic
    Functions}, Iwanami Shoten and Princeton 1971.

\bibitem{Stark} H.M.~Stark, {\em $L$-functions at $s=1$. IV. First
derivatives at $s=0$}, Adv. Math. 35 (1980) 197--235.

\bibitem{Steven} P.~Stevenhagen, {\em Hilbert's 12th problem, complex
multiplication and Shimura reciprocity}, Advanced Studies in Pure
Math. 30 (2001) ``Class Field Theory -- its centenary and prospect''
pp. 161--176.

\bibitem{Weil} A.~Weil, {\em Basic Number Theory}, Springer  1974.

\bibitem{Wells} R.O.~Wells, {\it Differential analysis on complex
manifolds}. Springer--Verlag, 1980.

\bibitem{Wer} A.~Werner, {\em Local heights on Mumford
curves}. Math. Ann. 306 (1996), no. 4, 819--831.

\bibitem{Zagier} D.~Zagier, {\em Modular forms and differential
operators}, in K. G. Ramanathan memorial issue.  Proc. Indian
Acad. Sci. Math. Sci.  104  (1994),  no. 1, 57--75.


\end{thebibliography}
\end{document}